%% file: paper.tex
\documentclass[final]{siamart220329}
\usepackage[utf8]{inputenc}
\usepackage[T1]{fontenc}
\usepackage{amsmath}
\usepackage{amssymb}
\usepackage{mathtools}
\usepackage{bm}
\usepackage{microtype}
\usepackage{algorithm}
\usepackage{algorithmicx}
\usepackage[noend]{algpseudocode}
\usepackage{tabularx}
\usepackage{overpic}
\usepackage{tikz}
\usepackage{todonotes}
\usepackage{enumitem}
\usepackage[binary-units=true]{siunitx}

\usepackage{matlab-prettifier}
\lstdefinestyle{matlab-custom}{
	language=Matlab,
	basicstyle=\footnotesize\ttfamily,
	keywordstyle=\bfseries\color{green!40!black},
	commentstyle=\itshape\color{purple!40!black},
	identifierstyle=\color{blue},
	stringstyle=\color{orange}
}

\makeatletter
\newcommand{\multiline}[1]{%
  \begin{tabularx}{\dimexpr\linewidth-\ALG@thistlm}[t]{@{}X@{}}
    #1
  \end{tabularx}
}
\makeatother

\graphicspath{{figures/}}
\makeatletter
\def\input@path{{figures/}}
\makeatother

\let\textv\v
\renewcommand{\vec}[1]{{\bm{#1}}}
\newcommand{\surf}{\Gamma}
\newcommand{\elem}{\mathcal{E}}

\newcommand{\lapbel}{\Delta_\surf}
\newcommand{\dt}{\Delta t}
\newcommand{\DtN}{\Sigma}
\newcommand{\ii}{\mathrm{i}}
\newcommand{\ee}{\mathrm{b}}
\renewcommand{\ss}{\mathrm{s}}

\newcommand{\disc}[1]{#1}
\newcommand{\dvec}[1]{\vec{\mathsf{#1}}}
\renewcommand{\u}{\xi}
\renewcommand{\v}{\TextOrMath\textv\eta}
\newcommand{\tang}{\surf}
\newcommand{\dxtang}{\partial_x^\tang}
\newcommand{\dytang}{\partial_y^\tang}
\newcommand{\dztang}{\partial_z^\tang}
\newcommand{\detg}{\lvert g \rvert}

\newcommand{\multmat}[1]{{\disc{M}[#1]}}
\newcommand{\definedas}{\coloneqq}

\newcommand{\pad}{\hspace{1pt}}
\newcommand{\indexset}{\mathcal{J}}
\newcommand{\tree}{\mathcal{T}}

\newcommand{\part}{{\text{P}}}
\newcommand{\ea}{\alpha}
\newcommand{\eb}{\beta}
\newcommand{\interface}{{\mathcal{I}_{\ea\eb}}}
\newcommand{\binorm}{{\vec{n}_b}}
\newcommand{\nelem}{N}
\newcommand{\nlevel}{L}
\newcommand{\cross}{\times}
\newcommand{\nudge}{\mkern 1mu}

\newsiamremark{remark}{Remark}

\title{A high-order fast direct solver for surface PDE\MakeLowercase{s}\thanks{Submitted to the editors \today.
\funding{This work was funded by the Simons Foundation.}}}

\author{Daniel Fortunato\thanks{Center for Computational Mathematics, Flatiron Institute, New York, NY 10010 (\email{dfortunato@flatironinstitute.org}).}}

\headers{A high-order fast direct solver for surface PDE\MakeLowercase{s}}{D. Fortunato}

\ifpdf
\hypersetup{
  pdftitle = {A high-order fast direct solver for surface PDEs},
  pdfauthor = {D. Fortunato}
}
\fi

\begin{document}

\maketitle

\begin{abstract}
We introduce a fast direct solver for variable-coefficient elliptic partial differential equations on surfaces based on the hierarchical Poincar\'{e}--Steklov method. The method takes as input an unstructured, high-order quadrilateral mesh of a surface and discretizes surface differential operators on each element using a high-order spectral collocation scheme. Elemental solution operators and Dirichlet-to-Neumann maps tangent to the surface are precomputed and merged in a pairwise fashion to yield a hierarchy of solution operators that may be applied in $\mathcal{O}(N \log N)$ operations for a mesh with $N$ elements. The resulting fast direct solver may be used to accelerate high-order implicit time-stepping schemes, as the precomputed operators can be reused for fast elliptic solves on surfaces. On a standard laptop, precomputation for a 12th-order surface mesh with over 1 million degrees of freedom takes 17 seconds, while subsequent solves take only 0.25 seconds. We apply the method to a range of problems on both smooth surfaces and surfaces with sharp corners and edges, including the static Laplace--Beltrami problem, the Hodge decomposition of a tangential vector field, and some time-dependent nonlinear reaction--diffusion systems.
\end{abstract}



\section{Introduction}

Surface-bound phenomena arise in a wide variety of applications, including electromagnetics~\cite{Epstein2010a, Epstein2013, Epstein2015, Colton2013}, plasma physics~\cite{Atanasiu1999, Garabedian2008, Malhotra2019a}, biological pattern formation~\cite{Jeong2017}, bulk--surface diffusion processes~\cite{Diegmiller2018, Miller2022}, biomechanics~\cite{Stoop2015, Miller2018}, and fluid dynamics~\cite{Saye2016b, Veerapaneni2011}. Modeling such phenomena often leads to surface partial differential equations (PDEs)---that is, partial differential equations whose differential operators are taken with respect to the on-surface metric. Thus, the numerical solution of surface PDEs is an important component in the simulation of surface phenomena.

High-order numerical methods have gained popularity in recent years due to their high levels of accuracy and efficiency per degree of freedom, and a panoply of numerical methods have been developed for discretizing and solving PDEs on surfaces using a variety of geometry representations, with both low-order and high-order accuracy. Broadly speaking, surface PDE solvers may be divided into three main categories: mesh-based methods (e.g., finite element~\cite{Dziuk2013} and integral-equation-based~\cite{ONeil2018} methods), meshless methods (e.g., radial basis function methods~\cite{Fuselier2013, Wendland2020}), and embedded methods (e.g., level set~\cite{Bertalmio2001} and closest point~\cite{vonGlehn2013} methods). Briefly, let us highlight each category with a particular focus on the fast and high-order accurate solvers available to each:

\vspace{1em}
\begin{itemize}[listparindent=\parindent]
\item \textbf{Mesh-based methods}.
Perhaps most popular among mesh-based methods is the surface finite element method (FEM)~\cite{Dziuk1988, Dziuk2007, Dziuk2007a, Dziuk2008, Demlow2009, Dziuk2013}, which weakly enforces the PDE on each element of a surface mesh, along with coupling conditions at element interfaces. Surface discretizations based on the high-order FEM have been used in the setting of isogeometric analysis, but typically employ at most cubic or quartic polynomials~\cite{Hughes2005, Jonsson2017}. In addition, little work has been done on developing efficient solvers for the resulting linear systems in the high-order context on surfaces, where metric quantities can lead to variable coefficients with strong anisotropies. Geometric multigrid methods have achieved mesh-independent convergence rates when solving surface PDEs discretized by FEM, but have traditionally been applied to discretizations using only linear~\cite{Landsberg2010, Bonito2012, Adler2018} or quadratic~\cite{Burger2009} basis functions. Sparse direct solvers such as UMFPACK~\cite{UMFPACK} may be applied to surface FEM discretizations for moderate problem sizes, but often suffer from dense fill-in at higher orders.

For simple surface PDEs such as the Laplace--Beltrami equation, second-kind integral equation formulations exist based on the Green's function for the three-dimensional Laplacian~\cite{ONeil2018}. These formulations can be discretized to high-order accuracy to yield well-conditioned matrices which are amenable to solution via an iterative method in only a constant number of iterations. Moreover, each iteration may be accelerated to linear time by the fast multipole method~\cite{Greengard1987}. However, the application of specialized near-singular quadrature rules can significantly hinder runtimes, even in modern implementations~\cite{Agarwal2021}. Finally, integral-equation-based approaches may not be applicable to the general variable-coefficient surface PDEs considered in this work.

For these mesh-based methods, high-order accuracy can be hampered by a surface mesh that is only low-order smooth. Recent work has shown that a flat surface mesh can be used to define a nearby high-order smooth surface mesh in a way that preserves multiscale features~\cite{Vico2020}.

\vspace{1em}
\item \textbf{Meshless methods}.
Meshless methods use scattered data (typically, a point cloud with neighbor information) to represent a surface via interpolation. The flexibility and robustness of the point cloud representation make such methods well suited to computing with noisy surface data, as might arise in imaging applications. Meshless methods based on radial basis functions  (RBFs) can yield high-order accurate solutions of surface PDEs~\cite{Fuselier2013,Shankar2020}, and result in sparse linear systems which may be inverted with a direct solver for moderate problem sizes. However, scaling can break down at high orders due to dense fill-in~\cite{Shankar2015}. Recently, preconditioned iterative methods~\cite{Lehto2017} and meshless multigrid methods~\cite{Wright2022} have been shown to yield ``mesh''-independent convergence rates at higher orders, at least up to order 7. RBF methods have also been applied to solve PDEs on evolving surfaces~\cite{Wendland2020}.

\vspace{1em}
\item \textbf{Embedded methods}.
Embedded methods recast a surface PDE as a volumetric PDE posed on all of $\mathbb{R}^3$ via extension, using a suitable mapping between points on the surface and points in the ambient space. Implicit surface representations based on the level set method have been used with Cartesian finite differences in the volume to yield second-order accurate discretizations for surface PDEs~\cite{Bertalmio2001, Greer2006, Greer2006a}, though solving in a narrow band can degrade the order of convergence~\cite{Greer2006}. In addition, standard iterative solvers can require many iterations to converge~\cite{Greer2006a}, and the effect of higher-order finite difference stencils on solver complexity has not been well studied in this context.

Closest point methods~\cite{Ruuth2008, Macdonald2010, vonGlehn2013, Macdonald2013} take a similar approach, but can handle a range surfaces with open ends, cusps, and corners through the use of a closest point extension. Higher-order finite difference stencils have been studied using narrow-banding with closest point extension, and the resulting linear systems are generally sparse. However, iterative methods can stall without sufficient preconditioning. To that end, narrow-band multigrid methods have been used in conjunction with the closest point method to yield mesh-independent iteration counts, though only for second-order stencils~\cite{Chen2015}.

Finally, it should be noted that embedded methods necessarily require physical quantities defined on the surface to be extended into the volume---a task that can be challenging to do with high-order accuracy~\cite{Askham2017}.

\end{itemize}
\vspace{1em}

In this work we propose a fast direct solver for a high-order spectral element discretization of a surface PDE. Given a collection of geometric patches defining the surface, the method directly discretizes the strong form of the surface PDE on each patch using classical Chebyshev spectral collocation~\cite{Trefethen2000, Boyd2001}. Continuity and continuity of the binormal derivative are imposed between neighboring patches, which are merged pairwise in a tree-like fashion using the hierarchical Poincar\'{e}--Steklov (HPS) scheme~\cite{Martinsson2013, Gillman2014, Gillman2015, Babb2018, Fortunato2021, Martinsson2019a}. The method can be viewed as an operator analogue of classical nested dissection~\cite{George1973}, with a factorization complexity of $\mathcal{O}(\nelem^{3/2})$ and subsequent solve complexity of $\mathcal{O}(\nelem \log \nelem)$ for a mesh with $\nelem$ elements.

Direct solvers offer a few key advantages in the high-order setting. First, the complexity constant associated with direct solvers is typically small (at least for two-dimensional problems). Second, the large condition numbers commonly associated with high-order discretizations do not affect the performance of a direct solver~\cite{Martinsson2019a}. Likewise, physical parameters such as variable coefficients and wavenumbers which could degrade the performance of an iterative method do not affect the runtime of a direct solver. Lastly, as the factorizations from a direct solver can be stored for reuse, repeated solves of the same PDE with different data are very fast. Because of this, direct solvers can be used to accelerate semi-implicit time-stepping schemes~\cite{Babb2020a, Babb2020b} or handle localized changes to geometry~\cite{Zhang2021a} while maintaining high-order accuracy in space. Though the surfaces we consider in this work are embedded in three dimensions, their surface PDEs are two-dimensional problems, enabling all the benefits of fast direct solvers in two dimensions.

The paper is structured as follows. In \cref{sec:formulation}, we define surface differential operators and spectral collocation on a single surface patch. \Cref{sec:dd} describes the Schur complement method for domain decomposition on surfaces, starting from the linear system of a single element and moving to two ``glued'' elements. In \cref{sec:hps}, we describe the hierarchical fast direct solver and its computational complexity. \Cref{sec:results} presents numerical examples demonstrating both the speed and high-order convergence of the method for a variety of static and time-dependent problems.

\section{Formulation and discretization}\label{sec:formulation}

\subsection{Surface differential operators}\label{sec:diffops}

Let $\surf$ be a smooth surface embedded in $\mathbb{R}^3$ locally parametrized by the mapping $\vec{x}(\u,\v) \colon \mathbb{R}^2 \to \mathbb{R}^3$, where $\vec{x}(\u,\v) = (x(\u,\v), y(\u,\v), z(\u,\v))$. The metric tensor $g$ on $\surf$ is given by
\[
g = \begin{bmatrix} g_{\u\u} & g_{\u\v} \\ g_{\v\u} & g_{\v\v} \end{bmatrix}
  = \begin{bmatrix}
\vec{x}_\u \cdot \vec{x}_\u & \vec{x}_\u \cdot \vec{x}_\v \\
\vec{x}_\v \cdot \vec{x}_\u & \vec{x}_\v \cdot \vec{x}_\v
\end{bmatrix},
\]
where $\vec{x}_\u = \partial \vec{x} / \partial \u$ and $\vec{x}_\v = \partial \vec{x} / \partial \v$. Let the components of the inverse of $g$ be denoted by
\[
g^{-1} = \begin{bmatrix} g^{\u\u} & g^{\u\v} \\ g^{\v\u} & g^{\v\v} \end{bmatrix}
       = \frac{1}{\detg} \begin{bmatrix} \phantom{-}g_{\v\v} & \!\!-g_{\u\v} \\ -g_{\v\u} & \!\!\phantom{-}g_{\u\u} \end{bmatrix},
\]
with $\detg = \det{g}$, and define the inverse metric quantities
\begin{equation}\label{eq:xi_eta_derivs}
\u_\vec{x} = g^{\u\u} \vec{x}_\u + g^{\v\u} \vec{x}_\v, \qquad \v_\vec{x} = g^{\u\v} \vec{x}_\u + g^{\v\v} \vec{x}_\v.
\end{equation}
The surface gradient $\nabla_\surf = (\dxtang, \dytang, \dztang)$ of a scalar-valued function $u = u(\u,\v) \definedas u(\vec{x}(\u,\v))$ defined on $\surf$ is a tangential vector field which can be written as~\cite{Frankel2012}
\begin{equation}\label{eq:cont_diff_xyz}
\nabla_\surf u = \begin{bmatrix}
\\[-1em] \dxtang u \\[0.4em] \dytang u \\[0.4em] \dztang u \vspace{0.1em}
\end{bmatrix}
= \begin{bmatrix}
\\[-1em]
\u_x \frac{\partial u}{\partial \u} + \v_x \frac{\partial u}{\partial \v} \\[0.4em]
\u_y \frac{\partial u}{\partial \u} + \v_y \frac{\partial u}{\partial \v} \\[0.4em]
\u_z \frac{\partial u}{\partial \u} + \v_z \frac{\partial u}{\partial \v} \vspace{0.1em}
\end{bmatrix}.
\end{equation}
The components $\dxtang$, $\dytang$, and $\dztang$ may be interpreted as tangential derivative operators on $\surf$ in the Cartesian directions $\vec{e}_x$, $\vec{e}_y$, and $\vec{e}_z$, respectively.
Similarly, the surface divergence of a vector field $\vec{u} = \vec{u}(\u,\v) = (u_1(\u,\v), u_2(\u,\v), \allowbreak u_3(\u,\v))$ tangent to $\surf$ can be written as
\[
\nabla_\surf \cdot \vec{u} = \dxtang u_1 + \dytang u_2 + \dztang u_3.
\]
The surface Laplacian $\lapbel$, or Laplace--Beltrami operator, is then given by
\[
\lapbel u = \nabla_\surf \cdot \nabla_\surf u = \dxtang \dxtang u + \dytang \dytang u + \dztang \dztang u.
\]
The Laplace--Beltrami operator is the generalization of the standard Laplacian to general manifolds. Just as the Poisson equation is the canonical elliptic PDE in Euclidean space, the Laplace--Beltrami equation is the canonical elliptic PDE on surfaces. Higher-order surface differential operators, such as the surface biharmonic operator, may be defined similarly. For simplicity of presentation, we restrict our attention to second-order surface differential operators in this work.

\begin{figure}[htb]
	\centering
	\hspace{1cm}
	\includegraphics[width=0.9\textwidth]{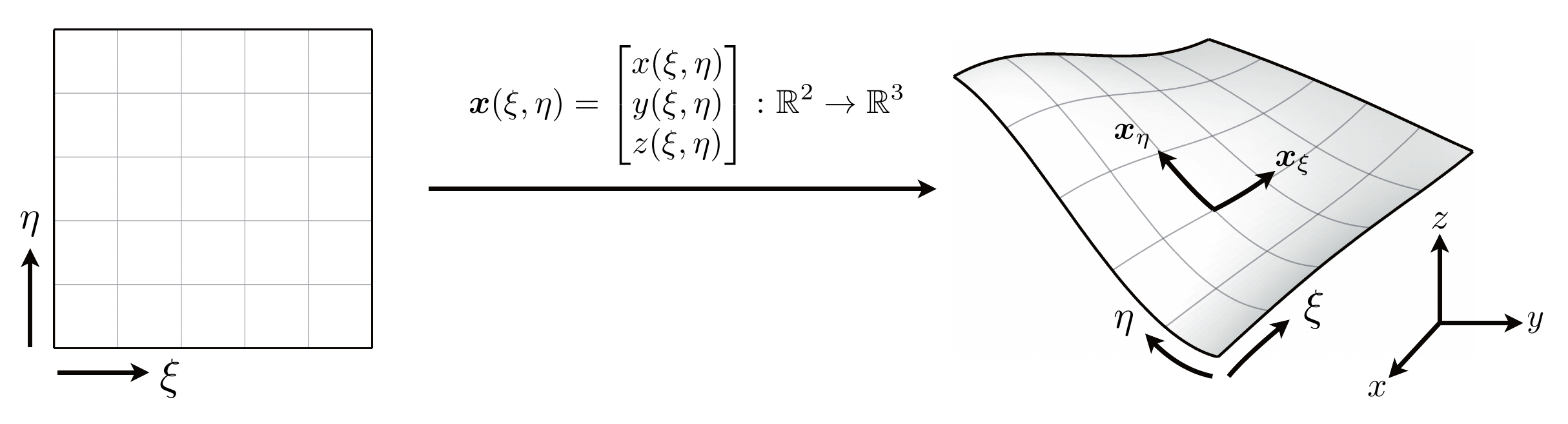}
	\vspace{-1em}
	\caption{The reference square with coordinates $(\u,\v) \in [-1,1]^2$ is mapped to a quadrilateral element via the mapping $\vec{x}(\u,\v) \colon [-1,1]^2 \to \mathbb{R}^3$. The partial derivatives $\vec{x}_\u$ and $\vec{x}_\v$ are used to form the metric tensor $g$, which encodes how lengths and angles change over the surface.}
	\label{fig:element_metric}
\end{figure}

We consider a general elliptic surface PDE on $\surf$,
\begin{equation}\label{eq:model_pde}
\mathcal{L}_\surf u(\vec{x}) = f(\vec{x}), \qquad \vec{x} \in \surf,
\end{equation}
where $f(\vec{x})$ is a smooth function on $\surf$ and $\mathcal{L}_\surf$ is a variable-coefficient, second-order, linear, elliptic partial differential operator of the form
\begin{equation}\label{eq:L_definition}
\mathcal{L}_\surf u(\vec{x}) = \sum_{i=1}^3 \sum_{j=i}^3 a_{ij}(\vec{x}) \, \partial_i^\tang \partial_j^\tang u(\vec{x}) + \sum_{i=1}^3 b_i(\vec{x}) \,\partial_i^\tang u(\vec{x}) + c(\vec{x}) u(\vec{x})
\end{equation}
with smooth coefficients $a_{11}(\vec{x})$, $a_{22}(\vec{x})$, $a_{33}(\vec{x})$, $a_{12}(\vec{x})$, $a_{23}(\vec{x})$, $a_{13}(\vec{x})$, $b_1(\vec{x})$, $b_2(\vec{x})$, $b_3(\vec{x})$, and $c(\vec{x})$ for $\vec{x} \in \surf$. Here, we identify $\partial_1^\tang$, $\partial_2^\tang$, $\partial_3^\tang$ with $\partial_x^\tang$, $\partial_y^\tang$, $\partial_z^\tang$, respectively. If $\surf$ is not a closed surface, \cref{eq:model_pde} may also be subject to boundary conditions, e.g., $u(\vec{x}) = g(\vec{x})$ for $\vec{x} \in \partial\surf$ and some function $g$.

\subsection{Discrete surfaces and spectral collocation}

To numerically compute with the surface $\Gamma$, we must discretize the geometry. While many choices of discrete surface representation exist, such as level sets~\cite{Sethian1999} and point clouds~\cite{Piret2012,Fuselier2013}, we choose to use surface meshes as they provide a flexible level of compatibility with computer-aided design (CAD) software. Specifically, we assume that $\surf$ has been provided as a high-order quadrilateral\footnote{We use quadrilateral elements throughout this work for simplicity and ease of implementation, but the algorithms described below can be modified to use triangular elements instead.} surface mesh consisting of a set of tensor-product elements $\{\elem_k\}_{k=1}^\nelem$. Such a mesh may be constructed through the use of high-order meshing software such as Gmsh~\cite{Geuzaine2009} or Rhinoceros~\cite{Rhino}.

\begin{figure}[htb]
	\centering
	\vspace{0.7em}
	\begin{overpic}[width=0.56\textwidth]{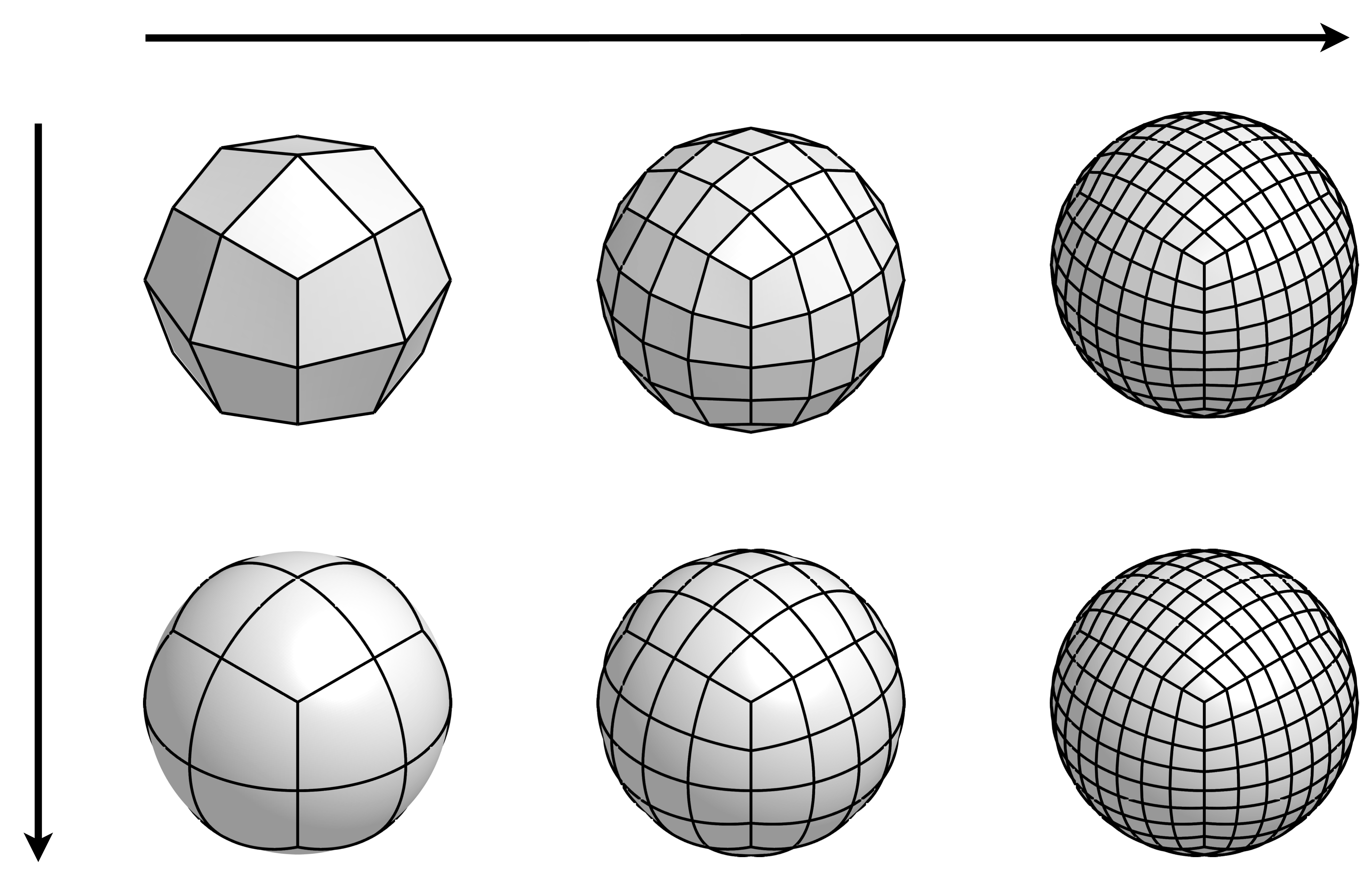}
		\put(-11,45) {\small\sffamily \parbox{1cm}{\centering Low \\ order}}
		\put(-11,10) {\small\sffamily \parbox{1cm}{\centering High \\ order}}
		\put(13,64) {\small\sffamily Coarse}
		\put(85,64) {\small\sffamily Fine}
	\end{overpic}
	~~
	\raisebox{0.8cm}{%
		\begin{overpic}[width=0.29\textwidth]{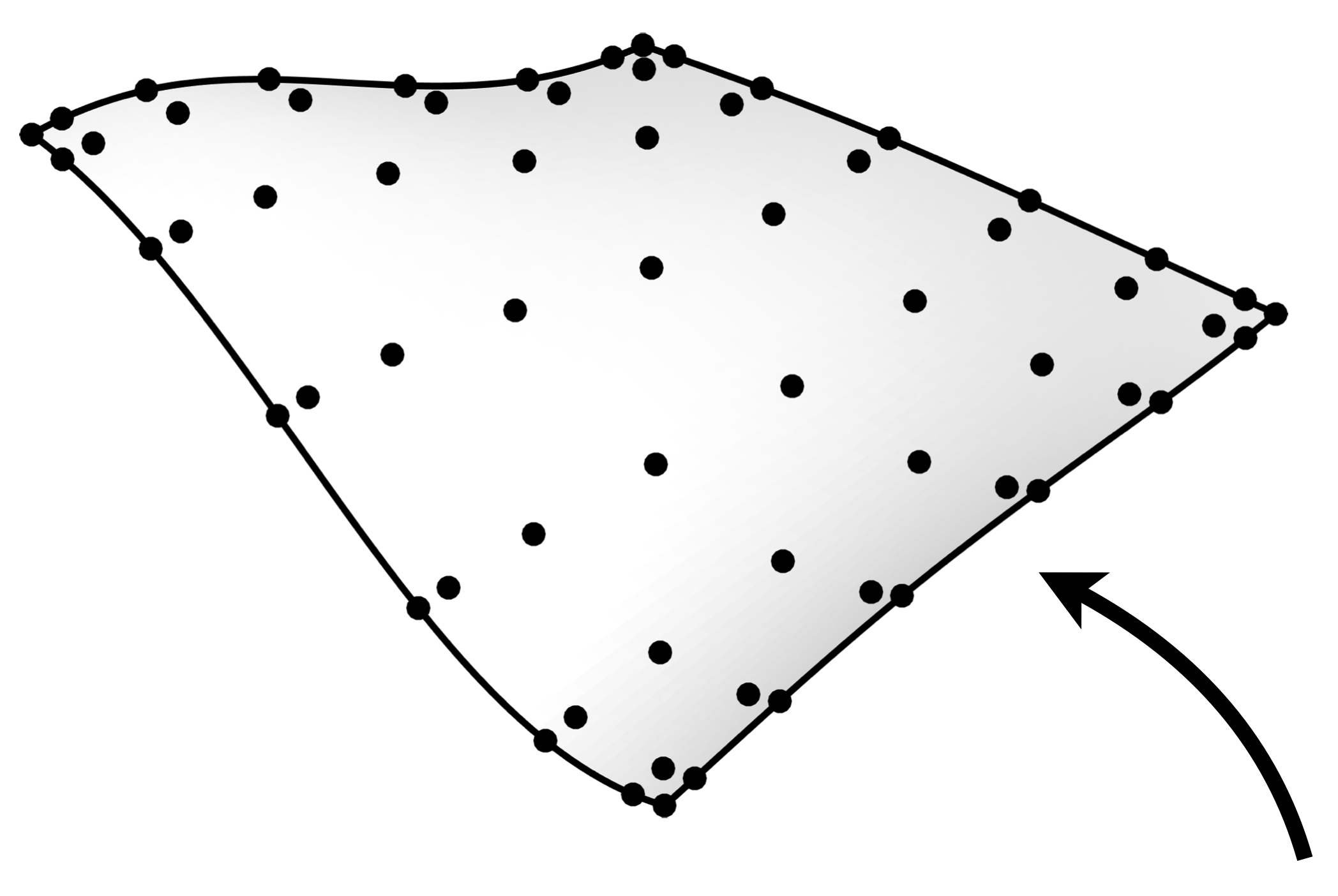}%
			\put(50,-10) {\small\sffamily \parbox{3cm}{\centering Tensor-product \\ Chebyshev nodes}}
		\end{overpic}%
	}
	\caption{(Left) Surface meshes may be constructed and refined along two axes: a coarse mesh may be subdivided into a fine mesh with many elements and the geometry may be represented using low-order or high-order polynomials. (Right) We represent each element using tensor product Chebyshev nodes of the second kind.}
	\label{fig:mesh_order}
\end{figure}
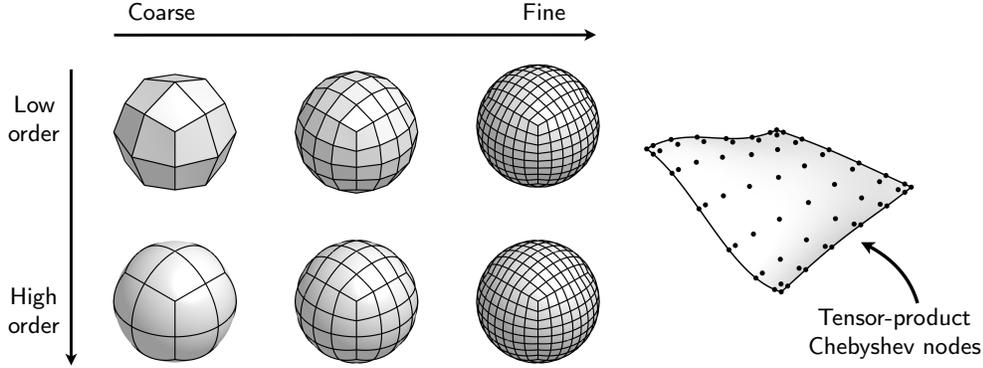

Our high-order format is based on tensor-product Chebyshev nodes, though other node sets could suffice. Let $\{\vec{\u}_{ij} = (\u_{ij},\v_{ij})\}_{i,j=1}^{p+1}$ be the set of tensor-product second-kind Chebyshev nodes of order $p$ over the reference square $[-1,1]^2$. We say that the mesh is order $p$ if each element has been provided as a set of nodes $\{\vec{x}_{ij} = (x_{ij}, y_{ij}, z_{ij})\}_{i,j=1}^{p+1}$ sampled at $\{\vec{\u}_{ij}\}_{i,j=1}^{p+1}$ on $\surf$. That is, the coordinate mapping for each element is a function $\vec{x}(\u,\v) : [-1,1]^2 \mapsto \mathbb{R}^3$ such that $\vec{x}(\u_{ij},\v_{ij}) = (x_{ij},y_{ij},z_{ij})$ for $i, j = 1, \ldots, p+1$ (see \cref{fig:element_metric}). Thus, on each element we can numerically approximate the coordinate mapping via interpolation through the nodes $\{\vec{x}_{ij}\}_{i,j=1}^{p+1}$:
\begin{equation}\label{eq:coord_approx}
\vec{x}(\u,\v) \approx \sum_{i=1}^{p+1} \sum_{j=1}^{p+1} \vec{x}_{ij} \, \ell_j(\u) \, \ell_i(\v), \qquad (\u,\v) \in [-1,1]^2,
\end{equation}
where $\ell_j$ is the $j$th Lagrange polynomial associated with the second-kind Chebyshev nodes. Partial derivatives of the elemental coordinate maps, $\vec{x}_\u$ and $\vec{x}_\v$, may then be computed through numerical spectral differentiation~\cite{Trefethen2000} and used to form an approximation to the metric tensor, $g$, on each element.

We now describe a spectral collocation method to discretize \cref{eq:model_pde} on a single element $\elem_k$. To approximate functions defined on each element we use an isoparametric representation, wherein the function is approximated on the same order-$p$ nodes as the element geometry. For a function $u(\vec{x})$ defined for $\vec{x} \in \elem_k$, denote by $u_{ij} \approx u(\vec{x}_{ij})$; that is, $u_{ij}$ is simply the function $u$ sampled on the grid. Then, just like \cref{eq:coord_approx}, we have that
\[
u(\u,\v) \approx \sum_{i=1}^{p+1} \sum_{j=1}^{p+1} u_{ij} \, \ell_j(\u) \, \ell_i(\v), \qquad (\u,\v) \in [-1,1]^2,
\]
where we have introduced the slight abuse of notation $u(\u,\v) \definedas u(\vec{x}(\u,\v))$.

To discretize $\mathcal{L}_\surf$ on the element $\elem_k$, we will compute discrete operators on the reference square $[-1,1]^2$ and then map them to $\elem_k$ using the numerical coordinate mapping. Let $\disc{D} \in \mathbb{C}^{(p+1) \times (p+1)}$ be the one-dimensional spectral differentiation matrix~\cite{Trefethen2000} associated with second-kind Chebyshev points on the interval $[-1,1]$ and let $\disc{I} \in \mathbb{C}^{(p+1) \times (p+1)}$ be the identity matrix. Two-dimensional differentiation matrices on the reference square $[-1,1]^2$ in the $\u$- and $\v$-directions can be constructed through the Kronecker products
\[
\disc{D}_\u = \disc{D} \otimes \disc{I}, \qquad \disc{D}_\v = \disc{I} \otimes \disc{D},
\]
respectively, where $\disc{D}_\u, \disc{D}_\v \in \mathbb{C}^{(p+1)^2 \times (p+1)^2}$. Let $\multmat{u} \in \mathbb{C}^{(p+1)^2 \times (p+1)^2}$ denote the diagonal multiplication matrix formed by placing the entries of $u_{ij}$ along the diagonal; we refer to $\multmat{u}$ as a multiplication matrix because $\multmat{u} v_{ij}$ computes the pointwise product of the functions $u$ and $v$ on the grid. Using \cref{eq:cont_diff_xyz}, differentiation matrices corresponding to the components of the surface gradient on $\elem_k$ are given by
\begin{equation}\label{eq:discrete_diff}
\begin{aligned}
\disc{D}_x^\tang = \multmat{\u_x} \, \disc{D}_\u + \multmat{\v_x} \, \disc{D}_\v, \\
\disc{D}_y^\tang = \multmat{\u_y} \, \disc{D}_\u + \multmat{\v_y} \, \disc{D}_\v, \\
\disc{D}_z^\tang = \multmat{\u_z} \, \disc{D}_\u + \multmat{\v_z} \, \disc{D}_\v,
\end{aligned}
\end{equation}
respectively, where the values $\{(\u_\vec{x})_{ij}\}$ and $\{(\v_\vec{x})_{ij}\}$ are tabulated using the numerical approximations to $\vec{x}_\u$, $\vec{x}_\v$, and $g$ in \cref{eq:xi_eta_derivs}. The operator $\mathcal{L}_\surf$ given in \cref{eq:L_definition} is then discretized on the element $\elem_k$ as the $(p+1)^2 \times (p+1)^2$ matrix
\begin{equation}\label{eq:discrete_L}
\begin{aligned}
\disc{L}_{\elem_k} = \sum_{i=1}^3 \sum_{j=i}^3 \multmat{a_{ij}} \disc{D}_i^\tang \disc{D}_j^\tang + \sum_{i=1}^3 \multmat{b_i} \disc{D}_i^\tang + \multmat{c},
\end{aligned}
\end{equation}
with all variable coefficients sampled on the grid. Again, we identify $D_1^\tang$, $D_2^\tang$, $D_3^\tang$ with $D_x^\tang$, $D_y^\tang$, $D_z^\tang$, respectively. Discrete surface gradient, surface divergence, and Laplace--Beltrami operators may be similarly defined. For instance, the discrete Laplace--Beltrami operator may be discretized as $\lapbel \approx \left(\disc{D}_x^\tang\right)^2 + \left(\disc{D}_y^\tang\right)^2 + \left(\disc{D}_z^\tang\right)^2$.

\section{Domain decomposition on surfaces}\label{sec:dd}

\subsection{A single element}
Consider the single-element boundary value problem on the element $\elem_k$,
\begin{subequations}
\begin{align}\label{eq:single_bvp}
\mathcal{L}_\surf u(\vec{x}) &= f(\vec{x}), \qquad \vec{x} \in \elem_k, \\
u(\vec{x}) &= g(\vec{x}), \qquad \vec{x} \in \partial\elem_k. \label{eq:single_bvp_bc}
\end{align}
\end{subequations}
Using~\cref{eq:discrete_L}, we discretize \cref{eq:single_bvp} as the $(p+1)^2 \times (p+1)^2$ linear system
\begin{equation}\label{eq:L_unordered}
L_{\elem_k} \dvec{u} = \dvec{f},
\end{equation}
where $\dvec{u}, \dvec{f} \in \mathbb{C}^{(p+1)^2 \times 1}$ are the vectors formed by stacking the values $\{u_{ij}\}$ and $\{f_{ij}\}$ column-wise. To impose the boundary condition~\cref{eq:single_bvp_bc} we use the approach of ``boundary bordering,'' whereby rows of the differential operator corresponding to boundary nodes are replaced by rows enforcing the boundary conditions. While the technique of rectangular spectral collocation~\cite{Driscoll2015} provides a robust way to impose general boundary conditions in spectral collocation methods without row deletion, a simple boundary bordering approach suffices for the Dirichlet problems considered here.

To illustrate which rows of~\cref{eq:L_unordered} correspond to boundary nodes, it is helpful to reorder the degrees of freedom. To that end, let $\indexset_\ii$ and $\indexset_\ee$ be the sets of (linear) indices for which the Chebyshev nodes $\{\vec{\u}_{ij}\}$ lie in the interior of the reference square or on its boundary, respectively, and denote their sizes by $n_\ii = \lvert \indexset_\ii \rvert = (p-1)^2$ and $n_\ee = \lvert \indexset_\ee \rvert = 4p$. For a matrix or vector, we will use superscripts to denote indices for row and column slicing; e.g., for a matrix $A$, $A^{\ii\ee}$ is the submatrix $A(\indexset_\ii, \indexset_\ee)$. Reordering the degrees of freedom in \cref{eq:L_unordered} in the order $\{\indexset_\ii, \indexset_\ee\}$ gives a block linear system,
\begin{equation}\label{eq:L_blocks}
\begin{bmatrix*}[l]
L_{\elem_k}^{\ii\ii} & L_{\elem_k}^{\ii\ee} \\[0.3em]
L_{\elem_k}^{\ee\ii} & L_{\elem_k}^{\ee\ee}
\end{bmatrix*}
\begin{bmatrix*}[l] \pad\dvec{u}^\ii \\[0.3em] \pad\dvec{u}^\ee \end{bmatrix*}
=
\begin{bmatrix*}[l] \pad\dvec{f}^\ii \\[0.3em] \pad\dvec{f}^\ee \end{bmatrix*}.
\end{equation}
We now proceed by replacing the last $n_\ee$ rows of \cref{eq:L_blocks} with the boundary conditions~\cref{eq:single_bvp_bc}, yielding the modified linear system
\begin{equation}\label{eq:L_blocks2}
\begin{bmatrix}
L_{\elem_k}^{\ii\ii} & L_{\elem_k}^{\ii\ee} \\[0.3em]
0 & I_{n_\ee}
\end{bmatrix}
\begin{bmatrix*}[l] \pad\dvec{u}^\ii \\[0.3em] \pad\dvec{u}^\ee \end{bmatrix*}
=
\begin{bmatrix} \pad\dvec{f}^\ii \\[0.3em] \pad\dvec{g} \end{bmatrix},
\end{equation}
where $\dvec{g} \in \mathbb{C}^{n_\ee \times 1}$ is the vector of boundary values of $g$ and $I_{n_\ee}$ is the $n_\ee \times n_\ee$ identity matrix. Next, we eliminate $\dvec{u}^\ee$ from this system by performing a Schur complement. In this case, as~\cref{eq:single_bvp_bc} imposes Dirichlet boundary conditions, we immediately obtain that $\dvec{u}^\ee = \dvec{g}$. Substituting this into \cref{eq:L_blocks2} yields a reduced $n_\ii \times n_\ii$ linear system for the interior unknowns $\dvec{u}^\ii$,
\begin{equation}\label{eq:L_schur}
L_{\elem_k}^{\ii\ii} \dvec{u}^\ii = \dvec{f}^\ii - L_{\elem_k}^{\ii\ee} \dvec{g}.
\end{equation}
Upon solving~\cref{eq:L_schur}, the solution values at all Chebyshev nodes on $\elem_k$ have now been computed. It is useful to note a few things:
\vspace{1em}
\begin{itemize}
\item First, we may write the solution $u$ as the sum of a homogeneous solution, $w$, and a particular solution, $v$,
\[
u = w + v,
\]
where $w$ satisfies
\begin{equation}\label{eq:homo}
\begin{aligned}
\mathcal{L}_\surf w(\vec{x}) &= 0, && \vec{x} \in \elem_k, \\
w(\vec{x}) &= g(\vec{x}), && \vec{x} \in \partial\elem_k,
\end{aligned}
\end{equation}
and $v$ satisfies
\begin{equation}\label{eq:part}
\begin{aligned}
\mathcal{L}_\surf v(\vec{x}) &= f(\vec{x}), && \vec{x} \in \elem_k, \\
v(\vec{x}) &= 0, && \vec{x} \in \partial\elem_k.
\end{aligned}
\end{equation}
\Cref{eq:L_schur} then becomes
\begin{equation}\label{eq:homo_and_part}
\begin{aligned}
\dvec{w}^\ii &= -\left(L_{\elem_k}^{\ii\ii}\right)^{-1}\!L_{\elem_k}^{\ii\ee} \dvec{g}, &\qquad \dvec{v}^\ii &= \left(L_{\elem_k}^{\ii\ii}\right)^{-1} \dvec{f}^\ii, \\
\dvec{w}^\ee &= \dvec{g}, & \dvec{v}^\ee &= 0,
\end{aligned}
\end{equation}
where $\dvec{u} = \dvec{w} + \dvec{v}$.

\vspace{1em}
\item Second, we may encode as a linear operator the action of solving~\cref{eq:homo} given any boundary data $g$. This operator, often called the ``solution operator,'' is discretized as the $(p+1)^2 \times n_\ee$ matrix
\begin{equation}\label{eq:sol_op1}
S_{\elem_k}(\indexset_\ii,:) = -\left(L_{\elem_k}^{\ii\ii}\right)^{-1}\!L_{\elem_k}^{\ii\ee},
\qquad
S_{\elem_k}(\indexset_\ee,:) = I_{n_\ee}.
\end{equation}
Given \emph{any} vector of boundary data $\dvec{g}$, the product $\dvec{w} = S_{\elem_k} \dvec{g}$ approximates $w$ satisfying~\cref{eq:homo}.
\end{itemize}

\begin{remark}\label{rem:corners}
From a domain decomposition point of view it is convenient to work on a boundary grid which does not include the four corner nodes, so operators that act on the boundary degrees of freedom (such as $S_{\elem_k}$) may be unambiguously separated into pieces that act on each side of an element~\cite{Gillman2014}. To achieve this effect, we choose the boundary grid on each edge to be the first-kind Chebyshev nodes and apply mappings to all boundary operators to re-interpolate to this node set. For stability, we use a boundary grid of degree $p-2$ on each side~\cite{Gillman2015}; thus, the number of boundary nodes is modified to be $n_\ee = 4(p-1)$.
\end{remark}

We now have all the ingredients necessary to locally solve a surface PDE on each element of the mesh. To couple elements together, we take a non-overlapping domain decomposition approach. We begin with the simplest possible domain decomposition problem, consisting of just two surface elements.

\subsection{Two ``glued'' elements}\label{sec:two_glued_patches}
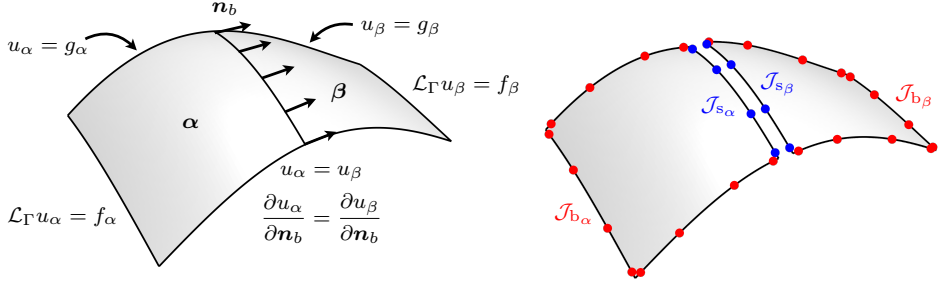
\begin{figure}[htb]
	\centering
	~~~~
	\begin{overpic}[width=0.435\textwidth]{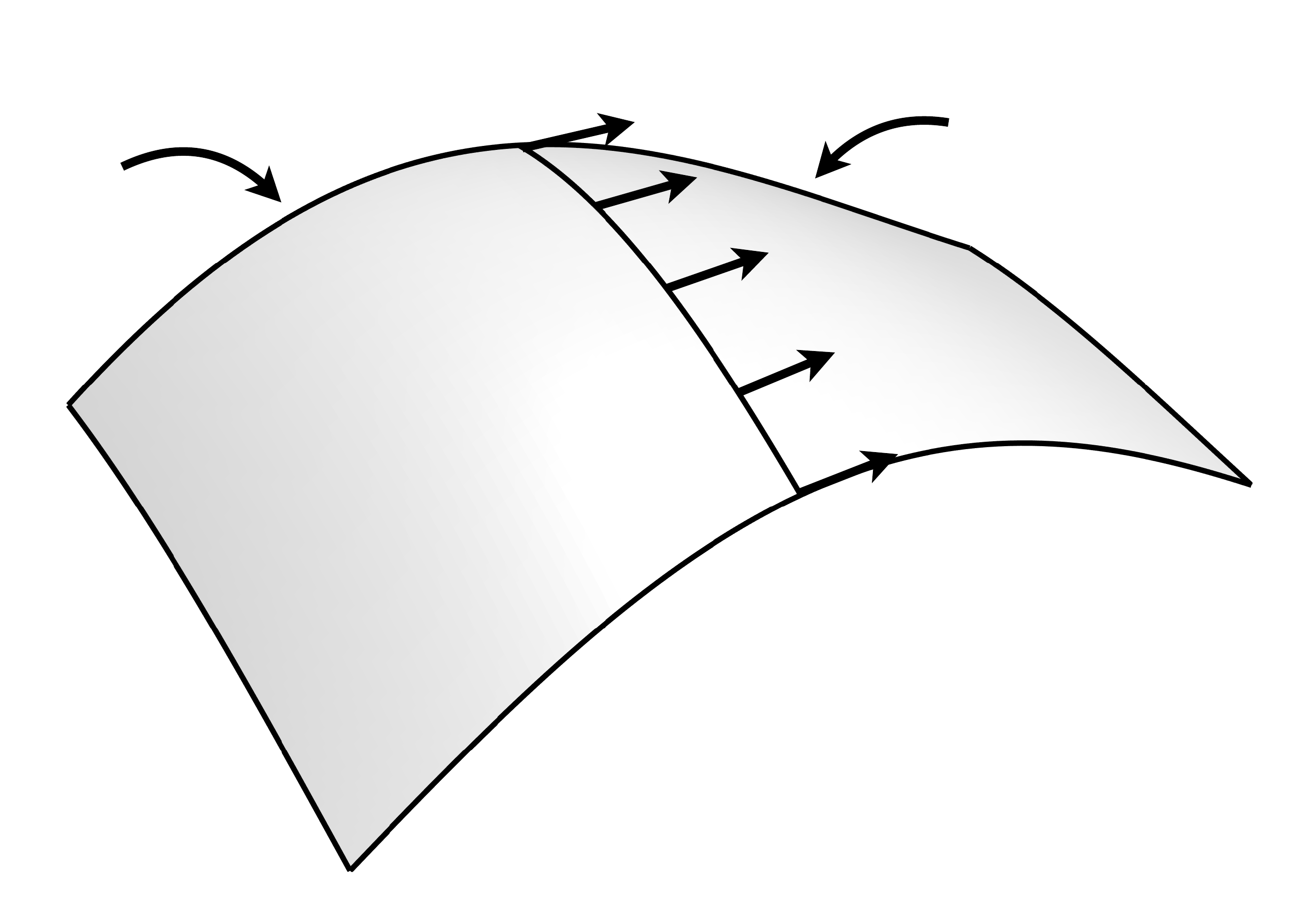}
	\put(32,36){\footnotesize$\bm{\ea}$}
	\put(67,43){\footnotesize$\bm{\eb}$}
	\put(-9,54){\footnotesize$u_\ea = g_\ea$}
	\put(74,59){\footnotesize$u_\eb = g_\eb$}
	\put(-9,14){\footnotesize$\mathcal{L}_\surf u_\ea = f_\ea$}
	\put(86,45){\footnotesize$\mathcal{L}_\surf u_\eb = f_\eb$}
	\put(39,63.3){\footnotesize$\binorm$}
	\put(55,25){\footnotesize$u_\ea = u_\eb$}
	\put(50.4,13.5){\footnotesize$\dfrac{\partial u_\ea}{\partial \binorm} = \dfrac{\partial u_\eb}{\partial \binorm}$}
	\end{overpic}
	~~~~~~
	\begin{overpic}[width=0.405\textwidth]{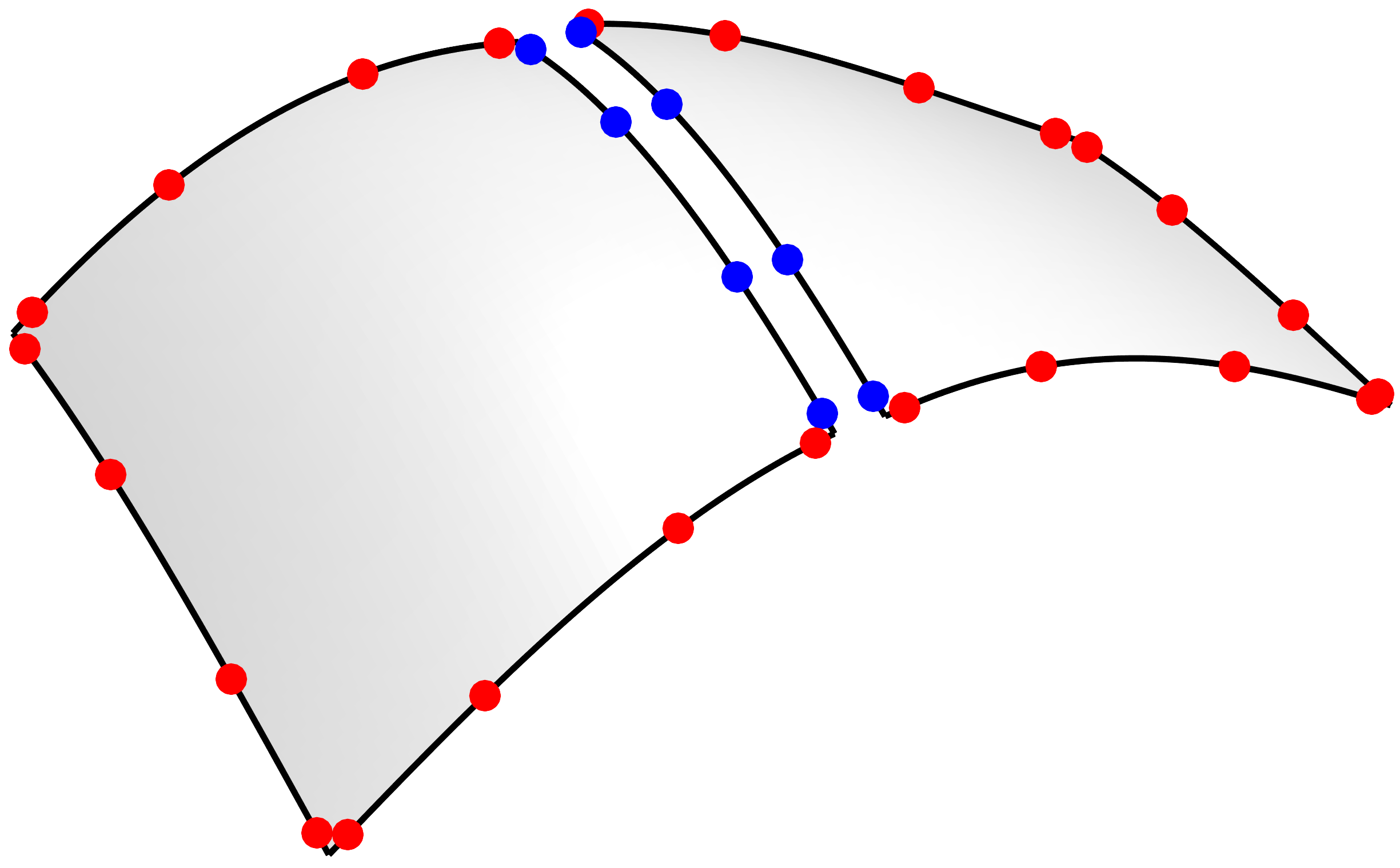}
	\put(3,16) {\footnotesize\color{red}$\indexset_{\ee_\ea}$}
	\put(89,46){\footnotesize\color{red}$\indexset_{\ee_\eb}$}
	\put(40,43){\footnotesize\color{blue}$\indexset_{\ss_\ea}$}
	\put(55,49){\footnotesize\color{blue}$\indexset_{\ss_\eb}$}
	\end{overpic}
	\caption{(Left) Two patches are ``glued'' together by imposing continuity and continuity of the binormal derivative across the interface between them. (Right) Indices corresponding to the shared and unshared boundary degrees of freedom are recorded for each patch. Values at the shared interface can be eliminated using the Schur complement method.}
	\label{fig:two_glued_patches}
\end{figure}

Consider the model problem of two ``glued'' surface elements, $\elem_\ea$ and $\elem_\eb$, separated by an interface $\interface = \elem_\ea \cap \elem_\eb$, depicted in \cref{fig:two_glued_patches}~(left):
\begin{equation}\label{eq:two_glued_patches}
\begin{aligned}
\mathcal{L}_\surf u_\ea(\vec{x}) &= f_\ea(\vec{x}), \qquad&&\vec{x} \in \elem_\ea, \\
\mathcal{L}_\surf u_\eb(\vec{x}) &= f_\eb(\vec{x}), &&\vec{x} \in \elem_\eb, \\
u_\ea(\vec{x}) &= g_\ea(\vec{x}), &&\vec{x} \in \partial\elem_\ea \setminus \interface, \\
u_\eb(\vec{x}) &= g_\eb(\vec{x}), &&\vec{x} \in \partial\elem_\eb \setminus \interface, \\
u_\ea(\vec{x}) &= u_\eb(\vec{x}), &&\vec{x} \in \interface, \\
\tfrac{\partial u_\ea}{\partial \binorm}(\vec{x}) &= \tfrac{\partial u_\eb}{\partial \binorm}(\vec{x}), &&\vec{x} \in \interface,
\end{aligned}
\end{equation}
where $\binorm$ is the binormal vector along the interface. As before, we write $u_\ea$ and $u_\eb$ as the sum of homogeneous and particular solutions, $u_\ea = w_\ea + v_\ea$ and $u_\eb = w_\eb + v_\eb$, with discrete versions $\dvec{u}_\ea = \dvec{w}_\ea + \dvec{v}_\ea$ and $\dvec{u}_\eb = \dvec{w}_\eb + \dvec{v}_\eb$, where $w_\ea$ and $w_\eb$ satisfy~\cref{eq:two_glued_patches} with $f_\ea \equiv f_\eb \equiv 0$ and $v_\ea$ and $v_\eb$ satisfy~\cref{eq:two_glued_patches} with $g_\ea \equiv g_\eb \equiv 0$.

Let $\indexset_{\ss_\ea}, \indexset_{\ss_\eb} \subset \{1, \ldots, n_\ee\}$ be index sets corresponding to the shared interface nodes (i.e., the nodes which are interior to the merged domain $\elem_\ea \cup \elem_\eb$) with respect to elements $\elem_\ea$ and $\elem_\eb$, and $\indexset_{\ee_\ea}$ and $\indexset_{\ee_\eb}$ the remaining indices corresponding to the boundary nodes. See \cref{fig:two_glued_patches} (right) for a diagram. As a reminder, we use superscripts to denote index slicing. As the nodes of $\elem_\ea$ and $\elem_\eb$ are identical along the shared interface $\interface$, continuity of the solution across $\interface$ simply means that $\dvec{u}_\ea^{\ss_\ea} = \dvec{u}_\eb^{\ss_\eb}$, so let us denote these solution values by $\dvec{u}_\interface$. To enforce continuity of the binormal derivative, let us evaluate the outward pointing binormal vectors of each element $\elem_k$ at the $n_\ee$ boundary nodes and collect the Cartesian components of each into the discrete vectors $\dvec{n}_{k,x}$, $\dvec{n}_{k,y}$, and $\dvec{n}_{k,z}$, for $k=\ea,\eb$. Then, we may define binormal derivative operators $D_{\elem_\ea}$ and $D_{\elem_\eb}$ that map function values in each element to the values of the outgoing flux on the boundary via
\begin{equation}\label{eq:normal_d}
D_{\elem_k} = \dvec{n}_{k,x} \circ D_x^\surf(\indexset_\ee, :) + \dvec{n}_{k,y} \circ D_y^\surf(\indexset_\ee, :) + \dvec{n}_{k,z} \circ D_z^\surf(\indexset_\ee, :),
\end{equation}
for $k=\ea,\eb$, where $\circ$ denotes the Hadamard product. The derivatives of $u_\ea$ and $u_\eb$ in the direction $\vec{n}_b$ along the interface can then be written as $(D_{\elem_\ea} \dvec{u}_\ea)^{\ss_\ea}$ and $-(D_{\elem_\eb} \dvec{u}_\eb)^{\ss_\eb}$, respectively. Therefore, continuity of the binormal derivative across the interface is equivalent to
\begin{equation}\label{eq:cont_d}
(D_{\elem_\ea} \dvec{u}_\ea)^{\ss_\ea} + (D_{\elem_\eb} \dvec{u}_\eb)^{\ss_\eb} = 0.
\end{equation}

Now, let $S_{\elem_\ea}$ and $S_{\elem_\eb}$ be the solution operators for elements $\elem_\ea$ and $\elem_\eb$, respectively. If we knew the value of $\dvec{u}_\interface$, then boundary data on all sides of each element would be specified and we could directly apply the solution operators to each element in an uncoupled fashion. So let $\tilde{\dvec{g}}_\ea, \tilde{\dvec{g}}_\eb \in \mathbb{C}^{n_\ee \times 1}$ be the vectors of boundary data on all sides of each element, with entries given by $\tilde{\dvec{g}}_\ea^{\smash{\ss_\ea}} = \tilde{\dvec{g}}_\eb^{\smash{\ss_\eb}} = \dvec{u}_\interface$, $\tilde{\dvec{g}}_\ea^{\smash{\ee_\ea}} = \dvec{g}_\ea$, and $\tilde{\dvec{g}}_\eb^{\smash{\ee_\eb}} = \dvec{g}_\eb$. With this definition, the homogeneous solutions on each element can be written as $\dvec{w}_\ea = S_{\elem_\ea} \tilde{\dvec{g}}_\ea$ and $\dvec{w}_\eb = S_{\elem_\eb} \tilde{\dvec{g}}_\eb$. Expanding $\dvec{u}_\ea = \dvec{w}_\ea + \dvec{v}_\ea$ and $\dvec{u}_\eb = \dvec{w}_\eb + \dvec{v}_\eb$ in~\cref{eq:cont_d}, we have
\[
(D_{\elem_\ea} S_{\elem_\ea} \tilde{\dvec{g}}_\ea + D_{\elem_\ea} \dvec{v}_\ea)^{\ss_\ea} + (D_{\elem_\eb} S_{\elem_\eb} \tilde{\dvec{g}}_\eb + D_{\elem_\eb} \dvec{v}_\eb)^{\ss_\eb} = 0.
\vspace{0.2em}
\]
The matrices $\DtN_{\elem_\ea} \definedas D_{\elem_\ea} S_{\elem_\ea}$ and $\DtN_{\elem_\eb} \definedas D_{\elem_\eb} S_{\elem_\eb}$ are discrete \emph{Dirichlet-to-Neumann} maps, so termed because they map Dirichlet data on the boundary of $\elem_k$ to the outward flux of the local solution to the PDE on the boundary of $\elem_k$. The vectors $\dvec{v}_\ea^{\prime} \definedas D_{\elem_\ea} \dvec{v}_\ea$ and $\dvec{v}_\eb^{\prime} \definedas D_{\elem_\eb} \dvec{v}_\eb$ are the ``particular fluxes''---the outgoing fluxes of the particular solution evaluated on the boundary of each element. Using this notation and recalling the definitions of $\tilde{\dvec{g}}_\ea$ and $\tilde{\dvec{g}}_\eb$, continuity of the binormal derivative across the interface yields an equation for the unknown vector of interface values $\dvec{u}_\interface$,
\[
\left(\DtN_{\elem_\ea}^{\nudge \ss_\ea \nudge \ss_\ea} \dvec{u}_\interface + \DtN_{\elem_\ea}^{\nudge \ss_\ea \nudge \ee_\ea} \dvec{g}_\ea + {\dvec{v}_\ea^{\prime}}^{\!\ss_\ea}\right)
+
\left(\DtN_{\elem_\eb}^{\nudge \ss_\eb \nudge \ss_\eb} \dvec{u}_\interface + \DtN_{\elem_\eb}^{\nudge \ss_\eb \nudge \ee_\eb} \dvec{g}_\eb + {\dvec{v}_\eb^{\prime}}^{\!\ss_\eb}\right) = 0.
\]
Writing $\dvec{u}_\interface$ as the sum of homogeneous and particular solutions, $\dvec{u}_\interface = \dvec{v}_\interface + \dvec{w}_\interface$, and substituting into the above, yields linear systems for the homogeneous and particular interface unknowns,
\begin{subequations}
\begin{alignat}{2}
\label{eq:w_interface}
&-\Bigl( \DtN_{\elem_\ea}^{\nudge \ss_\ea \nudge \ss_\ea} + \DtN_{\elem_\eb}^{\nudge \ss_\eb \nudge \ss_\eb} \Bigr) \, \dvec{w}_\interface &&=\, \DtN_{\elem_\ea}^{\nudge \ss_\ea \nudge \ee_\ea} \dvec{g}_\ea + \, \DtN_{\elem_\eb}^{\nudge \ss_\eb \nudge \ee_\eb} \dvec{g}_\eb, \\
\label{eq:v_interface}
&-\Bigl( \DtN_{\elem_\ea}^{\nudge \ss_\ea \nudge \ss_\ea} + \DtN_{\elem_\eb}^{\nudge \ss_\eb \nudge \ss_\eb} \Bigr) \, \dvec{v}_\interface &&=\, {\dvec{v}_\ea^{\prime}}^{\!\ss_\ea} + \, {\dvec{v}_\eb^{\prime}}^{\!\ss_\eb}.
\end{alignat}
\end{subequations}
Upon solution of \cref{eq:v_interface} and \cref{eq:w_interface}, boundary data is known on all sides of both elements $\elem_\ea$ and $\elem_\eb$, and the local solution operators $S_{\elem_\ea}$ and $S_{\elem_\eb}$ may be applied to recover the solutions $\dvec{u}_\ea$ and $\dvec{u}_\eb$ to \cref{eq:two_glued_patches} in the interiors of each element. As before, we may encode as a linear operator the action of solving~\cref{eq:w_interface} given any boundary data $\dvec{g}_\ea$ and $\dvec{g}_\eb$. This is the solution operator for the interface, $S_\interface$, and is given by the solution to
\begin{equation}\label{eq:S_interface}
-\Bigl( \DtN_{\elem_\ea}^{\nudge \ss_\ea \nudge \ss_\ea} + \DtN_{\elem_\eb}^{\nudge \ss_\eb \nudge \ss_\eb} \Bigr) S_\interface = \Bigl[ \, \DtN_{\elem_\ea}^{\nudge \ss_\ea \nudge \ee_\ea} \;\;\; \DtN_{\elem_\eb}^{\nudge \ss_\eb \nudge \ee_\eb} \Bigr].
\end{equation}
The operator $S_\interface$ takes in Dirichlet data on the boundary nodes of the merged domain $\elem_{\ea\eb} \definedas \elem_\ea \cup \elem_\eb$ and returns the values of the homogeneous solution along the interface, $\dvec{w}_\interface$.

Note that when ``gluing'' two elements, computing $S_\interface$ and $\dvec{v}_\interface$ only requires knowledge of the Dirichlet-to-Neumann map and particular flux for each element. Thus, obtaining the Dirichlet-to-Neumann map and particular flux of the merged domain will allow us to recursively apply this ``gluing'' procedure. The matrix $\DtN_{\elem_\ea}^{\nudge \ss_\ea \nudge \ss_\ea} + \DtN_{\elem_\eb}^{\nudge \ss_\eb \nudge \ss_\eb}$ may be viewed as a Schur complement of the block linear system corresponding to direct discretization of~\cref{eq:two_glued_patches}, taken with respect to the interface degrees of freedom. The Schur complement also allows us to construct the Dirichlet-to-Neumann map for the merged domain $\elem_{\ea\eb}$, using the Dirichlet-to-Neumann maps for the elements $\elem_\ea$ and $\elem_\eb$ and the solution operator for the interface, $S_\interface$, as
\begin{equation}\label{eq:DtN_interface}
\DtN_{\elem_{\ea\eb}} =
\left[\!
\begin{array}{cc}
\\[-1em]
\DtN_{\elem_\ea}^{\nudge\ee_\ea\nudge\ee_\ea} & 0 \\[0.4em]
0 & \DtN_{\elem_\eb}^{\nudge\ee_\eb\nudge\ee_\eb} \\[0.2em]
\end{array}
\!\!\right]
+
\left[\!
\begin{array}{c}
\\[-1em]
\DtN_{\elem_\ea}^{\nudge\ee_\ea\nudge\ss_\ea} \\[0.4em]
\DtN_{\elem_\eb}^{\nudge\ee_\eb\nudge\ss_\eb} \\[0.2em]
\end{array}
\!\!\right]
S_\interface.
\end{equation}
Similarly, the particular flux for the merged domain is given by
\begin{equation}\label{eq:vprime_interface}
\dvec{v}_{\ea\eb}^{\prime} = \left[\!
\begin{array}{c}
\\[-1em]
{\dvec{v}_\ea^{\prime}}^{\!\ee_\ea} \\[0.4em]
{\dvec{v}_\eb^{\prime}}^{\!\ee_\eb} \\ [0.2em]
\end{array}
\!\!\right]
+
\left[\!
\begin{array}{c}
\\[-1em]
\DtN_{\elem_\ea}^{\nudge\ee_\ea\nudge\ss_\ea} \\[0.4em]
\DtN_{\elem_\eb}^{\nudge\ee_\eb\nudge\ss_\eb} \\[0.2em]
\end{array}
\!\!\right]
\dvec{v}_\interface.
\end{equation}

\section{A fast direct solver}\label{sec:hps}

We now formulate the general hierarchical fast direct solver for a mesh of $\nelem$ elements, $\{\elem_k\}_{k=1}^\nelem$. The solver is based on the the HPS scheme~\cite{Martinsson2013, Gillman2014, Gillman2015, Babb2018, Fortunato2021, Martinsson2019a} and proceeds pairwise by recursively merging neighboring elements or groups of elements until a factorization of the surface PDE over entire surface mesh has been computed. We refer to this method as the surface HPS scheme.

\subsection{Local factorization stage}
The method begins with the local factorization stage, in which a solution operator $S_{\elem_k}$, Dirichlet-to-Neumann map $\DtN_{\elem_k}$, particular solution $\dvec{v}_k$, and particular flux $\dvec{v}_k^{\prime}$ are computed for every element $\elem_k$ in the mesh, $k = 1, \ldots, \nelem$. As the construction of these quantities is entirely local to each element, this stage can be trivially parallelized across the mesh. The local factorization stage is outlined is \cref{alg:init}.

\begin{algorithm}[htb]
\caption{Local factorization stage}
\begin{algorithmic}[1]
\Require{Mesh $\{\elem_k\}_{k=1}^\nelem$, surface differential operator $\mathcal{L}_\surf$, righthand side $f$}
\Ensure{Solution operators $\{S_{\elem_k}\}$, Dirichlet-to-Neumann maps $\{\DtN_{\elem_k}\}$, particular solutions $\{\dvec{v}_k\}$, and particular fluxes $\{\dvec{v}_k^{\prime}\}$ for every element $k = 1, \ldots, \nelem$.}%
\vspace{1em}
\State{Construct interior and boundary index sets $\indexset_\ii$, $\indexset_\ee$.}
\For{$k = 1, \ldots, \nelem$}
	\State{Evaluate $f$ at the nodes of $\elem_k$ to yield $\dvec{f}_k$.}
	\State{Compute discrete differentiation matrices, $D_x^\surf$, $D_y^\surf$, $D_z^\surf$, on $\elem_k$ using \cref{eq:discrete_diff}.}
	\State{Compute discrete operator $L_{\elem_k}$ using \cref{eq:discrete_L}.}
	\State{Compute the homogeneous solution operator $S_{\elem_k}$ such that \vspace{-0.4em}
	\[
		S_{\elem_k}(\indexset_\ii,:) = -\left(L_{\elem_k}^{\ii\ii}\right)^{-1}\!L_{\elem_k}^{\ii\ee},
\qquad
		S_{\elem_k}(\indexset_\ee,:) = I_{n_\ee}.\vspace{-0.3em}
	\]}
	\State{Compute the discrete particular solution such that $\dvec{v}_k^\ii = \left(L_{\elem_k}^{\ii\ii}\right)^{-1} \dvec{f}_k^\ii$, $\dvec{v}_k^\ee = 0$.}
	\State{Compute the Dirichlet-to-Neumann map $\DtN_{\elem_k} \definedas D_{\elem_k} S_{\elem_k}$.}
	\State{Compute the flux of the particular solution, $\dvec{v}_k^{\prime} \definedas D_{\elem_k} \dvec{v}_k$.}
	\State{Reinterpolate boundary grids of $S_{\elem_k}$, $\DtN_{\elem_k}$ according to \cref{rem:corners}.}
\EndFor
\State{\Return $\{S_{\elem_k}\}$, $\{\DtN_{\elem_k}\}$, $\{\dvec{v}_k\}$, $\{\dvec{v}_k^{\prime}\}$ for $k = 1, \ldots, \nelem$}
\end{algorithmic}
\label{alg:init}
\end{algorithm}

\subsection{Global factorization stage}

Once local operators have been computed for each element, the solver proceeds to the global factorization stage. Elements (or groups of elements) are merged pairwise in an upward pass according to a set of merge indices $\tree_\ell$, where $(\alpha,\beta) \in \tree_\ell$ indicates that elements $\elem_\alpha$ and $\elem_\beta$ should be merged at level $\ell$ in the hierarchy. Upon merging the two elements $\elem_\alpha$ and $\elem_\beta$, a new parent element $\elem_{\alpha \beta} \definedas \elem_\alpha \cup \elem_\beta$ is created which may be subsequently merged at higher levels in the tree. Thus, a complete specification of the tree is given by the level-ordered set of merge indices $\{\tree_\ell\}_{\ell=0}^\nlevel$, with $\ell=0$ corresponding to the leaf level (containing all the individual elements of the mesh) and $\ell=\nlevel$ corresponding to the top level (containing just two elements). For each merge, the interface linear systems are constructed as in \cref{sec:two_glued_patches} and inverted directly, yielding interfacial solution operators, interfacial particular solutions, merged Dirichlet-to-Neumann maps, and merged particular fluxes for the parent element. The global factorization stage is outlined in \cref{alg:build}.

\begin{remark}\label{rem:rank_def}
The Laplace--Beltrami problem $\lapbel u = f$ on a closed surface $\surf$ is rank-one deficient, but is uniquely solvable under the mean-zero conditions $\int_\surf u = \int_\surf \! f = 0$. In the present scheme this rank deficiency is only seen in the final merge at the top level in the hierarchy, i.e., $\dim( \operatorname{null} A ) = 1$ where $A = \DtN_{\elem_\alpha}^{\nudge \ss_\alpha \nudge \ss_\alpha} + \DtN_{\elem_\beta}^{\nudge \ss_\beta \nudge \ss_\beta}$ is the interface matrix at the top level. We fix this rank deficiency by adding the mean-zero constraint to the top-level linear system, yielding a modified matrix $A + \dvec{q} \dvec{q}^T$, where $\dvec{q}$ is the vector of scaled quadrature weights on the interface. Then $\dim( \operatorname{null}( A + \dvec{q} \dvec{q}^T ) ) = 0$~\cite{Sifuentes2015,ImbertGerard2017}.
\end{remark}

\begin{algorithm}[htb]
\caption{Global factorization stage (upward pass)}
\begin{algorithmic}[1]
\Require{\mbox{Mesh $\{\elem_k\}_{k=1}^\nelem$, Dirichlet-to-Neumann maps $\{\DtN_{\elem_k}\}_{k=1}^\nelem$, merge indices $\!\{\tree_\ell\}_{\ell=0}^\nlevel$}}
\Ensure{Solution operators $\{S_\interface\}$, Dirichlet-to-Neumann maps $\{\DtN_{\elem_{\alpha \beta}}\}$, particular solutions $\{\dvec{v}_\interface\}$, and particular fluxes $\{\dvec{v}_{\alpha \beta}^{\prime}\}$ for every merge $(\alpha, \beta) \in \tree_\ell$ for every level $\ell$}%
\vspace{1em}
\For{each level $\ell = 0, 1, \ldots, \nlevel$}
\For{each pair $(\alpha, \beta) \in \tree_\ell$}
	\State{Define the merged domain $\elem_{\alpha \beta} \definedas \elem_\alpha \cup \elem_\beta$.}
	\State{Define the shared interface $\interface \definedas \elem_\alpha \cap \elem_\beta$.}
	\State{Define indices $\indexset_{\ss_\alpha}$, $\indexset_{\ss_\beta}$ for the shared interface nodes on $\elem_\alpha$, $\elem_\beta$.}
	\State{Define indices $\indexset_{\ee_\alpha}$, $\indexset_{\ee_\beta}$ for the unshared boundary nodes on $\elem_\alpha$, $\elem_\beta$.}
	\State{\multiline{Solve the linear systems \cref{eq:S_interface} and \cref{eq:v_interface} for the solution operator $S_\interface$ and particular solution $\dvec{v}_\interface$ on the interface.}}
	\State{Compute the merged Dirichlet-to-Neumann map $\DtN_{\elem_{\alpha \beta}}$ via \cref{eq:DtN_interface}.}
	\State{Compute the merged particular flux $\dvec{v}_{\alpha\beta}^{\prime}$ via \cref{eq:vprime_interface}.}
\EndFor
\EndFor
\State{\Return $\{S_\interface\}$, $\{\DtN_{\elem_{\alpha \beta}}\}$, $\{\dvec{v}_\interface\}$, $\{\dvec{v}_{\alpha \beta}^{\prime}\}$ for all $(\alpha, \beta) \in \tree_\ell$ for every level $\ell$}
\end{algorithmic}
\label{alg:build}
\end{algorithm}

\subsection{Solve stage}

Once the global factorization stage is complete, the algorithm enters the solve stage. For a surface with a boundary, Dirichlet data $g$ is evaluated at the boundary nodes and passed to the top-level solution operator, where it is used recover the solution on the top-level interface. (For a closed surface, no boundary conditions are needed, and the solution on the top-level interface is simply given by the interfacial particular solution.) This data is passed to the child elements, where this procedure continues recursively in a downward pass until the values of the solution have been tabulated on all element interfaces. The local solution operators computed in the local factorization stage may then be applied to recover the solution in the interior of each element in a parallel fashion. The solve stage is outlined in \cref{alg:solve}.

The solve stage may be repeatedly executed with different boundary data $g$ with no need to revisit the local and global factorization stages. Moreover, the particular solutions and particular fluxes throughout the hierarchy maybe be efficiently updated for different righthand sides $f$, without the need to recompute solution operators and Dirichlet-to-Neumann maps. If factorizations of the interface linear systems are stored for every merge, then one need only apply those factorizations (e.g., using forward- and back-substitution) in an upwards pass through the hierarchy in order to update the righthand side $f$.

\begin{algorithm}[htb]
\caption{Solve stage (downward pass)}
\begin{algorithmic}[1]
\Require{Merged element $\elem_{\alpha \beta}$ (or leaf element $\elem_k$), Dirichlet data $\dvec{g}_{\alpha \beta}$}
\Ensure{Solutions $\{\dvec{u}_k\}_{k=1}^\nelem$ for every element}%
\vspace{1em}
\If{$\elem_k$ is a leaf}
	\State{Compute the local solution $\dvec{u}_k \definedas S_{\elem_k} \dvec{g}_k + \dvec{v}_k$.}
	\State{\Return $\dvec{u}_k$}
\Else
	\State{\multiline{Look up the elements $\elem_\alpha$, $\elem_\beta$ that were merged to make $\elem_{\alpha \beta}$.}}
	\State{Look up the shared interface $\interface$.}
	\State{Compute the solution on the shared interface, $\dvec{u}_\interface \definedas S_\interface \dvec{g}_{\alpha \beta} + \dvec{v}_\interface$.}
	\State{\multiline{Define boundary data on the children, $\dvec{g}_\alpha, \dvec{g}_\beta$, with
	$\dvec{g}_\alpha^{\smash{\ee_\alpha}} \!\definedas \dvec{g}_{\alpha \beta}^{\smash{\ee_\alpha}}$,
	$\dvec{g}_{\smash{\beta}}^{\smash{\text{\raisebox{0.1em}{$\ee_\beta$}}}} \!\definedas \dvec{g}_{\smash{\alpha \beta}}^{\smash{\text{\raisebox{0.1em}{$\ee_\beta$}}}}$, and
	$\dvec{g}_\alpha^{\smash{\ss_\alpha}} \!= \dvec{g}_{\smash{\beta}}^{\smash{\text{\raisebox{0.15em}{$\ss_\beta$}}}} \!\definedas \dvec{u}_\interface$.\vspace{0.1em}%
	}}
	\State{Recursively compute $\{\dvec{u}_\alpha\}$ using \cref{alg:solve} with Dirichlet data $\dvec{g}_\alpha$.}
	\State{Recursively compute $\{\dvec{u}_\beta\}$ using \cref{alg:solve} with Dirichlet data $\dvec{g}_\beta$.}
	\State{\Return $\{\dvec{u}_\alpha\} \cup \{\dvec{u}_\beta\}$}
\EndIf
\end{algorithmic}
\label{alg:solve}
\end{algorithm}

\subsection{Computational complexity}

The computational complexity of the surface HPS method is identical to that of the two-dimensional HPS scheme on unstructured meshes~\cite{Fortunato2021}. We briefly summarize it now.

The total number of degrees of freedom in the order-$p$ mesh $\{\elem_k\}_{k=1}^\nelem$ scales as $\nelem p^2$. On each element $\elem_k$, construction of the solution operator $S_{\elem_k}$ requires the solution of a $\mathcal{O}(p^2) \times \mathcal{O}(p^2)$ linear system with $\mathcal{O}(p)$ righthand sides, which can be computed in $\mathcal{O}(p^6)$ operations. The Dirichlet-to-Neumann map $\DtN_{\elem_k}$ can be computed as a matrix-matrix product in $\mathcal{O}(p^5)$ operations. Therefore, the local factorization stage (\cref{alg:init}) requires $\mathcal{O}(Np^6)$ operations.

The cost of the global factorization stage depends on the merge order defined by the hierarchical indices $\{\tree_\ell\}_{\ell=0}^\nlevel$. As in classical nested dissection, the hierarchy should be as balanced as possible so that the indices $\{\tree_\ell\}_{\ell=0}^\nlevel$ define an approximately balanced binary tree with $\nlevel \sim \log{\nelem}$. A balanced partitioning of the mesh may be automatically computed by conversion to a graph partitioning problem~\cite{Karypis1998}. We assume that the merge indices $\{\tree_\ell\}_{\ell=0}^\nlevel$ define a hierarchy of $\mathcal{O}(\log \nelem)$ levels, with level $\ell$ containing $\mathcal{O}(2^{-\ell} N)$ elements. Computing the merged solution operator for a merge in level $\ell$ requires solving an interface linear system of size $\mathcal{O}(2^{\ell/2} p) \times \mathcal{O}(2^{\ell/2} p)$. Hence, the cost of performing all merges on level $\ell$ scales as $(2^{-\ell} \nelem) \cdot (2^{\ell/2} p)^3 = 2^{\ell/2} N p^3$. Summing over all $\mathcal{O}(\log \nelem)$ levels yields an overall complexity for the global factorization stage (\cref{alg:build}) of $\mathcal{O}(N^{3/2} p^3)$.

At level $\ell$ of the solve stage, the solution on the shared interface is recovered via a matrix-vector product with an $\mathcal{O}(2^{\ell/2}p) \times \mathcal{O}(2^{\ell/2}p)$ matrix in $\mathcal{O}((2^{\ell/2}p)^2)$ operations. The cost of computing this product for all interfaces on level $\ell$ scales as $(2^{-\ell} \nelem) \cdot (2^{\ell/2}p)^2 = \nelem p^2$. Hence, the solution at all interfaces on all levels of the hierarchy can be obtained in $\mathcal{O}(p^2 \nelem \log \nelem)$ operations. At the bottom level, the solution is computed in the interior of each element via multiplication with the local solution operator, requiring a total of $\mathcal{O}(\nelem p^3)$ operations. Therefore, the overall cost of the solve stage (\cref{alg:solve}) is $\mathcal{O}(p^2 \nelem \log \nelem + \nelem p^3)$.

We take a fixed-$p$ perspective throughout this paper, choosing to refine the solution by increasing the number of surface elements $\nelem$ (i.e., $h$-refinement). While $p$-refinement and $hp$-adaptivity are powerful tools for solving PDEs in complex geometries where smoothness criteria can be considered in a localized fashion~\cite{Fortunato2021}, on discrete surfaces the neighborhood of every point is only approximately smooth, making it difficult to use $hp$-adaptivity in a localized way. From this perspective then, the overall complexity of the scheme is $\mathcal{O}(\nelem)$ for the local factorization stage, $\mathcal{O}(\nelem^{3/2})$ for the global factorization stage, and $\mathcal{O}(\nelem \log \nelem)$ for the solve stage. Updating the particular solution for a new righthand $f$ takes $\mathcal{O}(\nelem \log \nelem)$ time.

The memory complexity mirrors that of the solve stage, as every level of the hierarchy stores dense solution operators and Dirichlet-to-Neumann maps. Therefore, the memory cost scales as $\mathcal{O}(\nelem \log \nelem)$.

\section{Numerical examples}\label{sec:results}

We now demonstrate the convergence and performance of the surface HPS method. The method is implemented in an open-source MATLAB package that provides abstractions for computing with functions on surfaces~\cite{GithubRepo}. Codes for the numerical examples presented below are publicly available~\cite{GithubRepo}. All simulations were run in MATLAB R2021a on a MacBook Pro with a \SI{2.4}{\giga\hertz} 8-core Intel Core i9-9980HK CPU and \SI{64}{\giga\byte} of memory.

\subsection{Convergence}

\subsubsection{Laplace--Beltrami problem on a smooth surface}\label{sec:results_sphere}

\begin{figure}[htb]
	\centering
    \begin{minipage}{0.39\textwidth}
      \hspace{0.2cm}
      \raisebox{2em}{\includegraphics[width=0.76\textwidth]{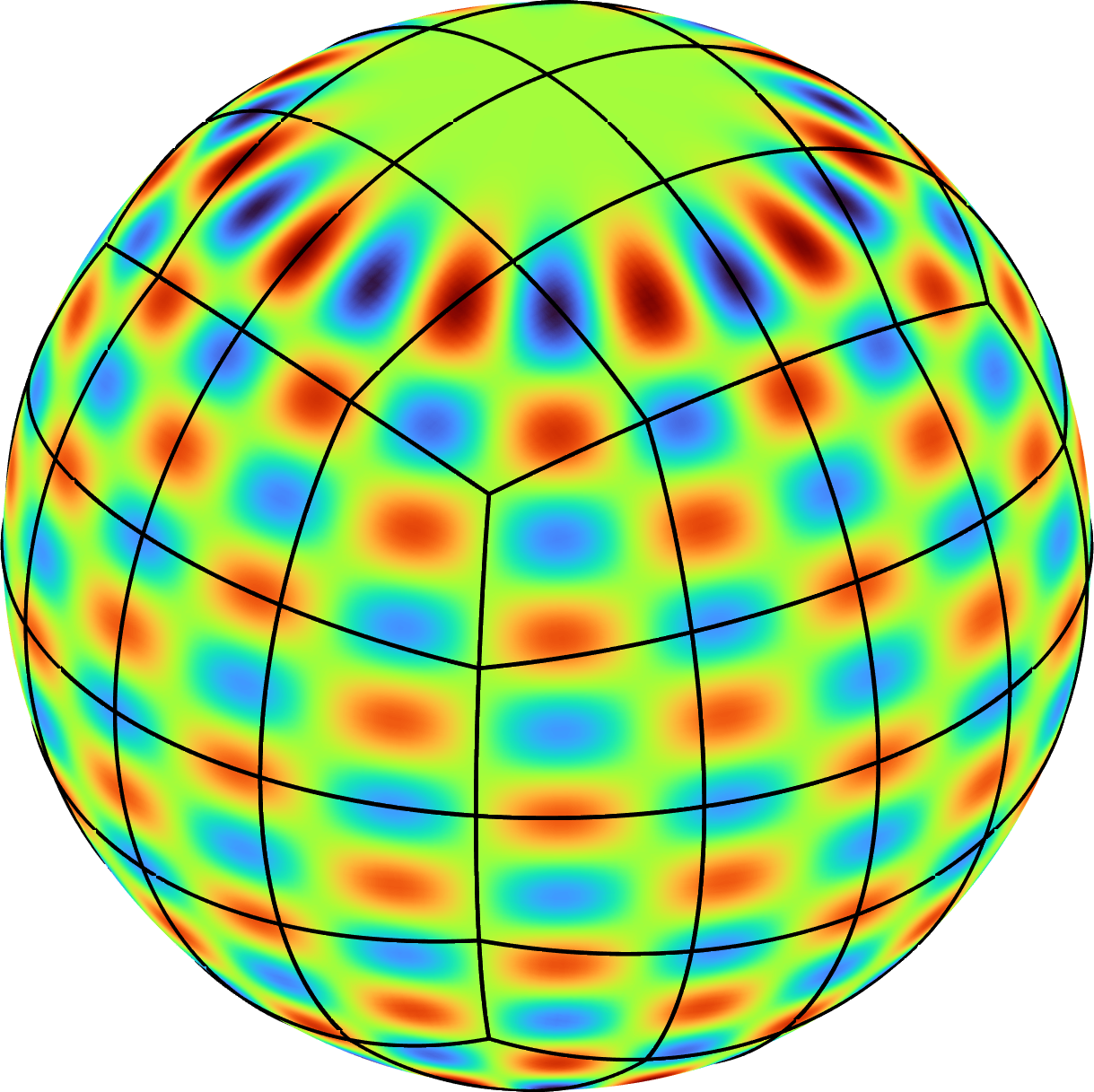}}
    \end{minipage}%
    \begin{minipage}{0.6\textwidth}
	    \centering
		\scalebox{0.97}{\input{conv_space.pdf_tex}}
    \end{minipage}%
    \vspace{-0.5em}
	\caption{(Left) A high-order cubed sphere mesh is used to approximate a spherical harmonic. (Right) The surface HPS scheme converges at a rate of $\mathcal{O}(h^{p-1})$.}
	\label{fig:convergence_sphere}
\end{figure}

We begin with a pure Laplace--Beltrami problem on a smooth surface,
\[
\lapbel u = f,
\]
where $u$ and $f$ are mean-zero functions over $\surf$, i.e., $\int_\surf f(\vec{x})\,d\vec{x} = 0$. To test the convergence of the method, we take $\surf$ to be the unit sphere; an exact solution to the Laplace--Beltrami problem is then given by the spherical harmonic $u(\vec{x}) = Y_\ell^m(\vec{x})$ when $f(\vec{x}) = -\ell (\ell+1) Y_\ell^m(\vec{x})$. \Cref{fig:convergence_sphere} (left) shows the $(\ell,m) = (20,10)$ spherical harmonic evaluated on a cubed sphere mesh. In \cref{fig:convergence_sphere} (right), we measure spatial convergence to the exact solution under mesh refinement for polynomial orders  $p = 4$, 8, 12, 16, and 20. As the mean mesh size $h \to 0$, the $L^\infty$ relative error decreases at a rate of $\mathcal{O}(h^{p-1})$, indicating high-order convergence.

\subsubsection{Laplace--Beltrami problem on surfaces with edges and corners}

In the present domain decomposition context, the interface conditions enforced between elements may still be enforced on a surface with sharp edges and corners. However, the binormal vector $\binorm$ along the interface between two elements may no longer be uniquely defined, as the elements on either side of a sharp edge will possess locally different binormal vectors. To remedy this, a natural choice is to enforce that the local fluxes on each element---defined via their respective local binormal vectors---must balance~\cite{Goodwill2021}.

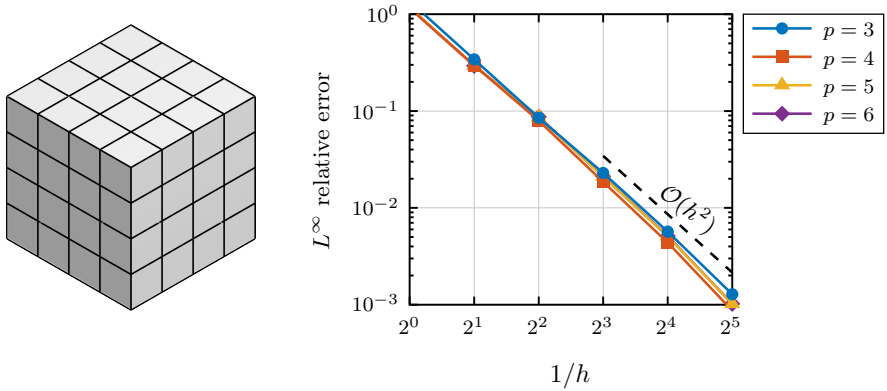
\begin{figure}[htb]
	\begin{minipage}{0.39\textwidth}
		\centering
		\raisebox{3em}{\includegraphics[width=0.65\textwidth]{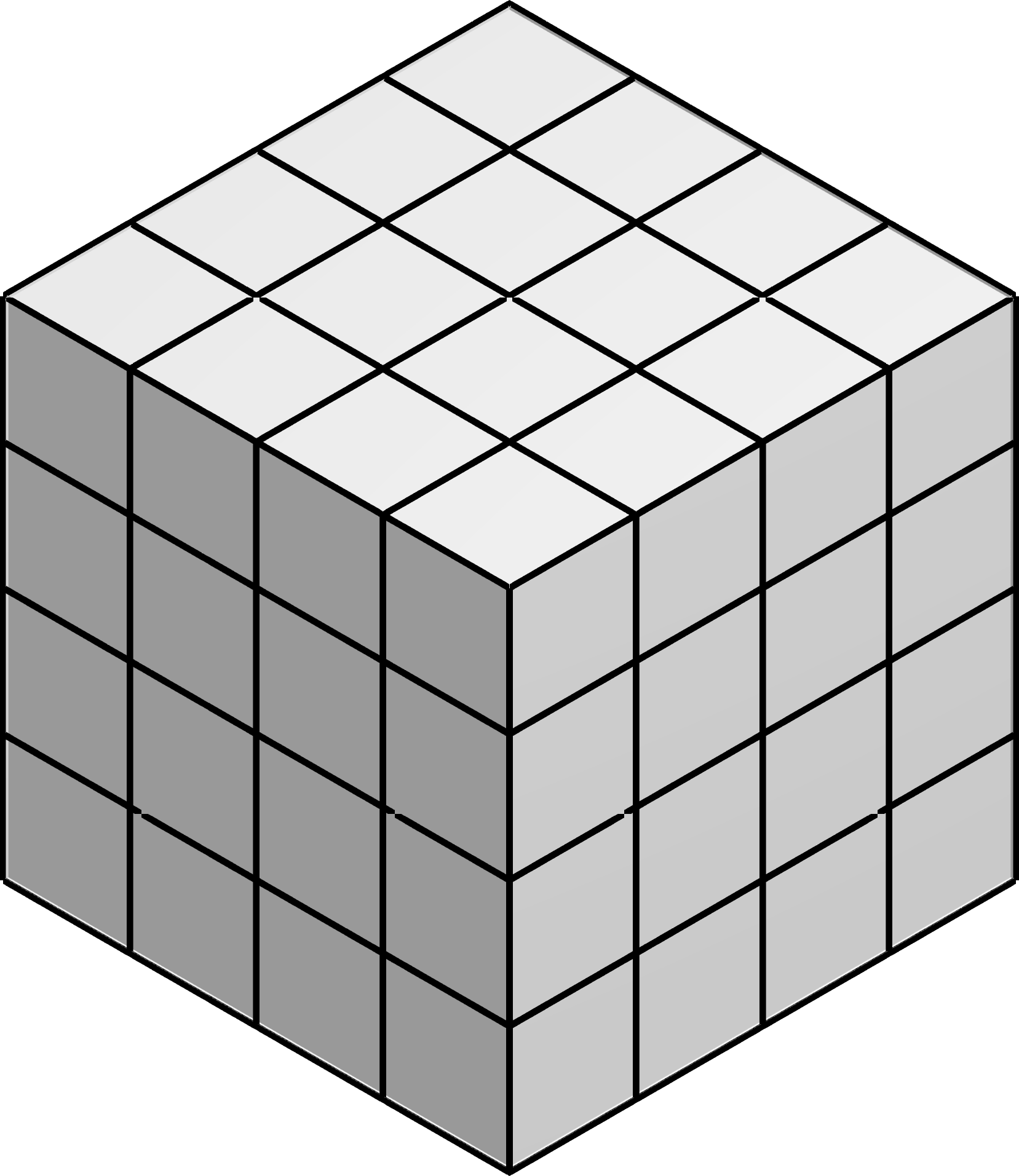}}
	\end{minipage}
	\begin{minipage}{0.6\textwidth}
		\begin{center}
		\scalebox{0.97}{\input{conv_cube.pdf_tex}}
		\begin{tikzpicture}[overlay]
			\put(11,88) {\rotatebox{-38}{\small$\mathcal{O}(h^2)$}};
		\end{tikzpicture}
		\end{center}
	\end{minipage}
	\vspace{-2em}
	\caption{Self-convergence study of the Laplace--Beltrami problem on the surface of a cube. High-order convergence is lost due to the weak vertex singularities induced by the corners of the cube. Polynomial orders $p = 3$, $4$, $5$, and $6$ all exhibit second-order convergence, as the solution $u$ is only in $H^2(\surf)$.}
	\label{fig:cube}
\end{figure}

Just as the solution to Laplace's equation in a domain with corners may possess weak corner singularities, solutions to the Laplace--Beltrami problem on a non-smooth surface may possess weak vertex singularities. Let $\surf$ be a non-smooth surface with a countable number of sharp vertices---that is, vertices whose conic angle does not equal $2 \pi$. Then in a polar neighborhood of vertex $i$, the solution $u$ to the Laplace--Beltrami problem $\lapbel u = f$ scales as $u \sim r^{2 \pi / \gamma_i}$, where $\gamma_i$ is the conic angle of the vertex~\cite[Lemma 1]{Goodwill2021}. \Cref{fig:cube} shows a self-convergence study of a Laplace--Beltrami problem on the surface of a cube, where $\gamma_i = 3\pi/2$ and $u \sim r^{4/3}$. Hence, $u \in H^2(\surf)$~\cite[Theorem 4]{Goodwill2021} and second-order convergence is expected~\cite{Brenner2008}. Using a uniform surface mesh, high-order convergence is lost in the presence of the vertex singularities regardless of the polynomial order employed. In this case, a refinement strategy based on $hp$-adaptivity may be used to recover super-algebraic convergence~\cite{Guo1986,Fortunato2021}.

\begin{figure}[htb]
	\begin{minipage}{0.39\textwidth}
		\raisebox{2em}{\includegraphics[width=0.85\textwidth]{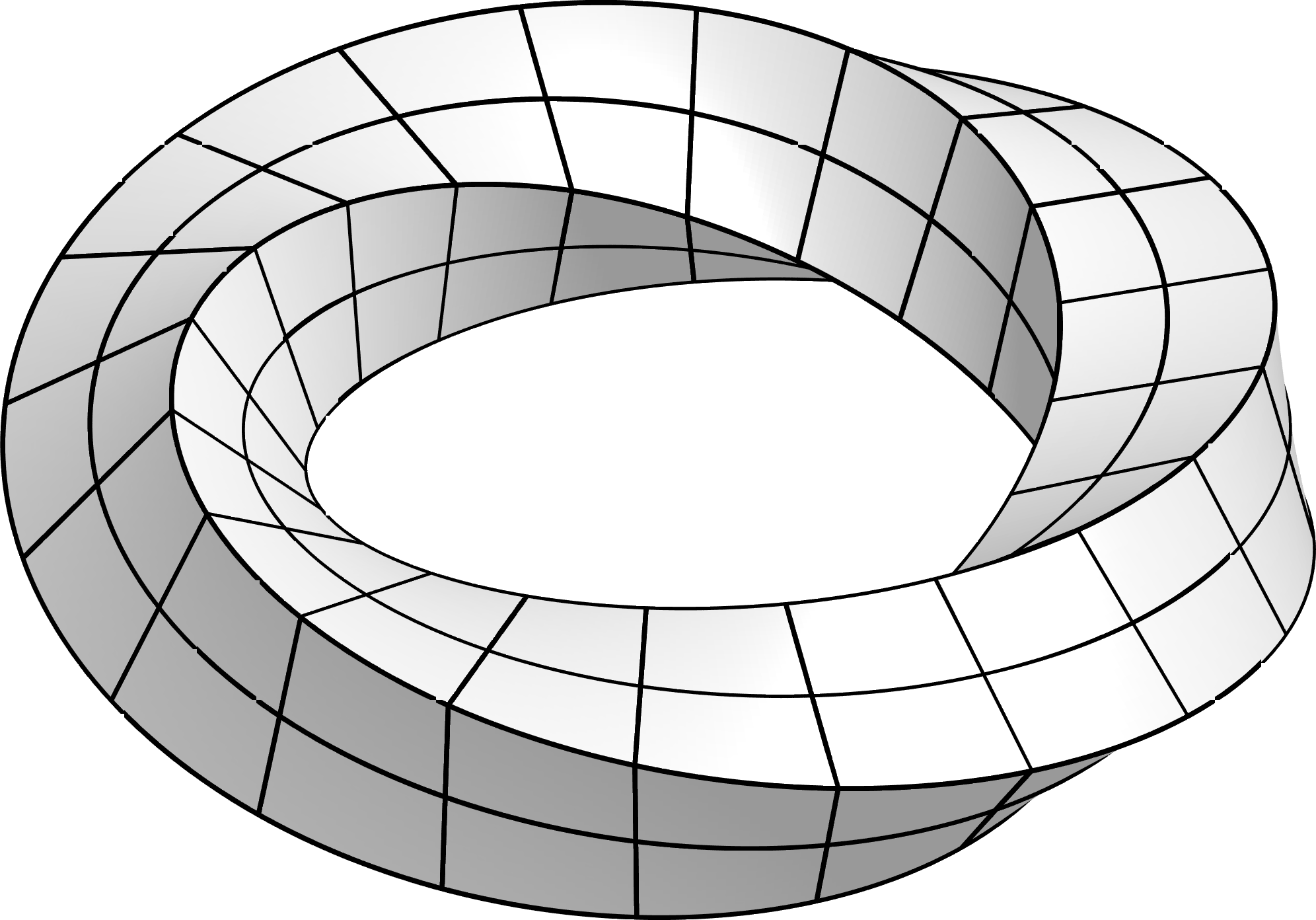}}
	\end{minipage}
	\begin{minipage}{0.6\textwidth}
		\begin{center}
		\scalebox{0.97}{\input{conv_twisted.pdf_tex}}
		\end{center}
	\end{minipage}
	\vspace{-1em}
	\caption{Self-convergence study of the Laplace--Beltrami problem on the surface of a twisted torus. High-order convergence can be obtained in the presence of sharp edges as long as no sharp corners are present.}
	\label{fig:twisted_torus}
\end{figure}

Note that in the case of a non-smooth surface $\surf$ with sharp edges but no sharp vertices, the solution to the Laplace--Beltrami problem no longer possesses geometrically-induced corner singularities and high-order convergence may still be obtained. \Cref{fig:twisted_torus} shows self-convergence of a Laplace--Beltrami problem on a twisted toroidal surface with sharp edges, where high-order convergence is recovered.

\subsection{Performance}

\begin{figure}[htb]
	\begin{center}
		\small
		\hspace{0.5cm}
		\input{complexity.pdf_tex}%
		\begin{tikzpicture}[overlay]
			\node[] at (-8.3,4.5) {\rotatebox{32}{$\mathcal{O}(\nelem)$}};
			\node[] at (-1.75,4.1) {\rotatebox{32}{$\mathcal{O}(\nelem)$}};
		\end{tikzpicture}
	\end{center}
	\vspace{-1em}
	\caption{Computational complexity of the surface HPS method. (Left) The runtime of the factorization stage is observed to scale as $\mathcal{O}(\nelem)$, even for a problem with $10^6$ degrees of freedom. For large enough problem sizes, the asymptotic scaling $\mathcal{O}(\nelem^{3/2})$ will eventually dominate. (Right) The total memory required to store all operators in the surface HPS method is also observed to scale linearly in this regime.}
	\label{fig:complexity}
\end{figure}
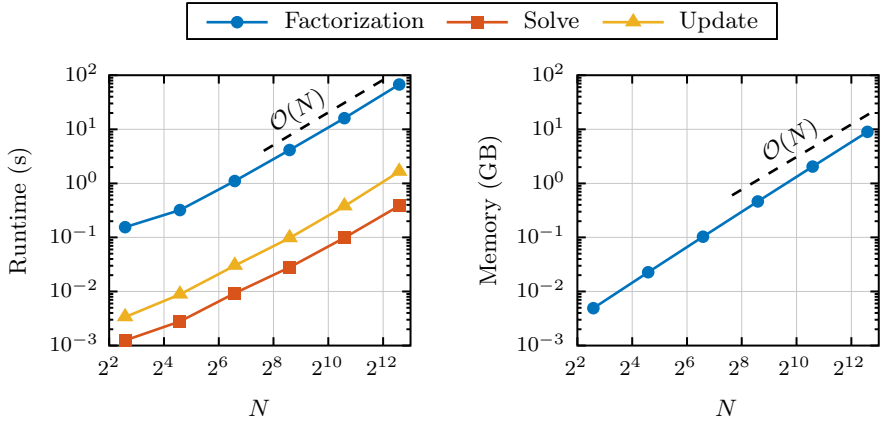

To illustrate the computational complexity of the surface HPS method, we measure the runtime and memory consumption of the factorization and solve stages as the number of elements in the mesh $\nelem$ is increased. \Cref{fig:complexity} (left) plots the runtime of each stage versus $\nelem$ for the simple Laplace--Beltrami problem considered in \cref{sec:results_sphere}, with $p=16$. Although the theoretical complexity of the factorization stage is $\mathcal{O}(\nelem^{3/2})$, here we observe $\mathcal{O}(\nelem)$ scaling even for a problem with over $10^6$ degrees of freedom. We expect that for a sufficiently fine mesh, the asymptotic complexity of $\mathcal{O}(\nelem^{3/2})$ will eventually dominate. Also plotted is the runtime complexity of updating the particular solutions given a new righthand side $f$. \Cref{fig:complexity} (right) depicts the memory complexity (i.e., total storage cost) of the method as $\nelem \to \infty$. Again, while $\mathcal{O}(\nelem \log \nelem)$ memory complexity is expected, in practice we observe $\mathcal{O}(\nelem)$ complexity for moderately-sized problems.

\begin{figure}[htb]
	\begin{minipage}{0.38\textwidth}
		\raisebox{1.2em}{\includegraphics[width=\textwidth]{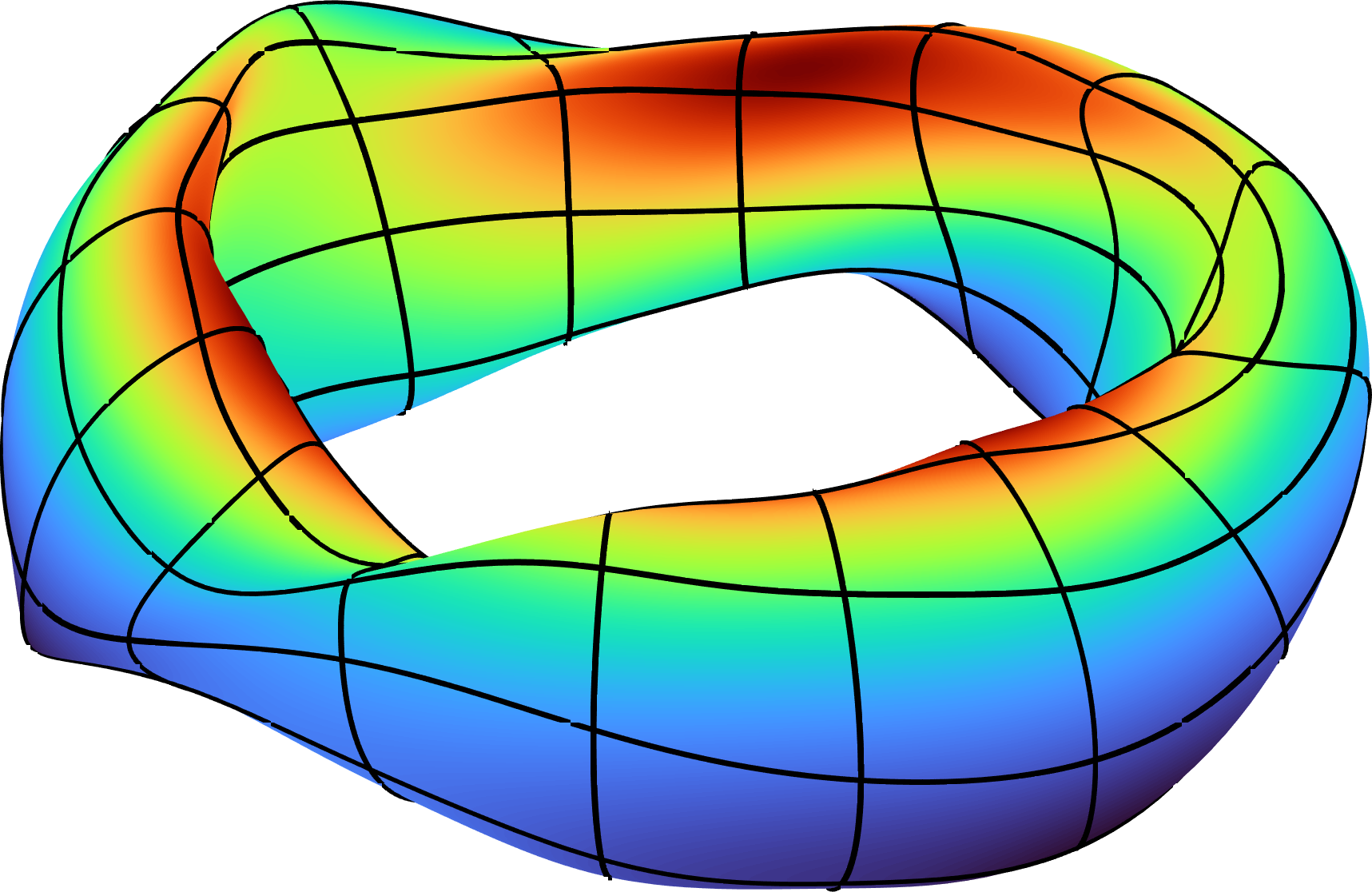}}
	\end{minipage}
	\hspace{-2.2cm}
	\rlap{%
	\begin{minipage}{0.9\textwidth}
		\begin{center}
		\scalebox{0.97}{\input{stellarator_benchmark.pdf_tex}}
		\end{center}
	\end{minipage}}
	\vspace{-0.3em}
	\caption{Relative accuracy versus runtime of the surface HPS compared to a high-order finite element method using MFEM~\cite{MFEM} and UMFPACK~\cite{UMFPACK}, an FFT-based spectral method~\cite{ImbertGerard2017}, and a layer potential approach~\cite{ONeil2018}. We use $p=20$ for both the surface HPS method and the finite element method. The solution to the benchmark problem is visualized on the left.}
	\label{fig:stellarator_benchmark}
\end{figure}

To better place the surface HPS method in the broader context of solvers for surface PDEs, we benchmarked our method against three alternative approaches for solving a reference Laplace--Beltrami problem $\lapbel u = f$ on a stellarator geometry~\cite{Garabedian2002, Malhotra2019a}. We manufacture a righthand side $f$ so that the solution $u$ is the Coulombic potential induced by a unit charge located at the off-surface point $(0,0,3)$. The solution is visualized in \cref{fig:stellarator_benchmark} (left). We compare the surface HPS approach to: a high-order finite element method discretized using the MFEM~\cite{MFEM} library and solved used UMFPACK~\cite{UMFPACK}; a pseudospectral method preconditioned by the flat Laplacian, which can be applied fast using the fast Fourier Transform (FFT) and is specifically designed for use on genus-one surfaces~\cite{ImbertGerard2017}; and a layer-potential preconditioner~\cite{ONeil2018} discretized using spectrally-accurate quadratures~\cite{Malhotra2019a}. For each method, we perform a sequence of mesh refinements and record both the runtime and the $L^\infty$ relative error between the computed and reference solutions. As the surface HPS method and the high-order finite element method are both based on unstructured meshes, we use the same sequence of mesh refinements for each, with $p=20$; for the spectral methods, we refine the number of grid points on the surface globally. \Cref{fig:stellarator_benchmark} plots the relative error versus the runtime for all the methods. While the FFT-based method is dominant for lower accuracies, the surface HPS method starts to narrow the gap around a relative accuracy of $10^{-5}$. However, the FFT-based approach is only applicable on toroidal surfaces of genus one, whereas the other approaches can be used on general surfaces. The high-order finite element method is competitive for lower orders, but when run at 20th order (shown here), it performs about 10 times slower than the surface HPS scheme for a given relative accuracy.

\begin{figure}[htb]
	\begin{center}
		\small
		\input{mfem_comparison.pdf_tex}
	\end{center}
	\vspace{-1em}
	\caption{Runtime comparison of the surface HPS solver and the UMFPACK~\cite{UMFPACK} solver applied to a high-order finite element discretization from the finite element package MFEM~\cite{MFEM}. The runtime of each method is plotted versus (left) number of elements $N$ under mesh refinement with a fixed polynomial order of $p=20$, and (right) polynomial order $p$ under $p$-refinement with a fixed number of elements $N=2048$.}
	\label{fig:mfem_comparison}
\end{figure}
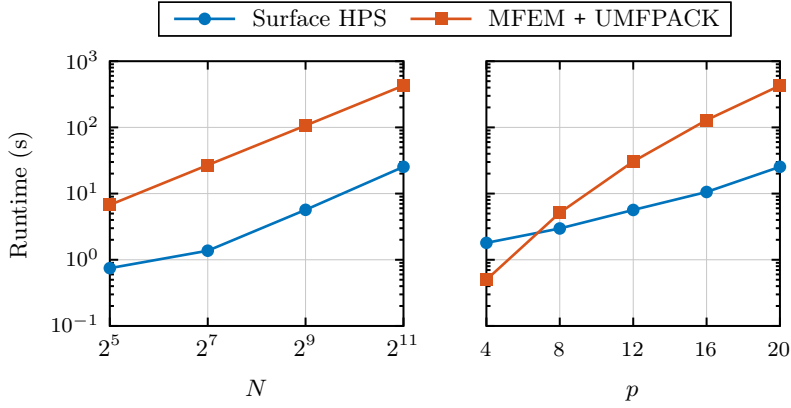

As both the surface HPS method and the high-order finite element method can use the same underlying high-order mesh, we further compare their runtimes using different polynomial degrees. \Cref{fig:mfem_comparison} shows the runtime of both methods plotted under mesh refinement (left) using a fixed polynomial order of $p=20$, and $p$-refinement (right) using a fixed mesh with $\nelem=2048$ elements. While the asymptotic complexity of UMFPACK matches that of the surface HPS scheme for a fixed $p$, the gap between the methods widens as $p \to \infty$. Finally, we note that while MFEM and UMFPACK are high-performance compiled codes that have been heavily optimized, our MATLAB implementation of the surface HPS scheme is relatively unoptimized.

\subsection{Hodge decomposition of a vector field}

The Laplace--Beltrami problem arises when computing the Hodge decomposition of tangential vector fields. For a vector field $\vec{f}$ tangent to the surface $\surf$, the Hodge decomposition~\cite{Schwarz1995} writes $\vec{f}$ as the sum of curl-free, divergence-free, and harmonic components,
\[
\vec{f} = \underbrace{\nabla_\surf u}_{\text{curl-free}} + \underbrace{\vec{n} \times \nabla_\surf v}_{\text{divergence-free}} + \underbrace{\vec{w}}_{\text{harmonic}}\hspace{-0.15cm},
\]
where $u$ and $v$ are scalar functions on $\surf$ and $\vec{w}$ is a harmonic vector field, i.e.,
\[
\nabla_\surf \cdot \vec{w} = 0, \qquad \nabla_\surf \cdot (\vec{n} \times \vec{w}) = 0.
\]
Such vector fields play an important role in integral-equation-based methods for computational electromagnetics~\cite{Epstein2015,Malhotra2019a}. To numerically compute such a decomposition, one may solve two Laplace--Beltrami problems for $u$ and $v$,
\[
\lapbel u = \nabla_\surf \cdot \vec{f} \quad \text{ and } \quad \lapbel v = -\nabla_\surf \cdot \left( \vec{n} \times \vec{f} \right),
\]
and then set $\vec{w} = \vec{f} - \nabla_\surf u - \vec{n} \times \nabla_\surf v$.

\begin{figure}[htb]
	\centering
	\begin{overpic}[width=0.95\textwidth, trim={0 3cm 0 0}, clip]{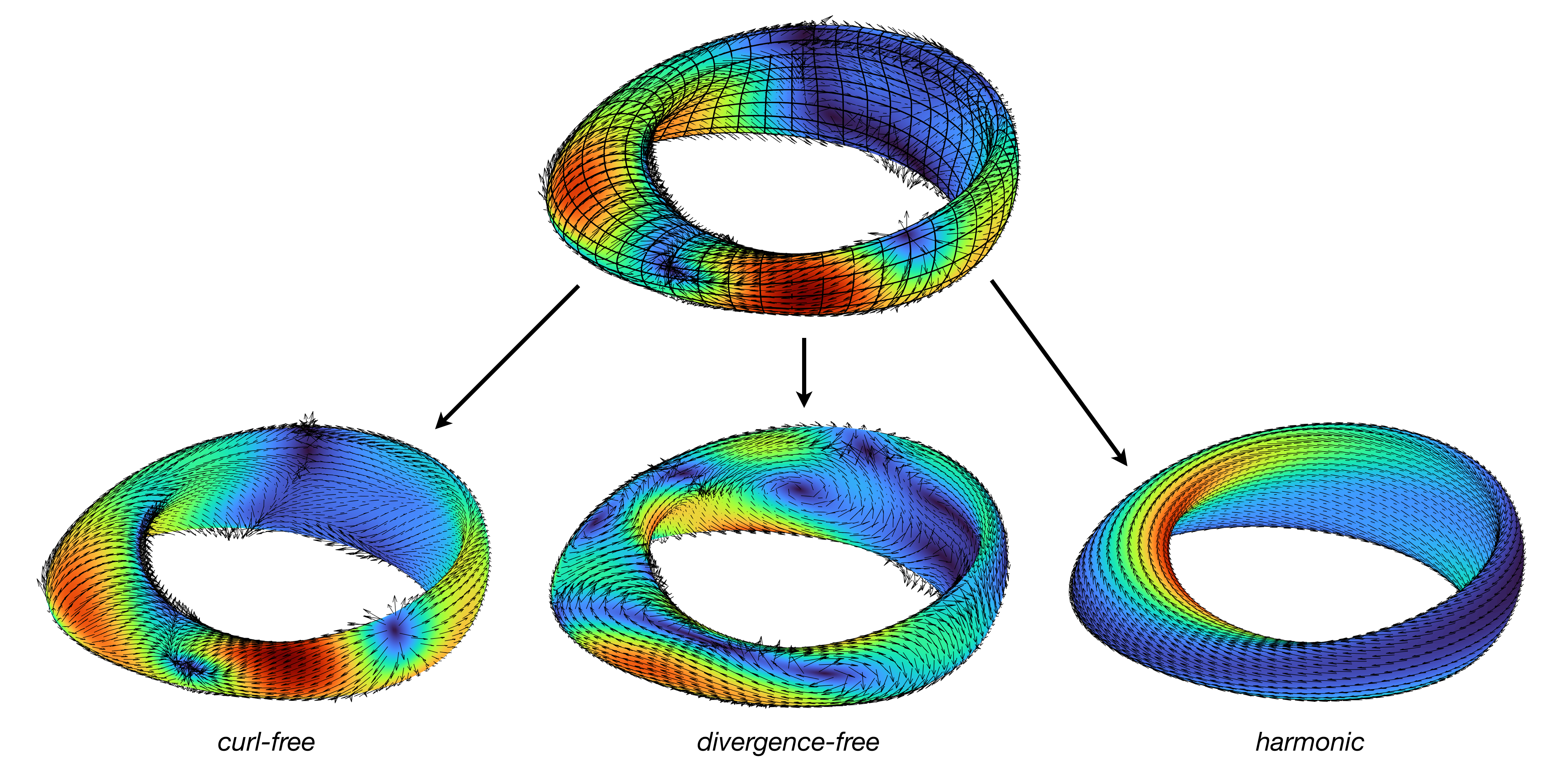}
		\put(12.5,-3) {\small\sffamily\textit{curl-free}}
		\put(42,-3) {\small\sffamily\textit{divergence-free}}
		\put(78.5,-3) {\small\sffamily\textit{harmonic}}
	\end{overpic}
	\vspace{1em}
	\caption{The Hodge decomposition of a random smooth vector field tangent to a deformed toroidal surface is numerically computed using a 16th-order mesh possessing roughly 200{,}000 degrees of freedom. The computation of the decomposition took 3.7 seconds. The resulting harmonic component $\vec{w}$ has residuals $\| \nabla_\surf \cdot \vec{w} \|_2 \approx 4 \cross 10^{-7}$ and $\| \nabla_\surf \cdot (\vec{n} \cross \vec{w}) \|_2 \approx 2 \times 10^{-6}$. Colors indicate magnitude.}
	\label{fig:hodge}
\end{figure}

\Cref{fig:hodge} shows the numerical Hodge decomposition of a random smooth vector field tangent to a deformed toroidal surface,computed on the 16th-order mesh depicted at the top, possessing roughly 200{,}000 degrees of freedom. The computation took 3.7 seconds. The residuals of the harmonic component are $\| \nabla_\surf \cdot \vec{w} \|_2 \approx 4 \cross 10^{-7}$ and $\| \nabla_\surf \cdot (\vec{n} \cross \vec{w}) \|_2 \approx 2 \times 10^{-6}$.

\subsection{Timestepping}\label{sec:timestepping}
Consider a nonlinear time-dependent PDE of the form
\[
\frac{\partial u}{\partial t} = \mathcal{L}_\surf u + \mathcal{N}(u),
\]
with $\mathcal{L}_\surf$ a linear surface differential operator and $\mathcal{N}$ a nonlinear, non-differential operator. Such PDEs appear as models for reaction--diffusion systems, where $\mathcal{L}_\surf$ contains the diffusion terms and $\mathcal{N}$ the reaction terms. As the timescales of diffusion and reaction are often orders of magnitude different, explicit timestepping schemes can suffer from a severe time step restriction. Implicit or semi-implicit timestepping schemes can alleviate stability issues, e.g., by treating the diffusion term $\mathcal{L}_\surf$ implicitly.

We use the implicit--explicit backward differentiation formula (IMEX-BDF) family of schemes~\cite{Ascher1995}, which are based on the backward differentiation formulae for the implicit term and the Adams--Bashforth for the explicit term. Fix a time step $\dt > 0$ and let $u^k(\vec{x}) \approx u(\vec{x}, k \dt)$ denote the approximate solution at time step $k$. Discretizing in time with the $K$th-order IMEX-BDF scheme results in a steady-state problem at each time step of the form
\begin{equation}\label{eq:imex_bdf}
\left( I - \omega \dt \mathcal{L}_\surf \right) u^{k+1} = \sum_{i=0}^{K-1} \mu_i u^{k-i} + \dt \sum_{i=0}^{K-1} \nu_i\,\mathcal{N}(u^{k-i}),
\end{equation}
where $\omega$, $\mu_i$, and $\nu_i$ are constants given by
\[
\begin{aligned}
&\text{IMEX-BDF1:} \; &&\omega = 1, &&\mu = 1, &&\nu = 1, \\
&\text{IMEX-BDF2:} \; &&\omega = \tfrac{2}{3}, &&\mu = \left(\tfrac43, -\tfrac13\right), &&\nu = \left(\tfrac43, -\tfrac23\right), \\
&\text{IMEX-BDF3:} \; &&\omega = \tfrac{6}{11}, &&\mu = \left(\tfrac{18}{11}, -\tfrac{9}{11}, \tfrac{2}{11}\right), &&\nu = \left(\tfrac{18}{11}, -\tfrac{18}{11}, \tfrac{6}{11}\right), \\
&\text{IMEX-BDF4:} \; &&\omega = \tfrac{12}{25}, &&\mu = \left(\tfrac{48}{25}, -\tfrac{36}{25}, \tfrac{16}{25}, -\tfrac{3}{25}\right), &&\nu = \left(\tfrac{48}{25}, -\tfrac{72}{25}, \tfrac{48}{25}, -\tfrac{12}{25}\right).
\end{aligned}
\]
For a fixed surface $\surf$, fixed time step $\dt$, and fixed linear operator $\mathcal{L}_\surf$, the inverse $(I - \omega \dt \mathcal{L}_\surf)^{-1}$ can be precomputed once for all time steps and applied fast using the HPS scheme, yielding a cost per time of $\mathcal{O}(\nelem \log \nelem)$.

\begin{figure}[htb]
  \begin{center}
    \small
    \hspace{1cm}
    \input{conv_time.pdf_tex}
    \begin{tikzpicture}[overlay]
      \node[] at (-3.7,5.28) {\rotatebox{-9}{\footnotesize$\mathcal{O}(\Delta t)$}};
      \node[] at (-3.7,4.28) {\rotatebox{-20}{\footnotesize$\mathcal{O}(\Delta t^2)$}};
      \node[] at (-3.7,3.45) {\rotatebox{-29}{\footnotesize$\mathcal{O}(\Delta t^3)$}};
      \node[] at (-3.7,2.65) {\rotatebox{-36}{\footnotesize$\mathcal{O}(\Delta t^4)$}};
    \end{tikzpicture}
  \end{center}
  \vspace{-1em}
  \caption{Self-convergence study of the reaction--diffusion system \cref{eq:cgle} on the cow model over a range of time steps. A reference solution is computed using a time step of $\Delta t = 2^{-10}$ to a final time of $t=1$.}
  \label{fig:t_convergence}
\end{figure}
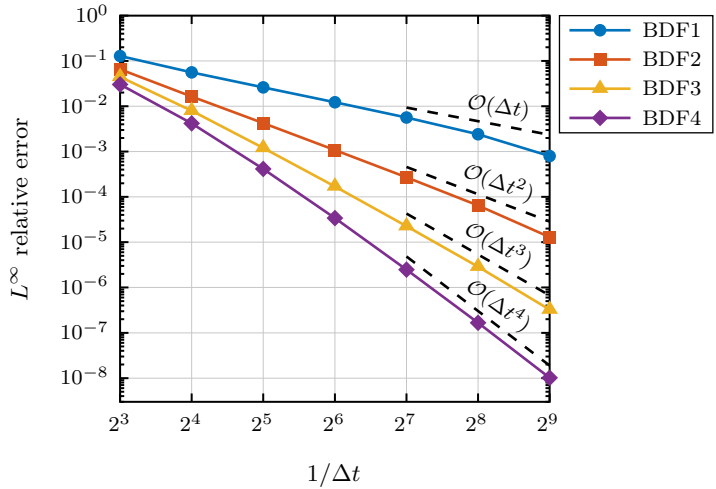

We now apply the IMEX-BDF methods to simulate time-dependent reaction--diffusion on surfaces. We choose three surface meshes of varying complexity: a blob shape defined through deformation of a sphere, a stellarator shape~\cite{Garabedian2002} with localized regions of high curvature, and a cow shape defined via a subdivision surface~\cite{CraneModel} and meshed to high order using Rhinoceros~\cite{Rhino}. \Cref{fig:bruss} depicts these meshes in the left column. On all three meshes, we use 16th-order elements to represent the solution. We perform a temporal self-convergence study on the cow mesh for the reaction--diffusion problem described in \cref{sec:cgle}. Using a reference solution computed with a time step of $\Delta t = 2^{-10}$, we simulate to a final time of $t=1$ over a range of time steps. \Cref{fig:t_convergence} depicts the results, demonstrating that $k$th-order accuracy is achieved for the IMEX-BDF$k$ method.

In all cases, the overall cost per time step is 0.05 seconds for the blob mesh with 24{,}576 degrees of freedom (\SI{0.04}{\second} for updating the particular solution and \SI{0.01}{\second} for the solve), 0.1 seconds for the stellarator mesh with 49{,}152 degrees of freedom (\SI{0.07}{\second} for updating the particular solution and \SI{0.03}{\second} for the solve), and 0.18 seconds for the cow mesh with 86{,}784 degrees of freedom (\SI{0.13}{\second} for updating the particular solution and \SI{0.05}{\second} for the solve).

\subsubsection{Turing system}\label{sec:bruss}

\begin{figure}
	\centering
	\begin{tabular}{ccc}
	~~\includegraphics[width=0.22\textwidth]{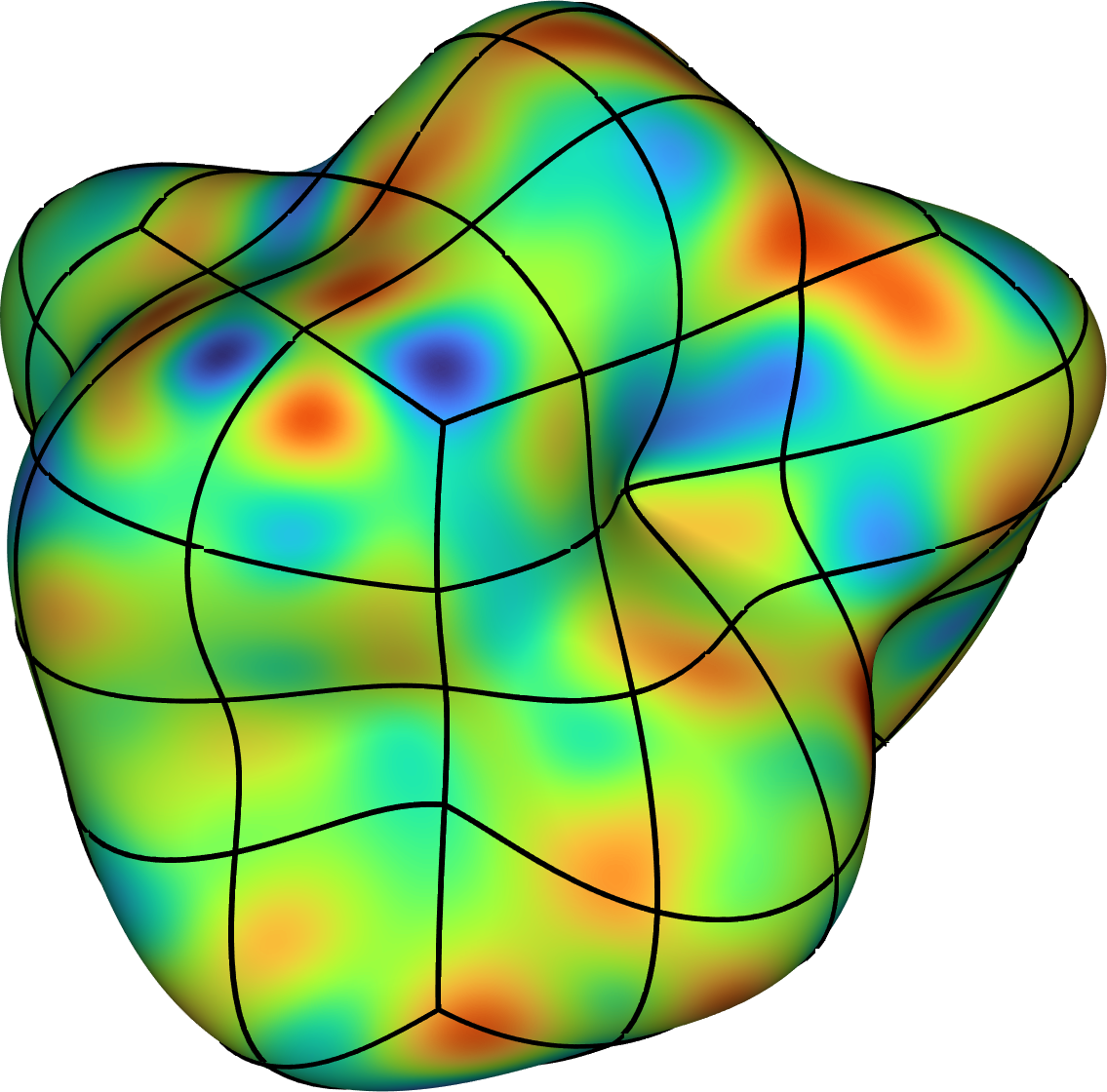} &
	~\includegraphics[width=0.22\textwidth]{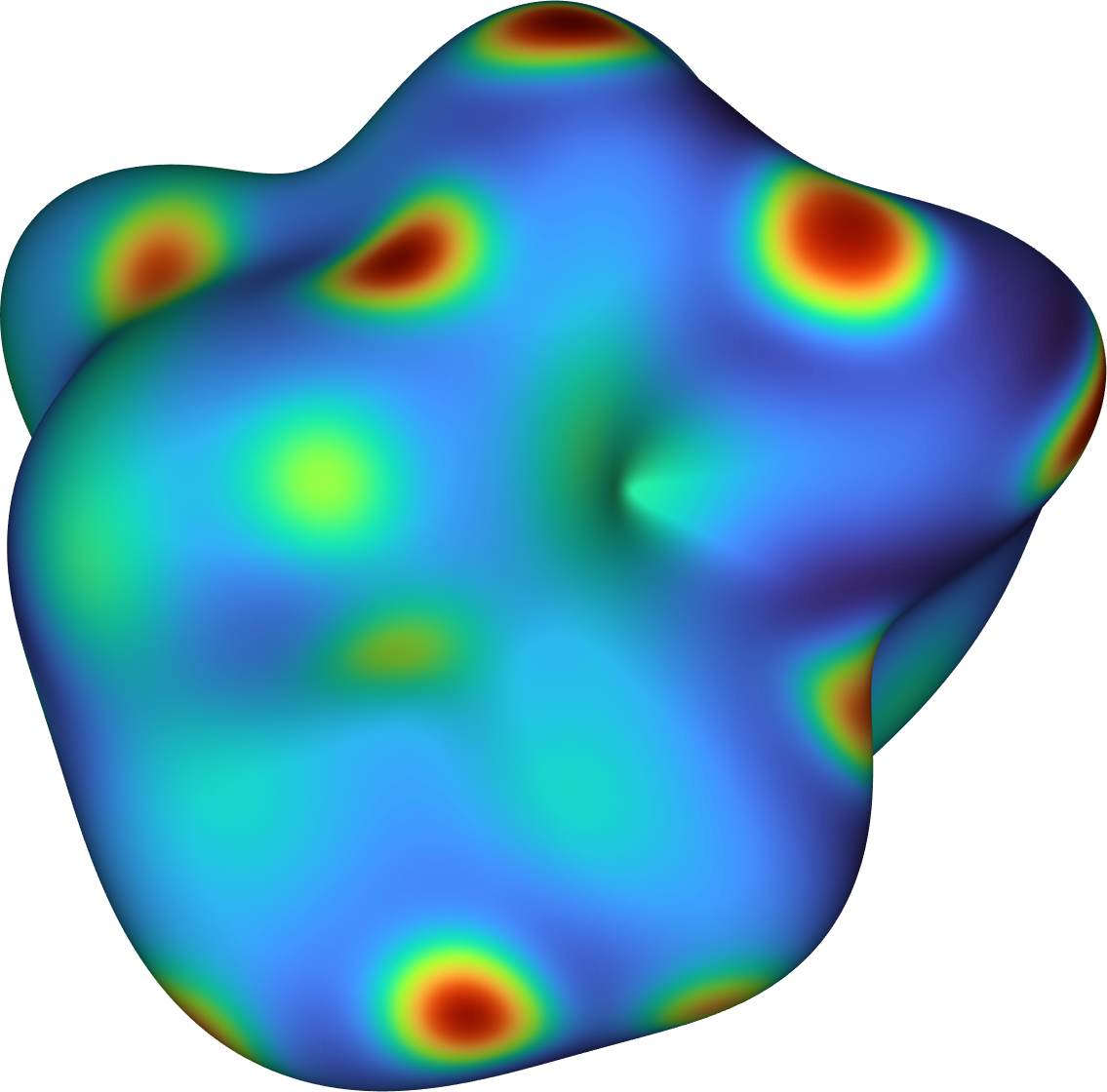} &
	~~\includegraphics[width=0.22\textwidth]{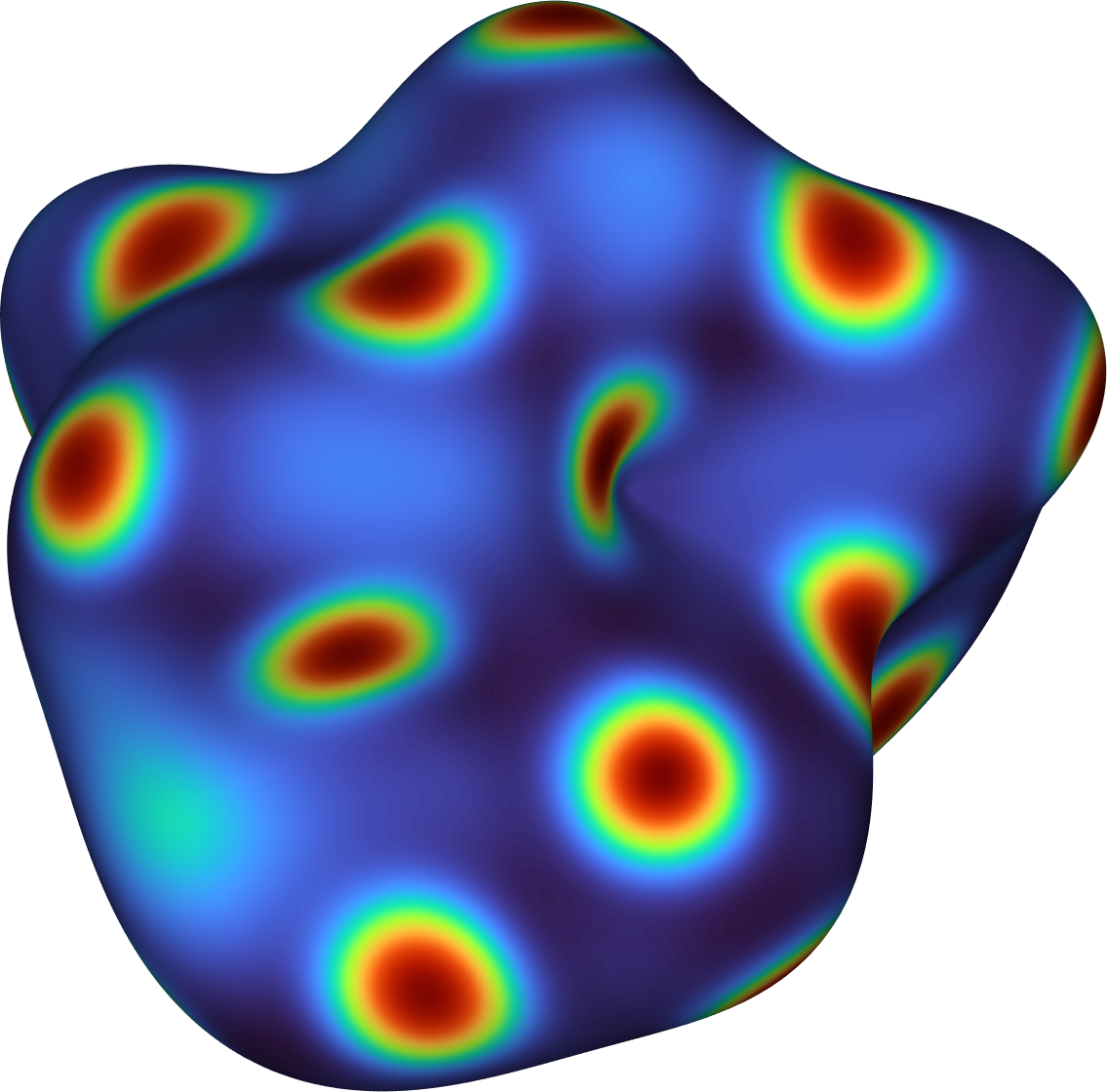} \\[1.2em]
	\includegraphics[width=0.3\textwidth]{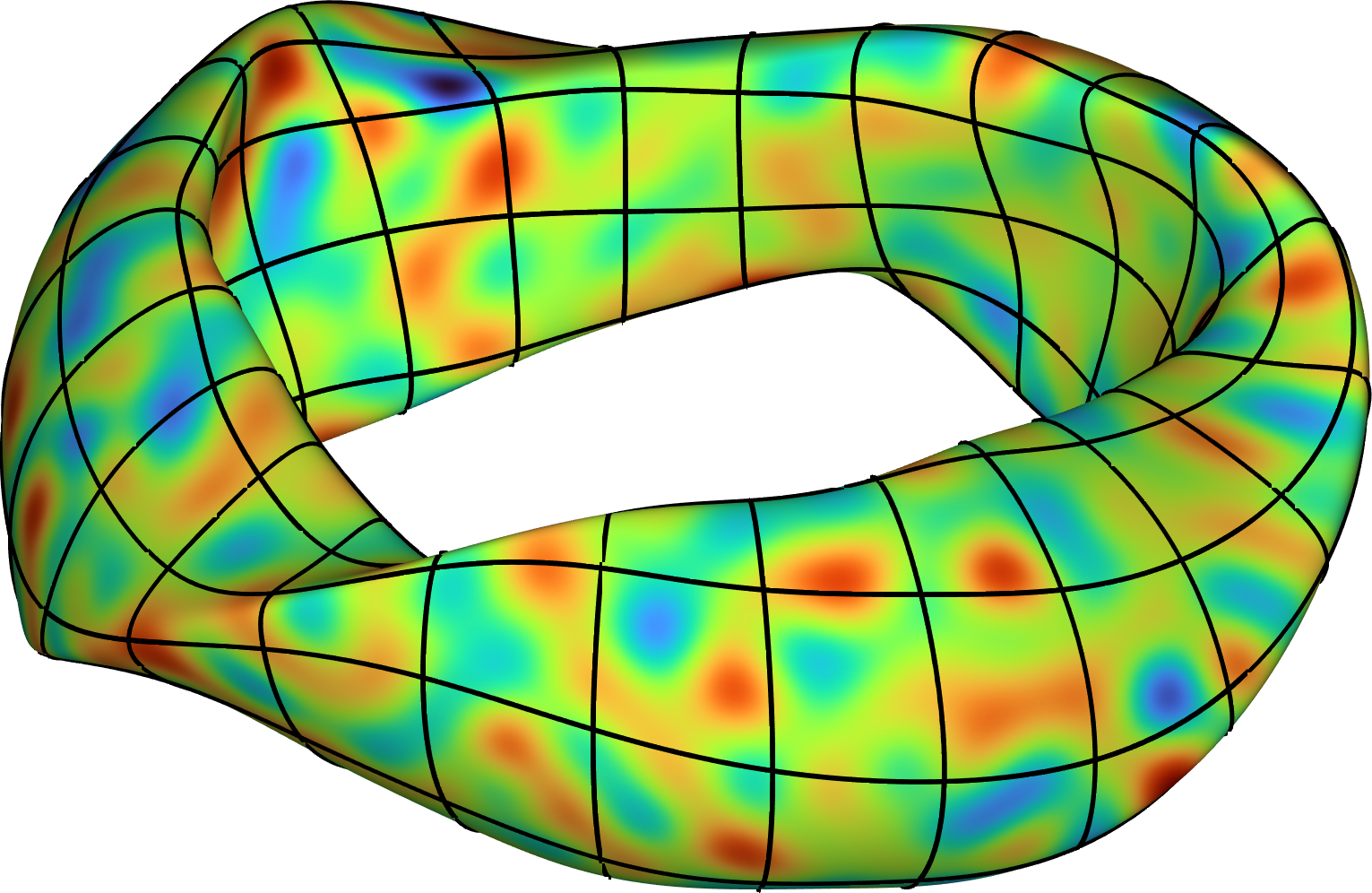} &
	\includegraphics[width=0.3\textwidth]{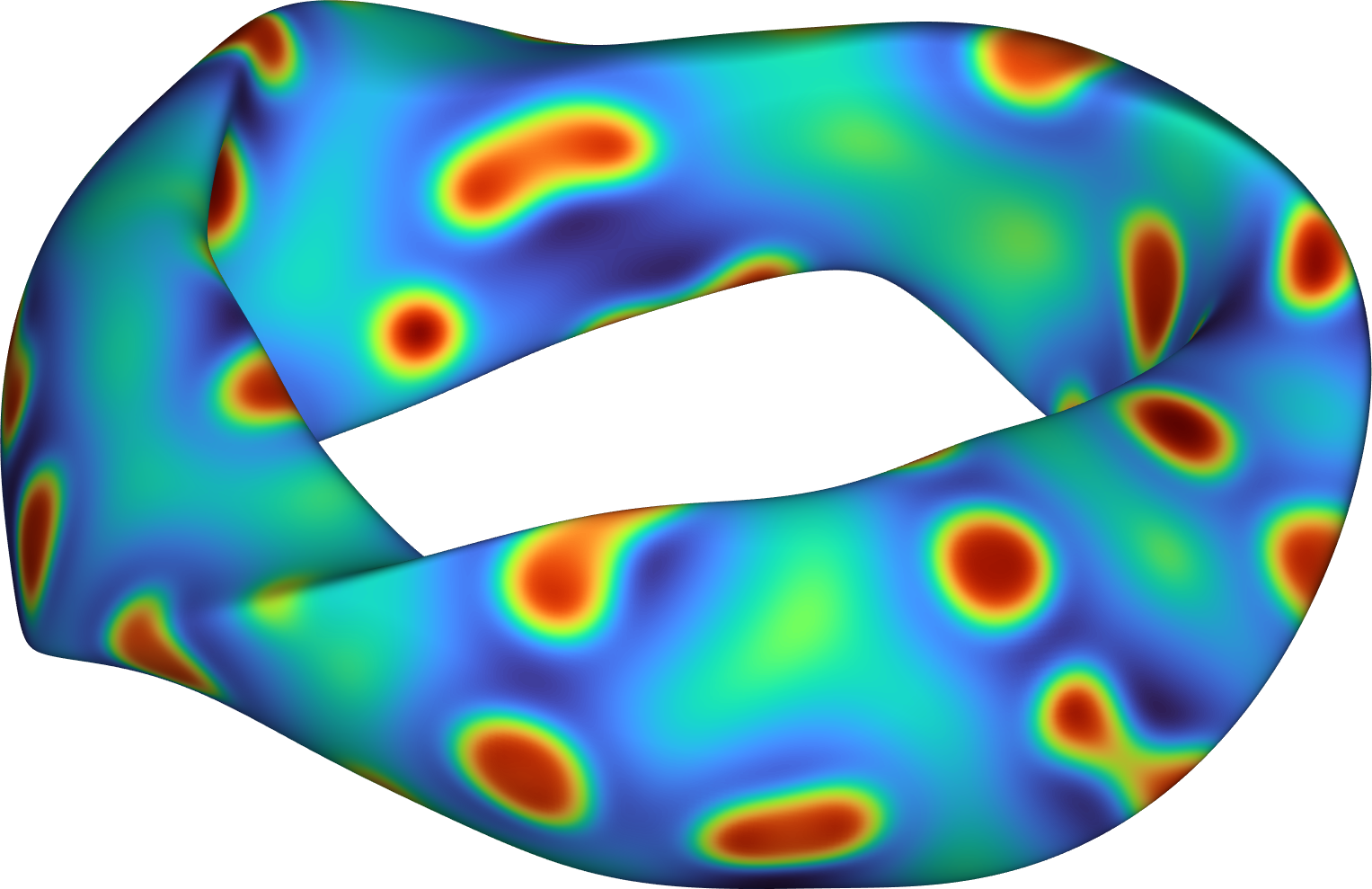} &
	\includegraphics[width=0.3\textwidth]{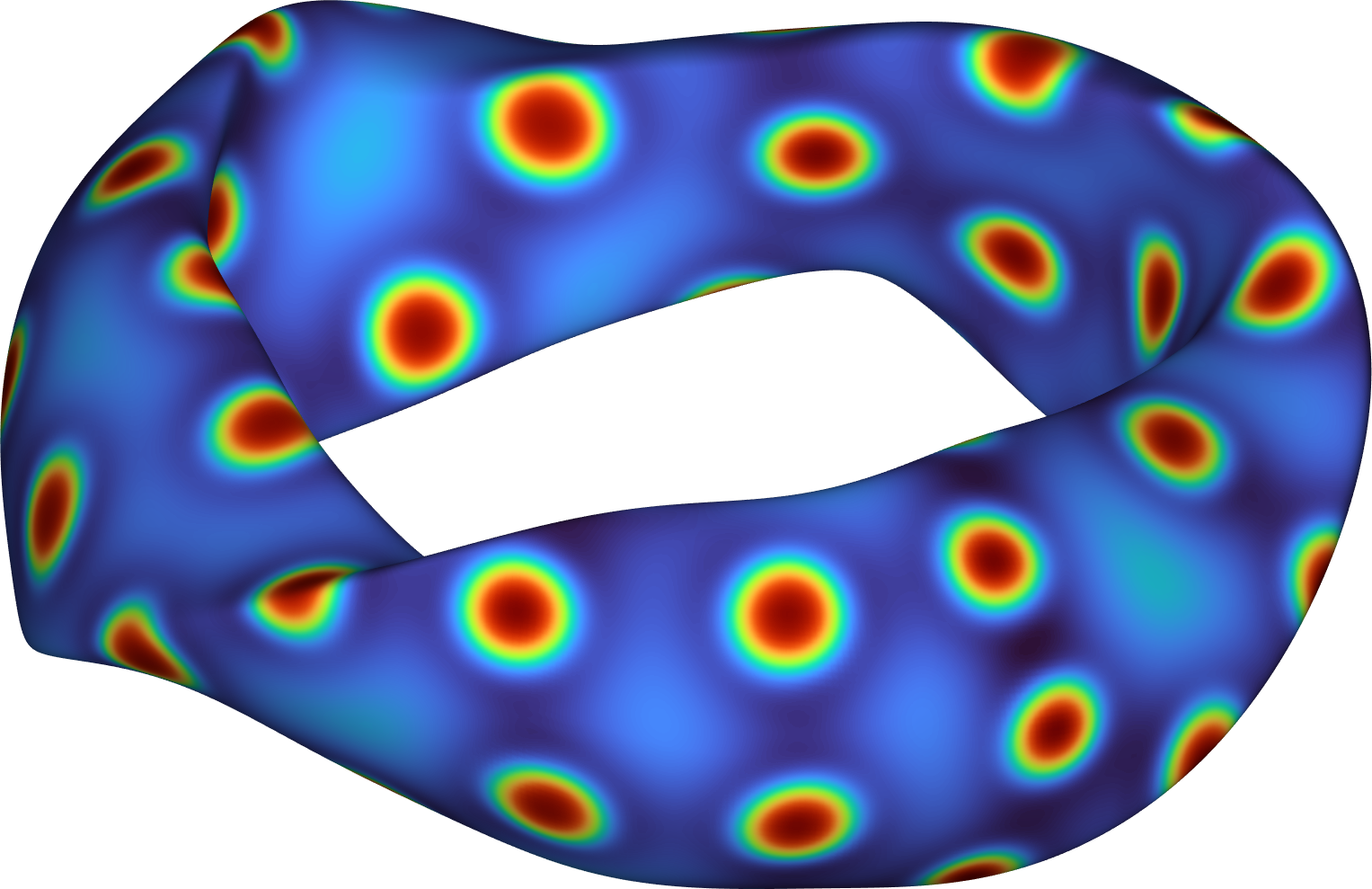} \\[1em]
	\includegraphics[width=0.29\textwidth]{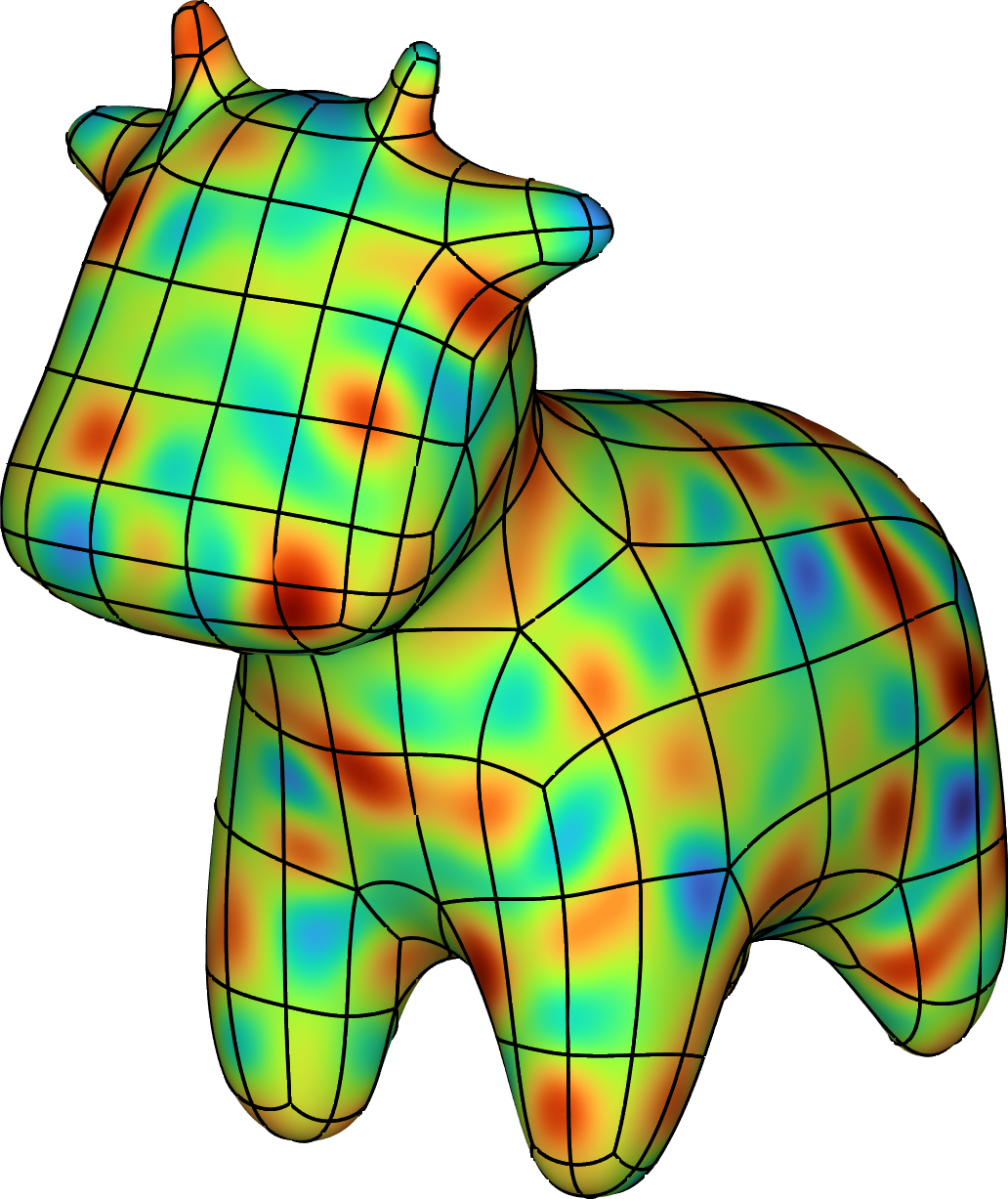} &
	\includegraphics[width=0.29\textwidth]{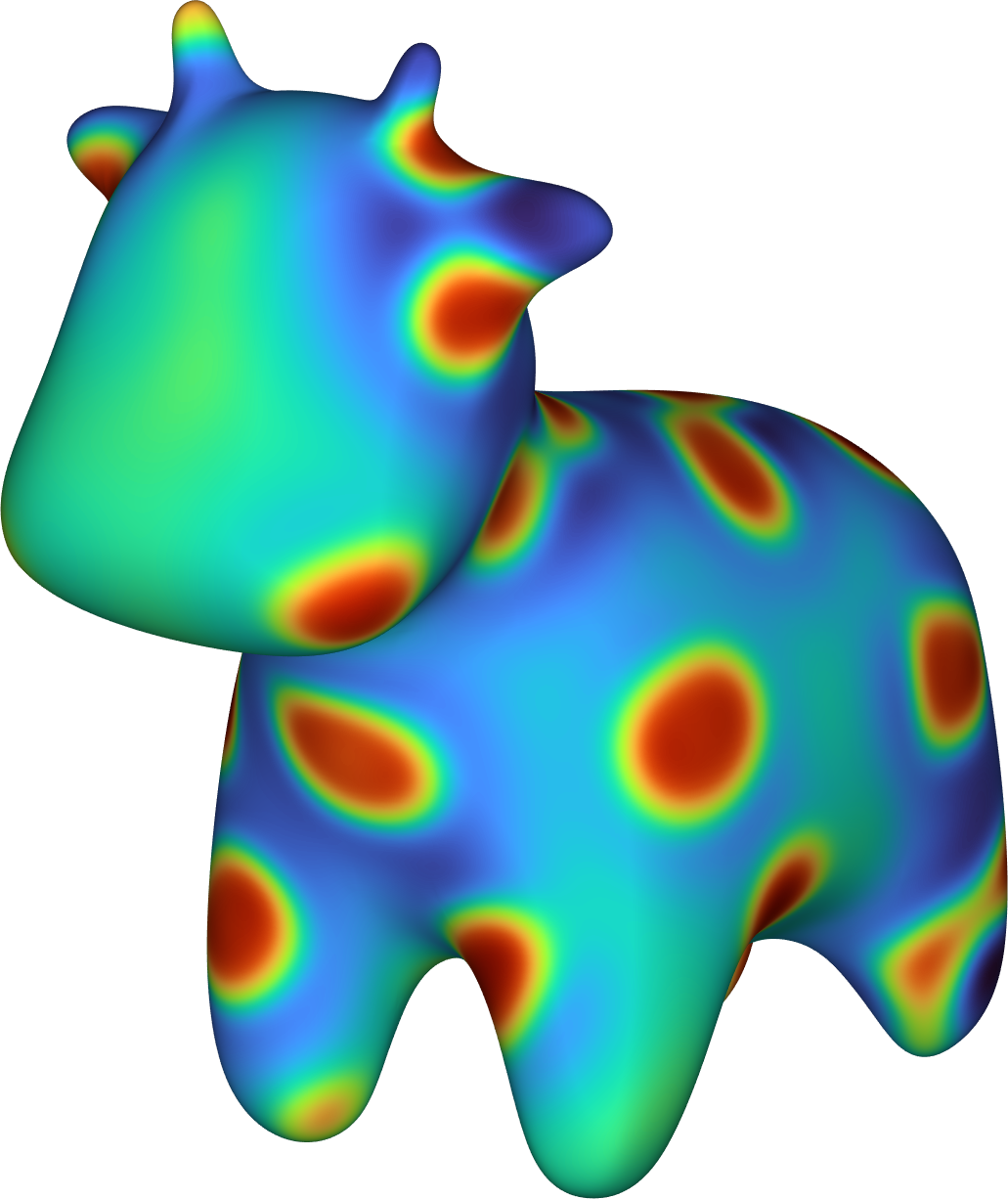} &
	\includegraphics[width=0.29\textwidth]{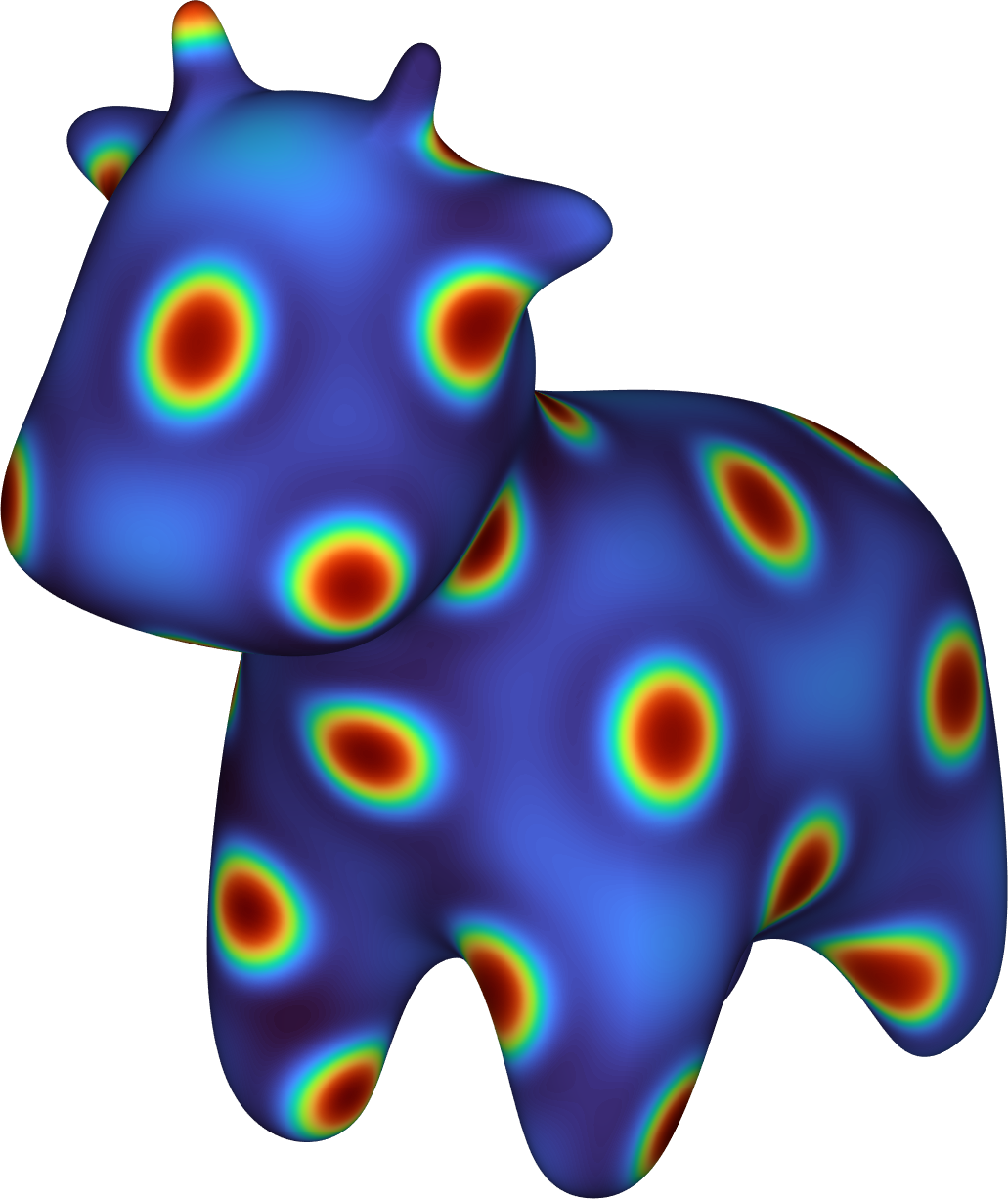}
	\end{tabular}
	\caption{Snapshots of a two-species Turing system simulated on three surfaces at times $t=0$, $t=20$, and $t=200$. We discretize in time with a fourth-order IMEX-BDF4 scheme and precompute a fast direct solver for the implicit problem, which computes the solution at each time step in $\mathcal{O}(N \log N)$ operations.}
	\label{fig:bruss}
\end{figure}

We first consider a two-species reaction--diffusion system on a surface $\surf$ given by~\cite{Barrio1999},
\[
\begin{aligned}
\frac{\partial u}{\partial t} &= \delta_u \lapbel u + \alpha u \left(1-\tau_1 v^2\right) + v \left(1-\tau_2 u\right), \\
\frac{\partial v}{\partial t} &= \delta_v \lapbel v + \beta v \left(1 + \frac{\alpha \tau_1}{\beta} u v \right) + u \left(\gamma + \tau_2 v\right).
\end{aligned}
\]
Solutions $u$ and $v$ to this system can exhibit Turing patterns---namely, spots and stripes---depending on the choice of parameters $\delta_u$, $\delta_v$, $\alpha$, $\beta$, $\gamma$, $\tau_1$, and $\tau_2$. Here, we take $\alpha = 0.899$, $\beta = -0.91$, $\gamma = -0.899$, $\tau_1 = 0.02$, and $\delta_u = 0.516 \delta_v$. On the blob, we take $\delta_v = 0.005$ and $\tau_2 = 0.15$. On the stellarator, we take $\delta_v = 0.02$ and $\tau_2 = 0.2$. On the cow, we take $\delta_v = 0.001$ and $\tau_2 = 0.2$. We simulate this system on all three surface meshes for 2000 time steps until a final time of $t = 200$ using a time step size of $\dt = 0.1$. Snapshots of the solutions at times $t = 0$, $t = 20$, and $t = 200$ are shown in \cref{fig:bruss}.

\subsubsection{Complex Ginzburg--Landau}\label{sec:cgle}

\begin{figure}
	\centering
	\begin{tabular}{ccc}
	~~\includegraphics[width=0.22\textwidth]{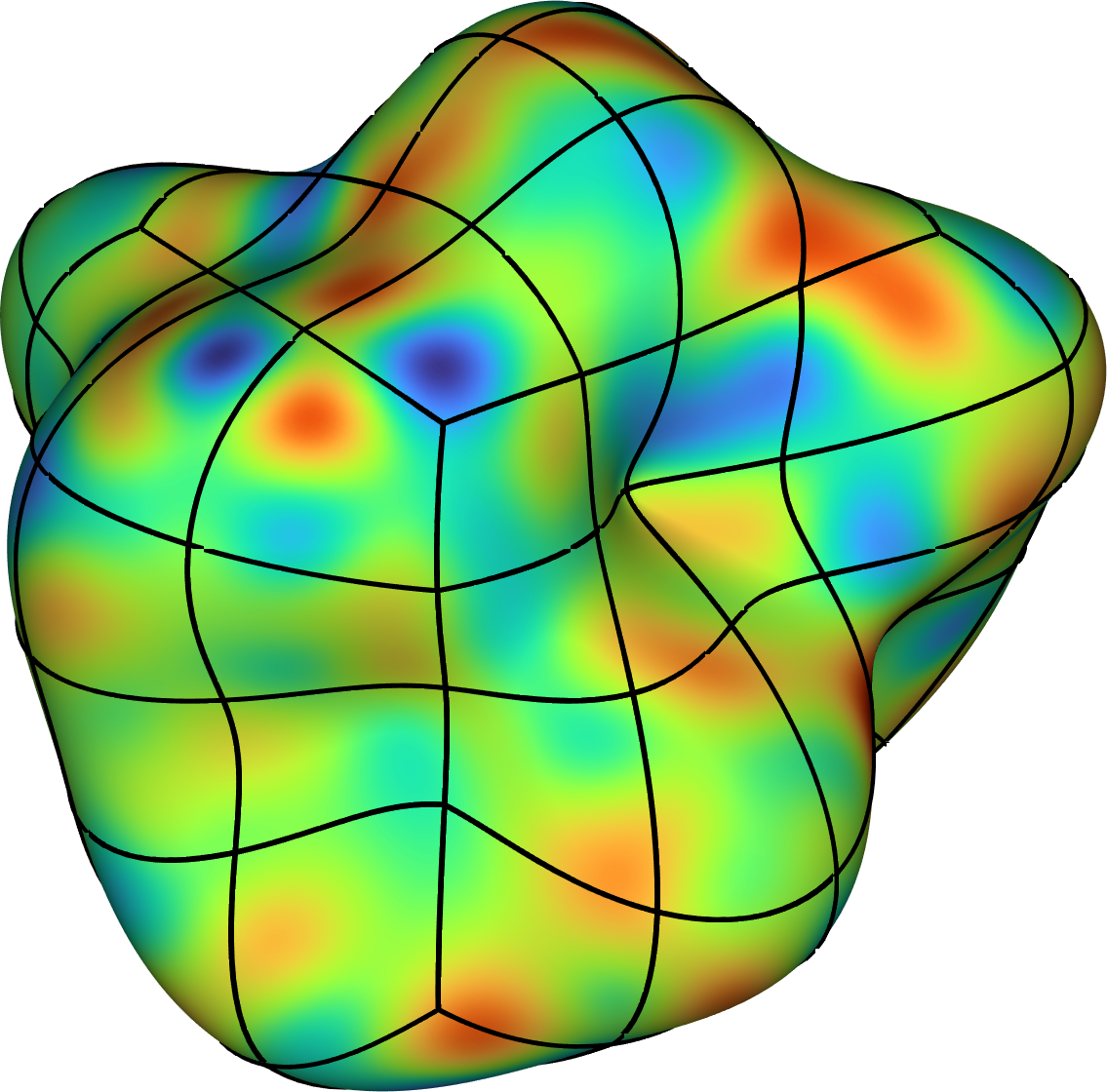} &
	~\includegraphics[width=0.22\textwidth]{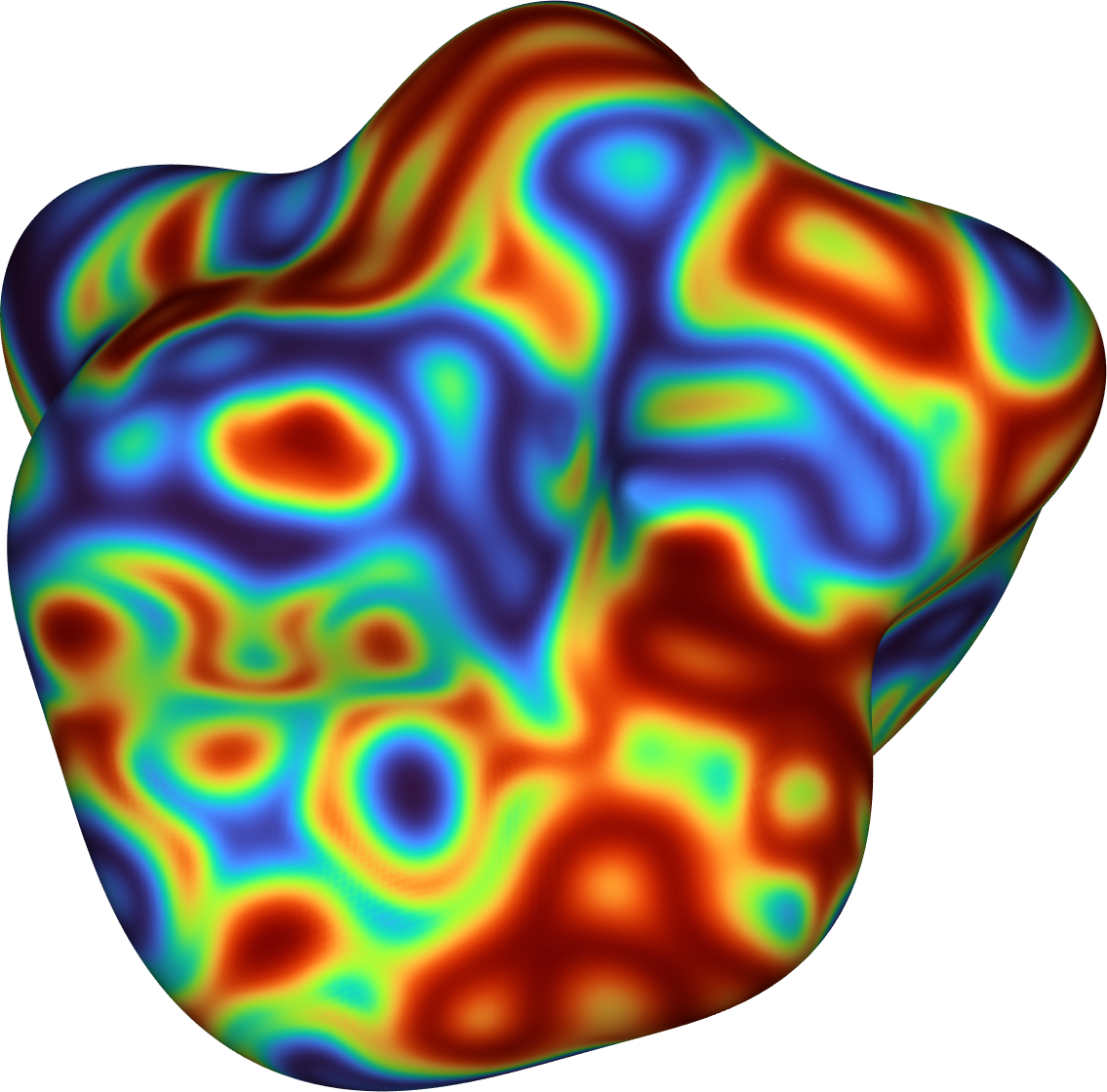} &
	~~\includegraphics[width=0.22\textwidth]{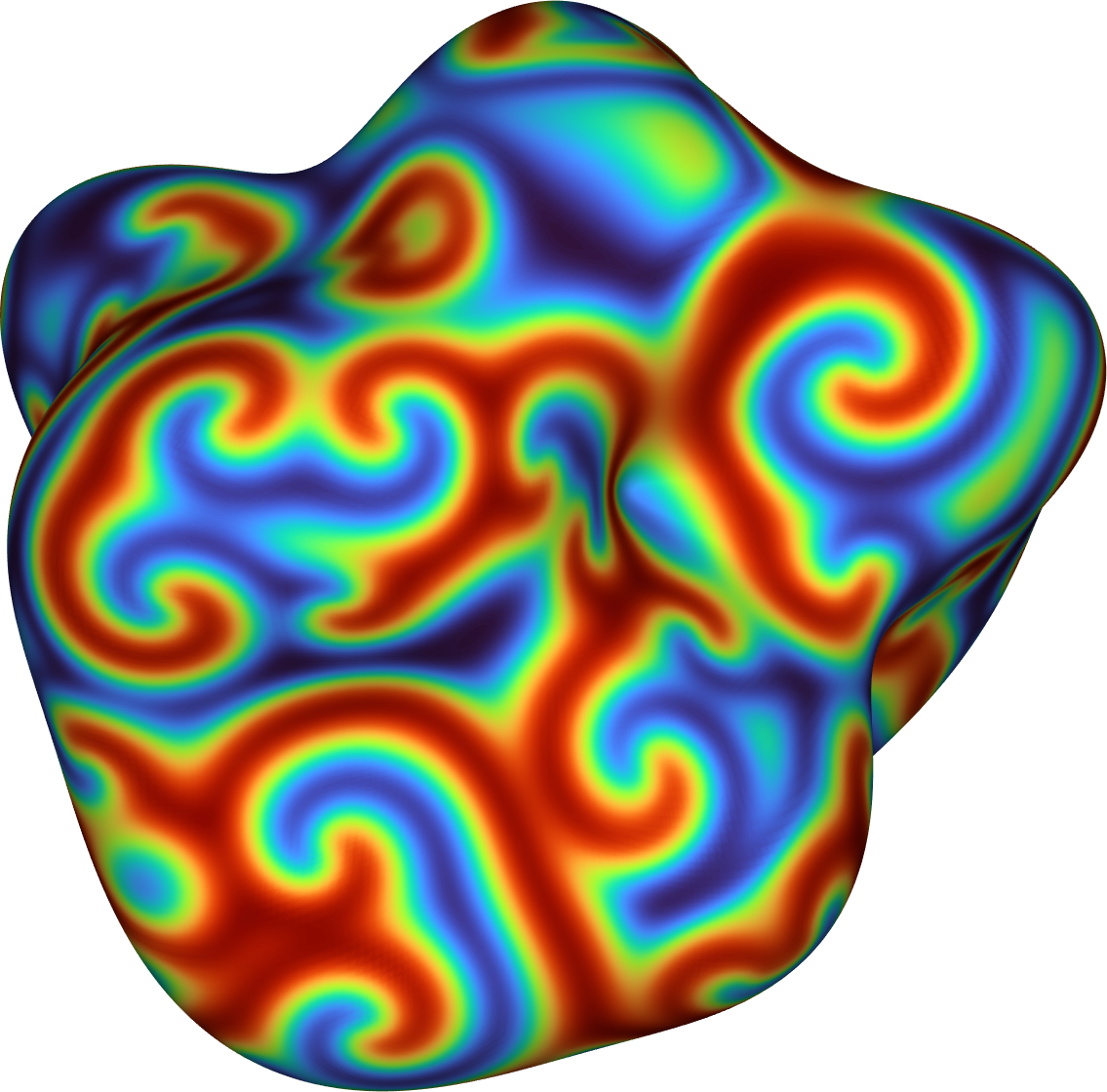} \\[1.2em]
	\includegraphics[width=0.3\textwidth]{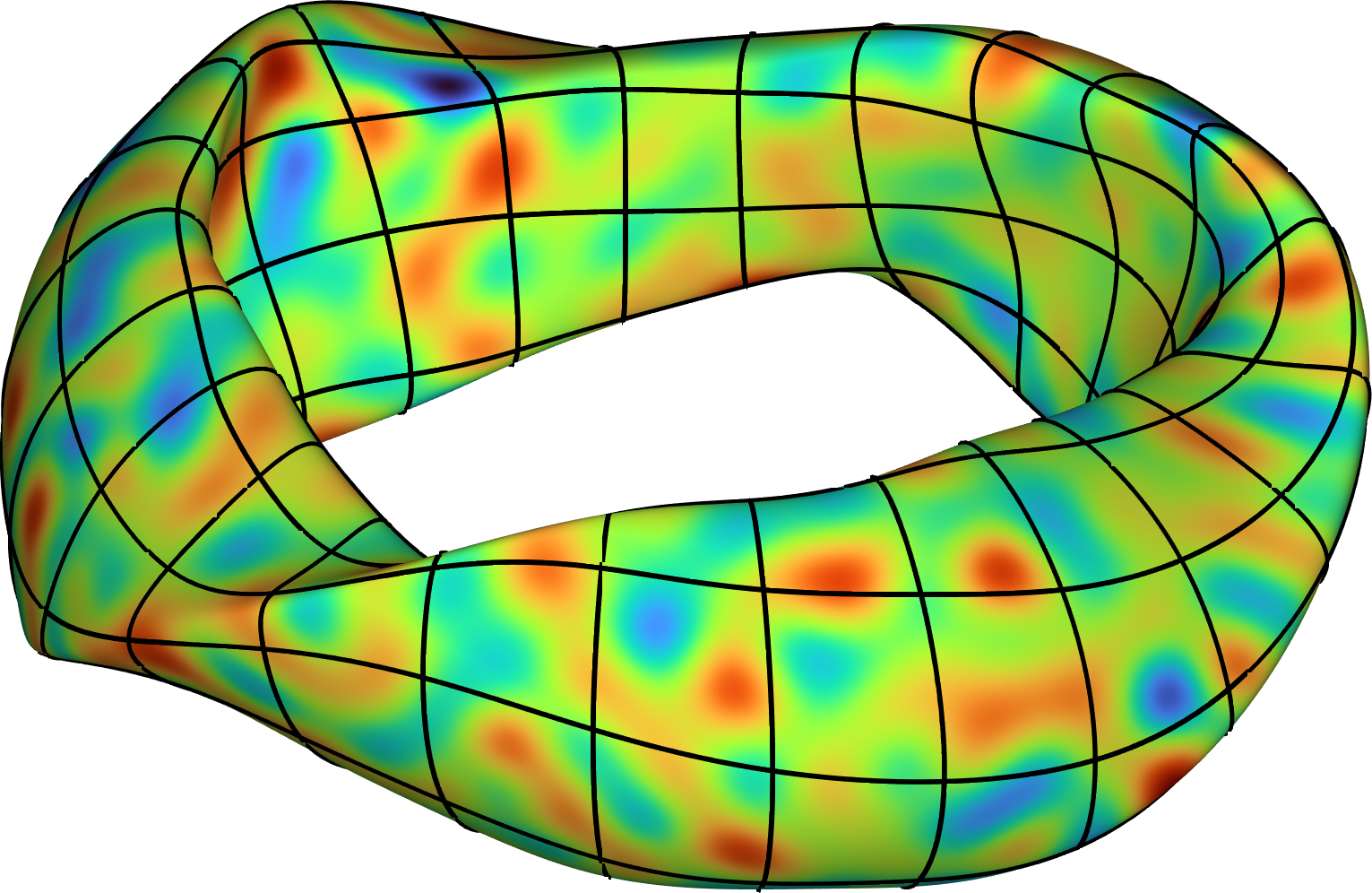} &
	\includegraphics[width=0.3\textwidth]{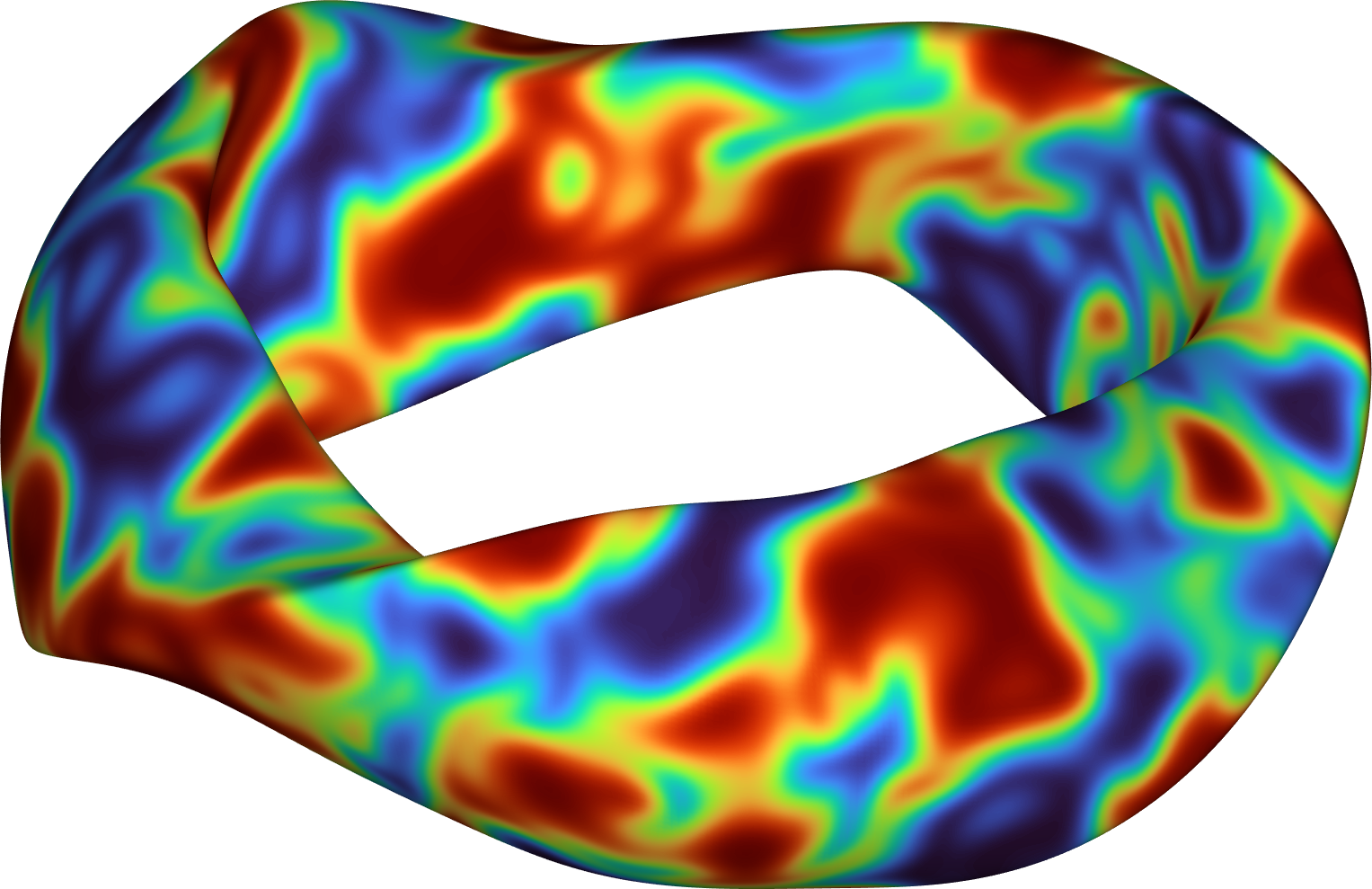} &
	\includegraphics[width=0.3\textwidth]{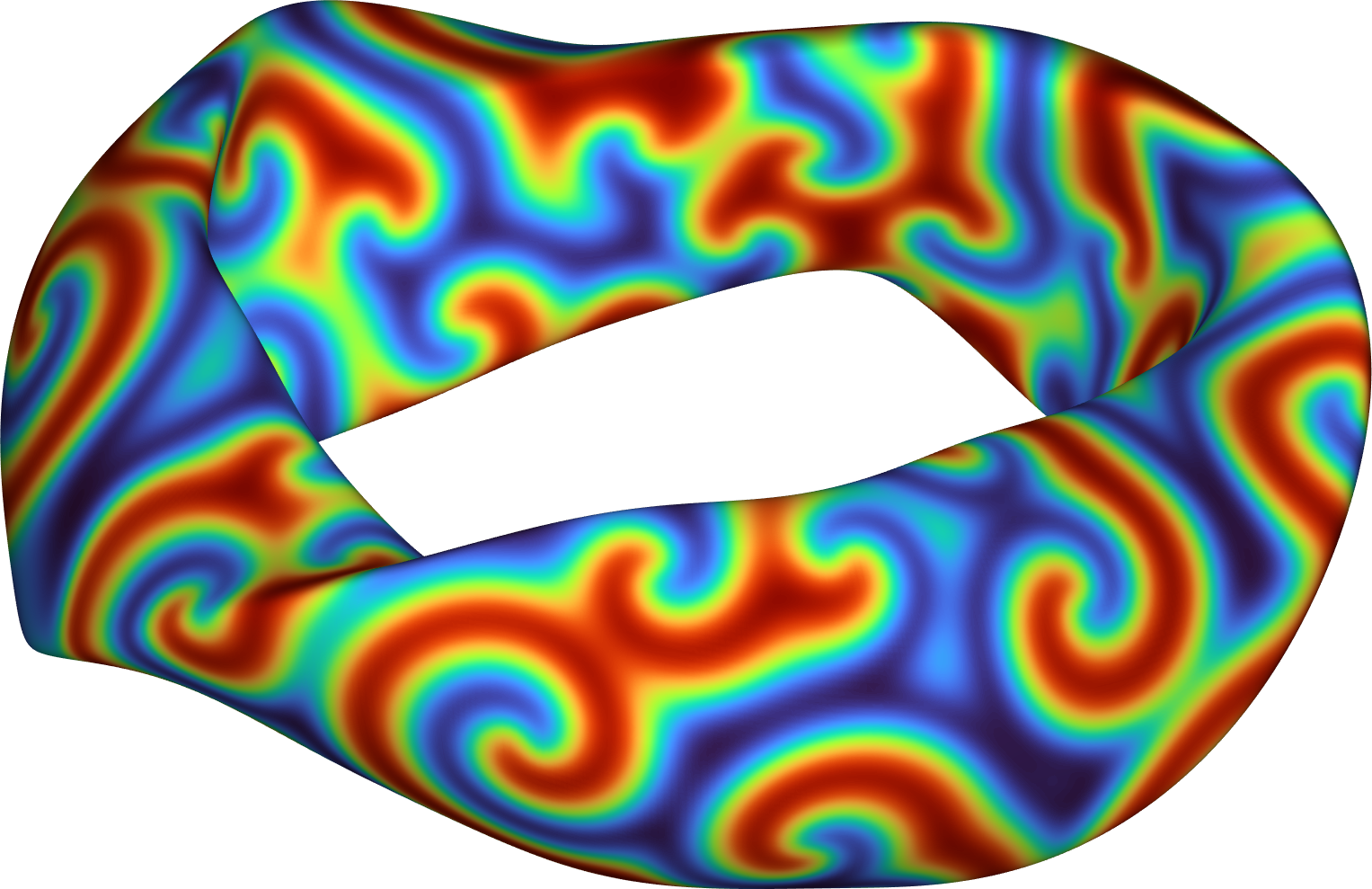} \\[1em]
	\includegraphics[width=0.29\textwidth]{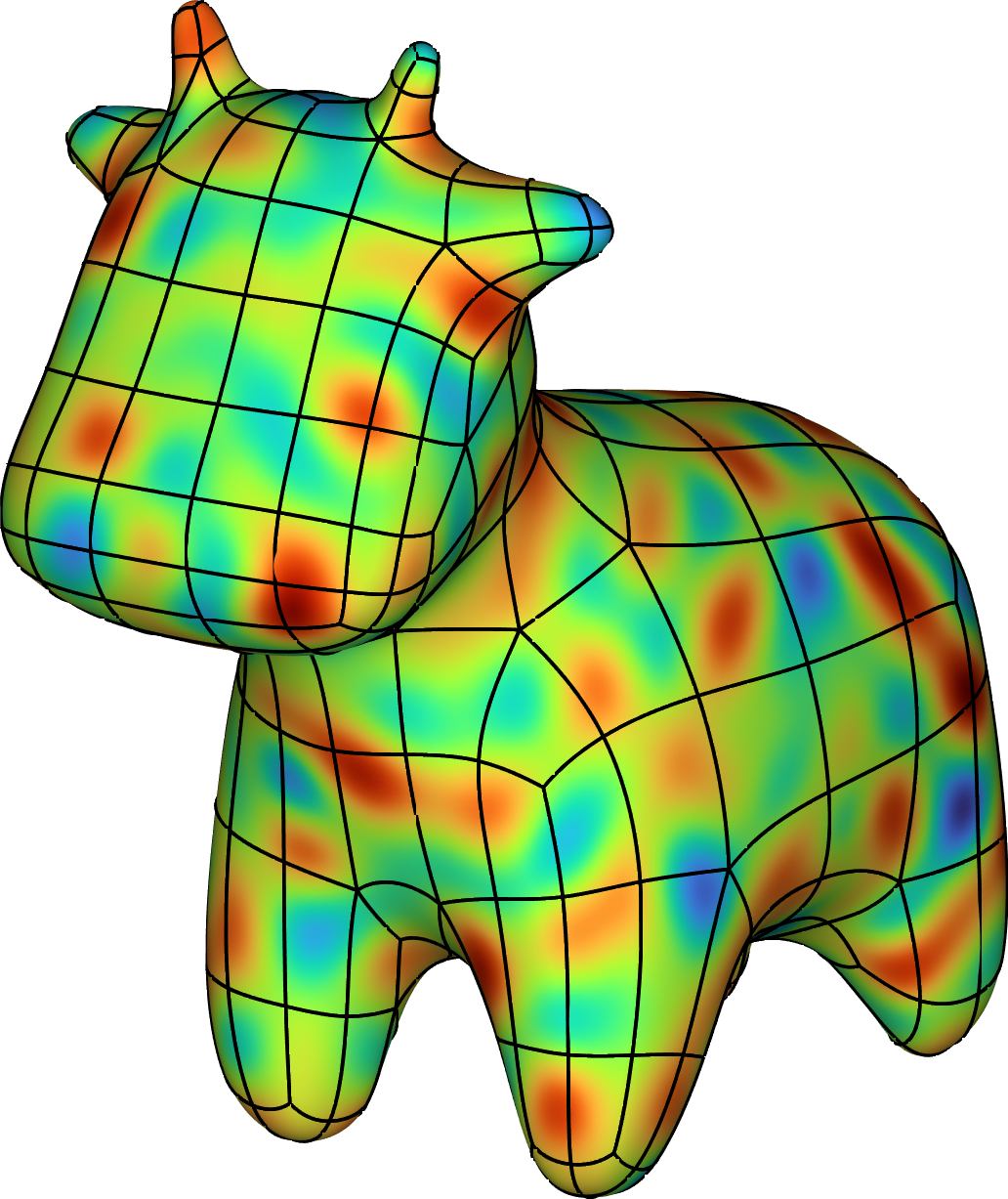} &
	\includegraphics[width=0.29\textwidth]{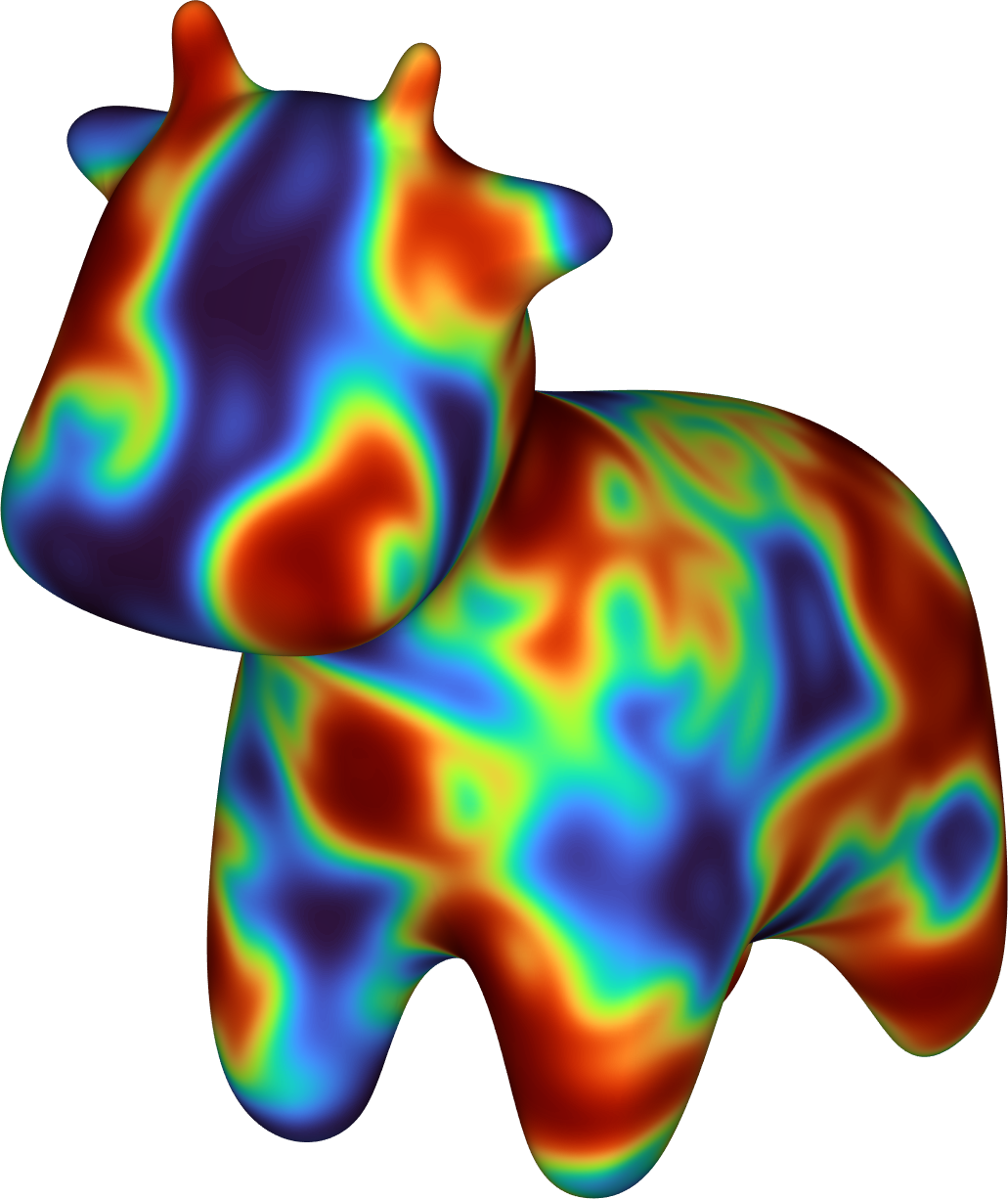} &
	\includegraphics[width=0.29\textwidth]{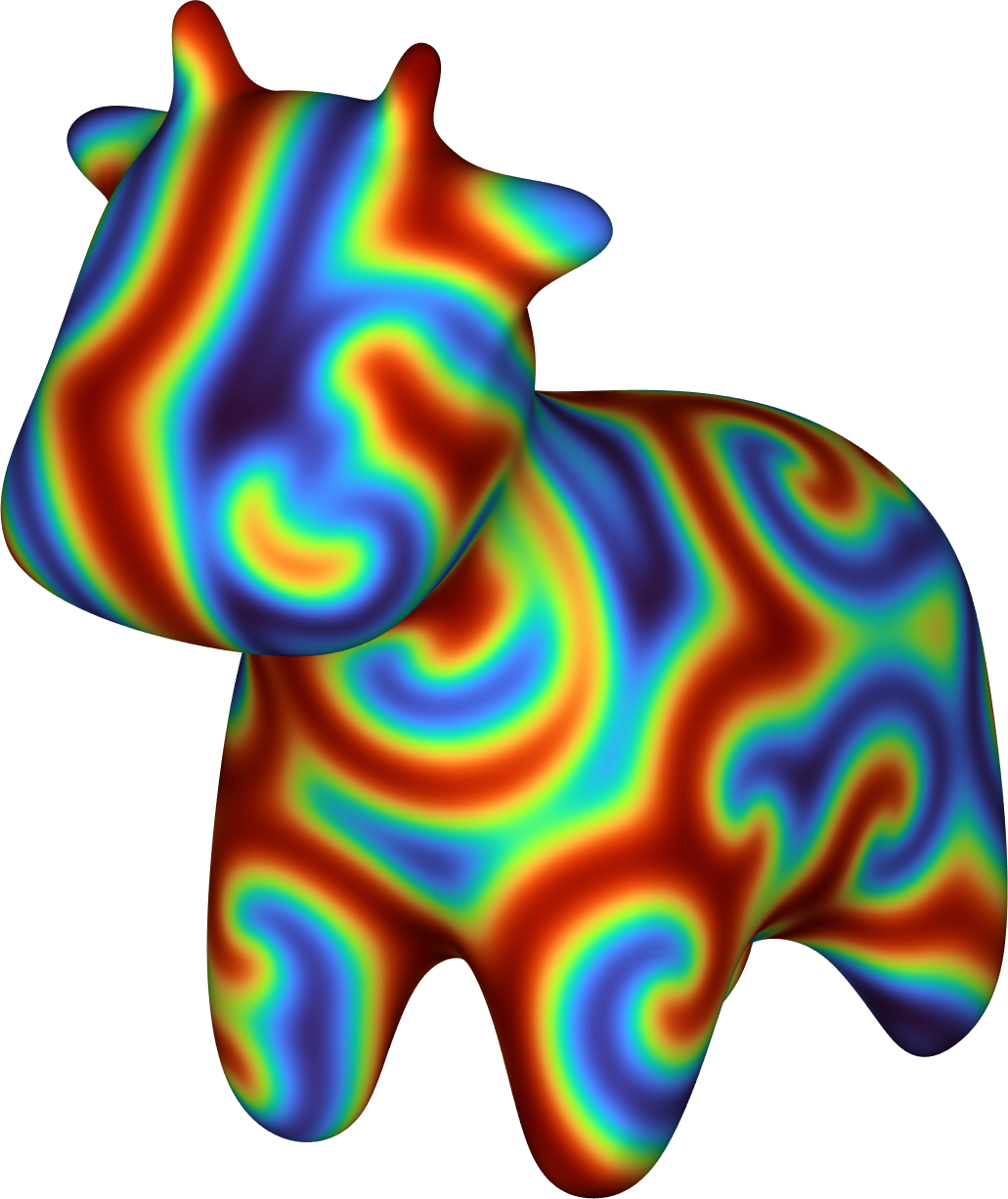}
	\end{tabular}
	\caption{Snapshots of the complex Ginzburg--Landau equation simulated on three surfaces at times $t=0$, $t=6$, and $t=60$. We discretize in time with a fourth-order IMEX-BDF4 scheme and precompute a fast direct solver for the implicit problem, which computes the solution at each time step in $\mathcal{O}(N \log N)$ operations.}
	\label{fig:cgle}
\end{figure}

The complex Ginzburg--Landau equation on a surface $\surf$ is given by
\begin{equation}\label{eq:cgle}
\frac{\partial u}{\partial t} = \delta (1+\alpha i) \lapbel u + u - (1+\beta i) u \lvert u \rvert^2,
\end{equation}
where $u = u(\vec{x}, t) \in \mathbb{C}$ and $\alpha$, $\beta$, and $\delta$ are parameters. On the blob, we take $\alpha = 0$, $\beta = -1.5$, and $\delta = 10^{-3}$. On the stellarator, we take $\alpha = 0$, $\beta = 1.5$, and $\delta = 10^{-2}$. On the cow, we take $\alpha = 0$, $\beta = 1.5$, and $\delta = 5 \cdot 10^{-4}$. We simulate this system on all three surface meshes for 2000 time steps until a final time of $t = 60$ using a time step size of $\dt = 0.03$. Snapshots of the solutions at times $t = 0$, $t = 6$, and $t = 60$ are shown in \cref{fig:cgle}.

\section{Conclusions and future work}\label{sec:conclusion}

We have presented a high-order accurate, fast direct solver for surface PDEs. The method is competitive with state-of-the-art solvers for elliptic PDEs on surfaces and can provide significant speedups when high polynomial degrees are employed. The precomputed factorizations stored by the solver can be used to accelerate implicit time-stepping schemes, allowing for long-time simulations with a low cost per time step.

Many avenues exist for extending the solver. Although we have restricted our presentation to second-order elliptic surface PDEs, the generalization to high-order differential operators (e.g., the surface biharmonic operator) is straightforward and simply requires the use of a different Poincar\'{e}--Steklov operator to impose higher-order continuity between patches. We also hope to extend the method to vector-valued PDEs, with the goal of solving the surface Stokes equations. Evolving and deformable surfaces pose an interesting challenge for the solver, as one factorization can no longer be used to accelerate all time steps. Instead, we plan to investigate the use of the fast direct solver as a preconditioner for nearby surfaces. Finally, we hope to couple this on-surface solver to a full bulk solver to accurately simulate coupled bulk--diffusion processes.

\section*{Acknowledgments}

We have greatly benefited from conversations with Leslie Greengard, Alex Barnett, Manas Rachh, Dhairya Malhotra, Mike O'Neil, Charlie Epstein, and Mengjian Hua. We thank Pearson Miller and Stas Shvartsman for sparking our interest in the development of fast methods for surface PDEs. Keenan Crane created the cow model used in \cref{fig:bruss,fig:cgle}.

\bibliographystyle{siamplain}
\bibliography{references}

\end{document}

%% file: figures/conv_space.pdf_tex.tex
\begingroup
  \makeatletter
  \providecommand\color[2][]{%
    \GenericError{(gnuplot) \space\space\space\@spaces}{%
      Package color not loaded in conjunction with
      terminal option `colourtext'%
    }{See the gnuplot documentation for explanation.%
    }{Either use 'blacktext' in gnuplot or load the package
      color.sty in LaTeX.}%
    \renewcommand\color[2][]{}%
  }%
  \providecommand\includegraphics[2][]{%
    \GenericError{(gnuplot) \space\space\space\@spaces}{%
      Package graphicx or graphics not loaded%
    }{See the gnuplot documentation for explanation.%
    }{The gnuplot epslatex terminal needs graphicx.sty or graphics.sty.}%
    \renewcommand\includegraphics[2][]{}%
  }%
  \providecommand\rotatebox[2]{#2}%
  \@ifundefined{ifGPcolor}{%
    \newif\ifGPcolor
    \GPcolortrue
  }{}%
  \@ifundefined{ifGPblacktext}{%
    \newif\ifGPblacktext
    \GPblacktextfalse
  }{}%
  \let\gplgaddtomacro\g@addto@macro
  \gdef\gplbacktext{}%
  \gdef\gplfronttext{}%
  \makeatother
  \ifGPblacktext
    \def\colorrgb#1{}%
    \def\colorgray#1{}%
  \else
    \ifGPcolor
      \def\colorrgb#1{\color[rgb]{#1}}%
      \def\colorgray#1{\color[gray]{#1}}%
      \expandafter\def\csname LTw\endcsname{\color{white}}%
      \expandafter\def\csname LTb\endcsname{\color{black}}%
      \expandafter\def\csname LTa\endcsname{\color{black}}%
      \expandafter\def\csname LT0\endcsname{\color[rgb]{1,0,0}}%
      \expandafter\def\csname LT1\endcsname{\color[rgb]{0,1,0}}%
      \expandafter\def\csname LT2\endcsname{\color[rgb]{0,0,1}}%
      \expandafter\def\csname LT3\endcsname{\color[rgb]{1,0,1}}%
      \expandafter\def\csname LT4\endcsname{\color[rgb]{0,1,1}}%
      \expandafter\def\csname LT5\endcsname{\color[rgb]{1,1,0}}%
      \expandafter\def\csname LT6\endcsname{\color[rgb]{0,0,0}}%
      \expandafter\def\csname LT7\endcsname{\color[rgb]{1,0.3,0}}%
      \expandafter\def\csname LT8\endcsname{\color[rgb]{0.5,0.5,0.5}}%
    \else
      \def\colorrgb#1{\color{black}}%
      \def\colorgray#1{\color[gray]{#1}}%
      \expandafter\def\csname LTw\endcsname{\color{white}}%
      \expandafter\def\csname LTb\endcsname{\color{black}}%
      \expandafter\def\csname LTa\endcsname{\color{black}}%
      \expandafter\def\csname LT0\endcsname{\color{black}}%
      \expandafter\def\csname LT1\endcsname{\color{black}}%
      \expandafter\def\csname LT2\endcsname{\color{black}}%
      \expandafter\def\csname LT3\endcsname{\color{black}}%
      \expandafter\def\csname LT4\endcsname{\color{black}}%
      \expandafter\def\csname LT5\endcsname{\color{black}}%
      \expandafter\def\csname LT6\endcsname{\color{black}}%
      \expandafter\def\csname LT7\endcsname{\color{black}}%
      \expandafter\def\csname LT8\endcsname{\color{black}}%
    \fi
  \fi
    \setlength{\unitlength}{0.0500bp}%
    \ifx\gptboxheight\undefined%
      \newlength{\gptboxheight}%
      \newlength{\gptboxwidth}%
      \newsavebox{\gptboxtext}%
    \fi%
    \setlength{\fboxrule}{0.5pt}%
    \setlength{\fboxsep}{1pt}%
    \definecolor{tbcol}{rgb}{1,1,1}%
\begin{picture}(4534.00,3174.00)%
    \gplgaddtomacro\gplbacktext{%
      \csname LTb\endcsname
      \put(508,704){\makebox(0,0)[r]{\strut{}\footnotesize $10^{-12}$}}%
      \csname LTb\endcsname
      \put(508,1079){\makebox(0,0)[r]{\strut{}\footnotesize $10^{-10}$}}%
      \csname LTb\endcsname
      \put(508,1454){\makebox(0,0)[r]{\strut{}\footnotesize $10^{-8}$}}%
      \csname LTb\endcsname
      \put(508,1829){\makebox(0,0)[r]{\strut{}\footnotesize $10^{-6}$}}%
      \csname LTb\endcsname
      \put(508,2203){\makebox(0,0)[r]{\strut{}\footnotesize $10^{-4}$}}%
      \csname LTb\endcsname
      \put(508,2578){\makebox(0,0)[r]{\strut{}\footnotesize $10^{-2}$}}%
      \csname LTb\endcsname
      \put(508,2953){\makebox(0,0)[r]{\strut{}\footnotesize $10^{0}$}}%
      \csname LTb\endcsname
      \put(588,528){\makebox(0,0){\strut{}\footnotesize $2^{1}$}}%
      \csname LTb\endcsname
      \put(1088,528){\makebox(0,0){\strut{}\footnotesize $2^{2}$}}%
      \csname LTb\endcsname
      \put(1588,528){\makebox(0,0){\strut{}\footnotesize $2^{3}$}}%
      \csname LTb\endcsname
      \put(2087,528){\makebox(0,0){\strut{}\footnotesize $2^{4}$}}%
      \csname LTb\endcsname
      \put(2587,528){\makebox(0,0){\strut{}\footnotesize $2^{5}$}}%
      \csname LTb\endcsname
      \put(3087,528){\makebox(0,0){\strut{}\footnotesize $2^{6}$}}%
    }%
    \gplgaddtomacro\gplfronttext{%
      \csname LTb\endcsname
      \put(-149,1828){\rotatebox{-270}{\makebox(0,0){\strut{}\small $L^\infty$ relative error}}}%
      \put(1837,154){\makebox(0,0){\strut{}$1/h$}}%
      \csname LTb\endcsname
      \put(3834,2843){\makebox(0,0)[l]{\strut{}\footnotesize \hspace{-0.15cm} $p=4$}}%
      \csname LTb\endcsname
      \put(3834,2623){\makebox(0,0)[l]{\strut{}\footnotesize \hspace{-0.15cm} $p=8$}}%
      \csname LTb\endcsname
      \put(3834,2403){\makebox(0,0)[l]{\strut{}\footnotesize \hspace{-0.15cm} $p=12$}}%
      \csname LTb\endcsname
      \put(3834,2183){\makebox(0,0)[l]{\strut{}\footnotesize \hspace{-0.15cm} $p=16$}}%
      \csname LTb\endcsname
      \put(3834,1963){\makebox(0,0)[l]{\strut{}\footnotesize \hspace{-0.15cm} $p=20$}}%
      \csname LTb\endcsname
      \put(3834,1743){\makebox(0,0)[l]{\strut{}\footnotesize \hspace{-0.15cm} $\mathcal{O}(h^{p-1})$}}%
    }%
    \gplbacktext
    \put(0,0){\includegraphics[width={226.70bp},height={158.70bp}]{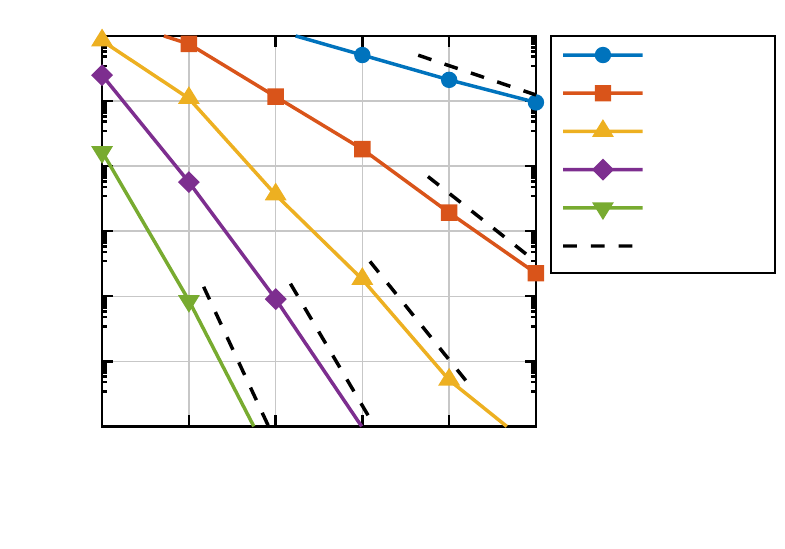}}%
    \gplfronttext
  \end{picture}%
\endgroup

%% file: figures/conv_cube.pdf_tex.tex
\begingroup
  \makeatletter
  \providecommand\color[2][]{%
    \GenericError{(gnuplot) \space\space\space\@spaces}{%
      Package color not loaded in conjunction with
      terminal option `colourtext'%
    }{See the gnuplot documentation for explanation.%
    }{Either use 'blacktext' in gnuplot or load the package
      color.sty in LaTeX.}%
    \renewcommand\color[2][]{}%
  }%
  \providecommand\includegraphics[2][]{%
    \GenericError{(gnuplot) \space\space\space\@spaces}{%
      Package graphicx or graphics not loaded%
    }{See the gnuplot documentation for explanation.%
    }{The gnuplot epslatex terminal needs graphicx.sty or graphics.sty.}%
    \renewcommand\includegraphics[2][]{}%
  }%
  \providecommand\rotatebox[2]{#2}%
  \@ifundefined{ifGPcolor}{%
    \newif\ifGPcolor
    \GPcolortrue
  }{}%
  \@ifundefined{ifGPblacktext}{%
    \newif\ifGPblacktext
    \GPblacktextfalse
  }{}%
  \let\gplgaddtomacro\g@addto@macro
  \gdef\gplbacktext{}%
  \gdef\gplfronttext{}%
  \makeatother
  \ifGPblacktext
    \def\colorrgb#1{}%
    \def\colorgray#1{}%
  \else
    \ifGPcolor
      \def\colorrgb#1{\color[rgb]{#1}}%
      \def\colorgray#1{\color[gray]{#1}}%
      \expandafter\def\csname LTw\endcsname{\color{white}}%
      \expandafter\def\csname LTb\endcsname{\color{black}}%
      \expandafter\def\csname LTa\endcsname{\color{black}}%
      \expandafter\def\csname LT0\endcsname{\color[rgb]{1,0,0}}%
      \expandafter\def\csname LT1\endcsname{\color[rgb]{0,1,0}}%
      \expandafter\def\csname LT2\endcsname{\color[rgb]{0,0,1}}%
      \expandafter\def\csname LT3\endcsname{\color[rgb]{1,0,1}}%
      \expandafter\def\csname LT4\endcsname{\color[rgb]{0,1,1}}%
      \expandafter\def\csname LT5\endcsname{\color[rgb]{1,1,0}}%
      \expandafter\def\csname LT6\endcsname{\color[rgb]{0,0,0}}%
      \expandafter\def\csname LT7\endcsname{\color[rgb]{1,0.3,0}}%
      \expandafter\def\csname LT8\endcsname{\color[rgb]{0.5,0.5,0.5}}%
    \else
      \def\colorrgb#1{\color{black}}%
      \def\colorgray#1{\color[gray]{#1}}%
      \expandafter\def\csname LTw\endcsname{\color{white}}%
      \expandafter\def\csname LTb\endcsname{\color{black}}%
      \expandafter\def\csname LTa\endcsname{\color{black}}%
      \expandafter\def\csname LT0\endcsname{\color{black}}%
      \expandafter\def\csname LT1\endcsname{\color{black}}%
      \expandafter\def\csname LT2\endcsname{\color{black}}%
      \expandafter\def\csname LT3\endcsname{\color{black}}%
      \expandafter\def\csname LT4\endcsname{\color{black}}%
      \expandafter\def\csname LT5\endcsname{\color{black}}%
      \expandafter\def\csname LT6\endcsname{\color{black}}%
      \expandafter\def\csname LT7\endcsname{\color{black}}%
      \expandafter\def\csname LT8\endcsname{\color{black}}%
    \fi
  \fi
    \setlength{\unitlength}{0.0500bp}%
    \ifx\gptboxheight\undefined%
      \newlength{\gptboxheight}%
      \newlength{\gptboxwidth}%
      \newsavebox{\gptboxtext}%
    \fi%
    \setlength{\fboxrule}{0.5pt}%
    \setlength{\fboxsep}{1pt}%
    \definecolor{tbcol}{rgb}{1,1,1}%
\begin{picture}(4534.00,3174.00)%
    \gplgaddtomacro\gplbacktext{%
      \csname LTb\endcsname
      \put(508,704){\makebox(0,0)[r]{\strut{}\footnotesize $10^{-3}$}}%
      \csname LTb\endcsname
      \put(508,1454){\makebox(0,0)[r]{\strut{}\footnotesize $10^{-2}$}}%
      \csname LTb\endcsname
      \put(508,2203){\makebox(0,0)[r]{\strut{}\footnotesize $10^{-1}$}}%
      \csname LTb\endcsname
      \put(508,2953){\makebox(0,0)[r]{\strut{}\footnotesize $10^{0}$}}%
      \csname LTb\endcsname
      \put(588,528){\makebox(0,0){\strut{}\footnotesize $2^{0}$}}%
      \csname LTb\endcsname
      \put(1088,528){\makebox(0,0){\strut{}\footnotesize $2^{1}$}}%
      \csname LTb\endcsname
      \put(1588,528){\makebox(0,0){\strut{}\footnotesize $2^{2}$}}%
      \csname LTb\endcsname
      \put(2087,528){\makebox(0,0){\strut{}\footnotesize $2^{3}$}}%
      \csname LTb\endcsname
      \put(2587,528){\makebox(0,0){\strut{}\footnotesize $2^{4}$}}%
      \csname LTb\endcsname
      \put(3087,528){\makebox(0,0){\strut{}\footnotesize $2^{5}$}}%
    }%
    \gplgaddtomacro\gplfronttext{%
      \csname LTb\endcsname
      \put(-83,1828){\rotatebox{-270}{\makebox(0,0){\strut{}\small $L^\infty$ relative error}}}%
      \put(1837,154){\makebox(0,0){\strut{}$1/h$}}%
      \csname LTb\endcsname
      \put(3834,2183){\makebox(0,0)[l]{\strut{}\footnotesize \hspace{-0.15cm} $p=6$}}%
      \csname LTb\endcsname
      \put(3834,2403){\makebox(0,0)[l]{\strut{}\footnotesize \hspace{-0.15cm} $p=5$}}%
      \csname LTb\endcsname
      \put(3834,2623){\makebox(0,0)[l]{\strut{}\footnotesize \hspace{-0.15cm} $p=4$}}%
      \csname LTb\endcsname
      \put(3834,2843){\makebox(0,0)[l]{\strut{}\footnotesize \hspace{-0.15cm} $p=3$}}%
    }%
    \gplbacktext
    \put(0,0){\includegraphics[width={226.70bp},height={158.70bp}]{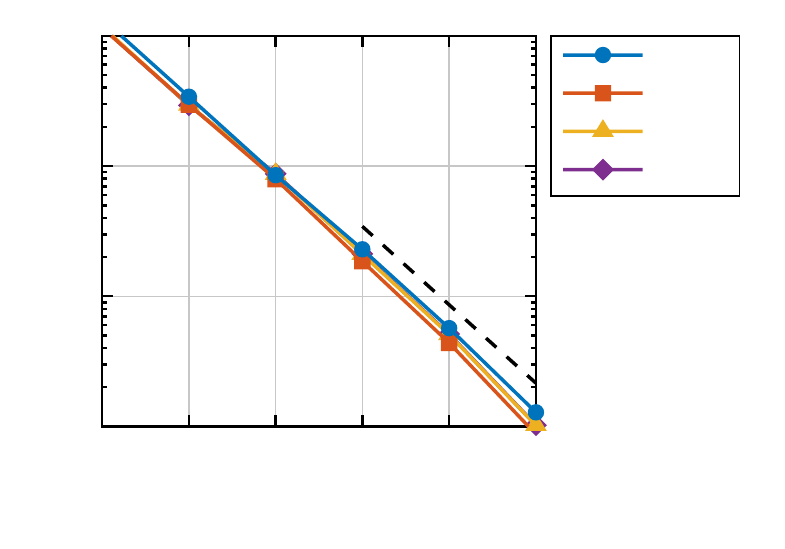}}%
    \gplfronttext
  \end{picture}%
\endgroup

%% file: figures/conv_twisted.pdf_tex.tex
\begingroup
  \makeatletter
  \providecommand\color[2][]{%
    \GenericError{(gnuplot) \space\space\space\@spaces}{%
      Package color not loaded in conjunction with
      terminal option `colourtext'%
    }{See the gnuplot documentation for explanation.%
    }{Either use 'blacktext' in gnuplot or load the package
      color.sty in LaTeX.}%
    \renewcommand\color[2][]{}%
  }%
  \providecommand\includegraphics[2][]{%
    \GenericError{(gnuplot) \space\space\space\@spaces}{%
      Package graphicx or graphics not loaded%
    }{See the gnuplot documentation for explanation.%
    }{The gnuplot epslatex terminal needs graphicx.sty or graphics.sty.}%
    \renewcommand\includegraphics[2][]{}%
  }%
  \providecommand\rotatebox[2]{#2}%
  \@ifundefined{ifGPcolor}{%
    \newif\ifGPcolor
    \GPcolortrue
  }{}%
  \@ifundefined{ifGPblacktext}{%
    \newif\ifGPblacktext
    \GPblacktextfalse
  }{}%
  \let\gplgaddtomacro\g@addto@macro
  \gdef\gplbacktext{}%
  \gdef\gplfronttext{}%
  \makeatother
  \ifGPblacktext
    \def\colorrgb#1{}%
    \def\colorgray#1{}%
  \else
    \ifGPcolor
      \def\colorrgb#1{\color[rgb]{#1}}%
      \def\colorgray#1{\color[gray]{#1}}%
      \expandafter\def\csname LTw\endcsname{\color{white}}%
      \expandafter\def\csname LTb\endcsname{\color{black}}%
      \expandafter\def\csname LTa\endcsname{\color{black}}%
      \expandafter\def\csname LT0\endcsname{\color[rgb]{1,0,0}}%
      \expandafter\def\csname LT1\endcsname{\color[rgb]{0,1,0}}%
      \expandafter\def\csname LT2\endcsname{\color[rgb]{0,0,1}}%
      \expandafter\def\csname LT3\endcsname{\color[rgb]{1,0,1}}%
      \expandafter\def\csname LT4\endcsname{\color[rgb]{0,1,1}}%
      \expandafter\def\csname LT5\endcsname{\color[rgb]{1,1,0}}%
      \expandafter\def\csname LT6\endcsname{\color[rgb]{0,0,0}}%
      \expandafter\def\csname LT7\endcsname{\color[rgb]{1,0.3,0}}%
      \expandafter\def\csname LT8\endcsname{\color[rgb]{0.5,0.5,0.5}}%
    \else
      \def\colorrgb#1{\color{black}}%
      \def\colorgray#1{\color[gray]{#1}}%
      \expandafter\def\csname LTw\endcsname{\color{white}}%
      \expandafter\def\csname LTb\endcsname{\color{black}}%
      \expandafter\def\csname LTa\endcsname{\color{black}}%
      \expandafter\def\csname LT0\endcsname{\color{black}}%
      \expandafter\def\csname LT1\endcsname{\color{black}}%
      \expandafter\def\csname LT2\endcsname{\color{black}}%
      \expandafter\def\csname LT3\endcsname{\color{black}}%
      \expandafter\def\csname LT4\endcsname{\color{black}}%
      \expandafter\def\csname LT5\endcsname{\color{black}}%
      \expandafter\def\csname LT6\endcsname{\color{black}}%
      \expandafter\def\csname LT7\endcsname{\color{black}}%
      \expandafter\def\csname LT8\endcsname{\color{black}}%
    \fi
  \fi
    \setlength{\unitlength}{0.0500bp}%
    \ifx\gptboxheight\undefined%
      \newlength{\gptboxheight}%
      \newlength{\gptboxwidth}%
      \newsavebox{\gptboxtext}%
    \fi%
    \setlength{\fboxrule}{0.5pt}%
    \setlength{\fboxsep}{1pt}%
    \definecolor{tbcol}{rgb}{1,1,1}%
\begin{picture}(4534.00,3174.00)%
    \gplgaddtomacro\gplbacktext{%
      \csname LTb\endcsname
      \put(508,704){\makebox(0,0)[r]{\strut{}\footnotesize $10^{-12}$}}%
      \csname LTb\endcsname
      \put(508,1079){\makebox(0,0)[r]{\strut{}\footnotesize $10^{-10}$}}%
      \csname LTb\endcsname
      \put(508,1454){\makebox(0,0)[r]{\strut{}\footnotesize $10^{-8}$}}%
      \csname LTb\endcsname
      \put(508,1829){\makebox(0,0)[r]{\strut{}\footnotesize $10^{-6}$}}%
      \csname LTb\endcsname
      \put(508,2203){\makebox(0,0)[r]{\strut{}\footnotesize $10^{-4}$}}%
      \csname LTb\endcsname
      \put(508,2578){\makebox(0,0)[r]{\strut{}\footnotesize $10^{-2}$}}%
      \csname LTb\endcsname
      \put(508,2953){\makebox(0,0)[r]{\strut{}\footnotesize $10^{0}$}}%
      \csname LTb\endcsname
      \put(588,528){\makebox(0,0){\strut{}\footnotesize $2^{2}$}}%
      \csname LTb\endcsname
      \put(1421,528){\makebox(0,0){\strut{}\footnotesize $2^{3}$}}%
      \csname LTb\endcsname
      \put(2254,528){\makebox(0,0){\strut{}\footnotesize $2^{4}$}}%
      \csname LTb\endcsname
      \put(3087,528){\makebox(0,0){\strut{}\footnotesize $2^{5}$}}%
    }%
    \gplgaddtomacro\gplfronttext{%
      \csname LTb\endcsname
      \put(-215,1828){\rotatebox{-270}{\makebox(0,0){\strut{}\small $L^\infty$ relative error}}}%
      \put(1837,154){\makebox(0,0){\strut{}$1/h$}}%
      \csname LTb\endcsname
      \put(3834,2843){\makebox(0,0)[l]{\strut{}\footnotesize \hspace{-0.15cm} $p=4$}}%
      \csname LTb\endcsname
      \put(3834,2623){\makebox(0,0)[l]{\strut{}\footnotesize \hspace{-0.15cm} $p=8$}}%
      \csname LTb\endcsname
      \put(3834,2403){\makebox(0,0)[l]{\strut{}\footnotesize \hspace{-0.15cm} $p=12$}}%
      \csname LTb\endcsname
      \put(3834,2183){\makebox(0,0)[l]{\strut{}\footnotesize \hspace{-0.15cm} $\mathcal{O}(h^{p-1})$}}%
    }%
    \gplbacktext
    \put(0,0){\includegraphics[width={226.70bp},height={158.70bp}]{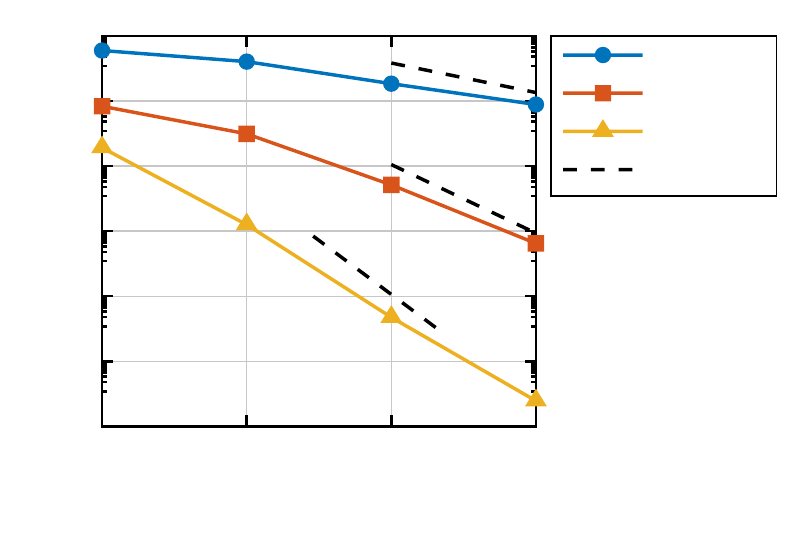}}%
    \gplfronttext
  \end{picture}%
\endgroup

%% file: figures/complexity.pdf_tex.tex
\begingroup
  \makeatletter
  \providecommand\color[2][]{%
    \GenericError{(gnuplot) \space\space\space\@spaces}{%
      Package color not loaded in conjunction with
      terminal option `colourtext'%
    }{See the gnuplot documentation for explanation.%
    }{Either use 'blacktext' in gnuplot or load the package
      color.sty in LaTeX.}%
    \renewcommand\color[2][]{}%
  }%
  \providecommand\includegraphics[2][]{%
    \GenericError{(gnuplot) \space\space\space\@spaces}{%
      Package graphicx or graphics not loaded%
    }{See the gnuplot documentation for explanation.%
    }{The gnuplot epslatex terminal needs graphicx.sty or graphics.sty.}%
    \renewcommand\includegraphics[2][]{}%
  }%
  \providecommand\rotatebox[2]{#2}%
  \@ifundefined{ifGPcolor}{%
    \newif\ifGPcolor
    \GPcolortrue
  }{}%
  \@ifundefined{ifGPblacktext}{%
    \newif\ifGPblacktext
    \GPblacktextfalse
  }{}%
  \let\gplgaddtomacro\g@addto@macro
  \gdef\gplbacktext{}%
  \gdef\gplfronttext{}%
  \makeatother
  \ifGPblacktext
    \def\colorrgb#1{}%
    \def\colorgray#1{}%
  \else
    \ifGPcolor
      \def\colorrgb#1{\color[rgb]{#1}}%
      \def\colorgray#1{\color[gray]{#1}}%
      \expandafter\def\csname LTw\endcsname{\color{white}}%
      \expandafter\def\csname LTb\endcsname{\color{black}}%
      \expandafter\def\csname LTa\endcsname{\color{black}}%
      \expandafter\def\csname LT0\endcsname{\color[rgb]{1,0,0}}%
      \expandafter\def\csname LT1\endcsname{\color[rgb]{0,1,0}}%
      \expandafter\def\csname LT2\endcsname{\color[rgb]{0,0,1}}%
      \expandafter\def\csname LT3\endcsname{\color[rgb]{1,0,1}}%
      \expandafter\def\csname LT4\endcsname{\color[rgb]{0,1,1}}%
      \expandafter\def\csname LT5\endcsname{\color[rgb]{1,1,0}}%
      \expandafter\def\csname LT6\endcsname{\color[rgb]{0,0,0}}%
      \expandafter\def\csname LT7\endcsname{\color[rgb]{1,0.3,0}}%
      \expandafter\def\csname LT8\endcsname{\color[rgb]{0.5,0.5,0.5}}%
    \else
      \def\colorrgb#1{\color{black}}%
      \def\colorgray#1{\color[gray]{#1}}%
      \expandafter\def\csname LTw\endcsname{\color{white}}%
      \expandafter\def\csname LTb\endcsname{\color{black}}%
      \expandafter\def\csname LTa\endcsname{\color{black}}%
      \expandafter\def\csname LT0\endcsname{\color{black}}%
      \expandafter\def\csname LT1\endcsname{\color{black}}%
      \expandafter\def\csname LT2\endcsname{\color{black}}%
      \expandafter\def\csname LT3\endcsname{\color{black}}%
      \expandafter\def\csname LT4\endcsname{\color{black}}%
      \expandafter\def\csname LT5\endcsname{\color{black}}%
      \expandafter\def\csname LT6\endcsname{\color{black}}%
      \expandafter\def\csname LT7\endcsname{\color{black}}%
      \expandafter\def\csname LT8\endcsname{\color{black}}%
    \fi
  \fi
    \setlength{\unitlength}{0.0500bp}%
    \ifx\gptboxheight\undefined%
      \newlength{\gptboxheight}%
      \newlength{\gptboxwidth}%
      \newsavebox{\gptboxtext}%
    \fi%
    \setlength{\fboxrule}{0.5pt}%
    \setlength{\fboxsep}{1pt}%
    \definecolor{tbcol}{rgb}{1,1,1}%
\begin{picture}(6518.00,3400.00)%
    \gplgaddtomacro\gplbacktext{%
      \csname LTb\endcsname
      \put(311,782){\makebox(0,0)[r]{\strut{}\footnotesize $10^{-3}$}}%
      \csname LTb\endcsname
      \put(311,1190){\makebox(0,0)[r]{\strut{}\footnotesize $10^{-2}$}}%
      \csname LTb\endcsname
      \put(311,1598){\makebox(0,0)[r]{\strut{}\footnotesize $10^{-1}$}}%
      \csname LTb\endcsname
      \put(311,2005){\makebox(0,0)[r]{\strut{}\footnotesize $10^{0}$}}%
      \csname LTb\endcsname
      \put(311,2413){\makebox(0,0)[r]{\strut{}\footnotesize $10^{1}$}}%
      \csname LTb\endcsname
      \put(311,2821){\makebox(0,0)[r]{\strut{}\footnotesize $10^{2}$}}%
      \csname LTb\endcsname
      \put(391,606){\makebox(0,0){\strut{}\footnotesize $2^{2}$}}%
      \csname LTb\endcsname
      \put(803,606){\makebox(0,0){\strut{}\footnotesize $2^{4}$}}%
      \csname LTb\endcsname
      \put(1215,606){\makebox(0,0){\strut{}\footnotesize $2^{6}$}}%
      \csname LTb\endcsname
      \put(1626,606){\makebox(0,0){\strut{}\footnotesize $2^{8}$}}%
      \csname LTb\endcsname
      \put(2038,606){\makebox(0,0){\strut{}\footnotesize $2^{10}$}}%
      \csname LTb\endcsname
      \put(2450,606){\makebox(0,0){\strut{}\footnotesize $2^{12}$}}%
    }%
    \gplgaddtomacro\gplfronttext{%
      \csname LTb\endcsname
      \put(-280,1801){\rotatebox{-270}{\makebox(0,0){\strut{}Runtime (s)}}}%
      \put(1523,298){\makebox(0,0){\strut{}$N$}}%
      \csname LTb\endcsname
      \put(1714,3221){\makebox(0,0)[l]{\strut{}Factorization}}%
    }%
    \gplgaddtomacro\gplbacktext{%
    }%
    \gplgaddtomacro\gplfronttext{%
      \csname LTb\endcsname
      \put(3519,3221){\makebox(0,0)[l]{\strut{}Solve}}%
    }%
    \gplgaddtomacro\gplbacktext{%
    }%
    \gplgaddtomacro\gplfronttext{%
      \csname LTb\endcsname
      \put(4691,3221){\makebox(0,0)[l]{\strut{}Update}}%
    }%
    \gplgaddtomacro\gplbacktext{%
      \csname LTb\endcsname
      \put(3830,782){\makebox(0,0)[r]{\strut{}\footnotesize $10^{-3}$}}%
      \csname LTb\endcsname
      \put(3830,1190){\makebox(0,0)[r]{\strut{}\footnotesize $10^{-2}$}}%
      \csname LTb\endcsname
      \put(3830,1598){\makebox(0,0)[r]{\strut{}\footnotesize $10^{-1}$}}%
      \csname LTb\endcsname
      \put(3830,2005){\makebox(0,0)[r]{\strut{}\footnotesize $10^{0}$}}%
      \csname LTb\endcsname
      \put(3830,2413){\makebox(0,0)[r]{\strut{}\footnotesize $10^{1}$}}%
      \csname LTb\endcsname
      \put(3830,2821){\makebox(0,0)[r]{\strut{}\footnotesize $10^{2}$}}%
      \csname LTb\endcsname
      \put(3910,606){\makebox(0,0){\strut{}\footnotesize $2^{2}$}}%
      \csname LTb\endcsname
      \put(4322,606){\makebox(0,0){\strut{}\footnotesize $2^{4}$}}%
      \csname LTb\endcsname
      \put(4734,606){\makebox(0,0){\strut{}\footnotesize $2^{6}$}}%
      \csname LTb\endcsname
      \put(5145,606){\makebox(0,0){\strut{}\footnotesize $2^{8}$}}%
      \csname LTb\endcsname
      \put(5557,606){\makebox(0,0){\strut{}\footnotesize $2^{10}$}}%
      \csname LTb\endcsname
      \put(5969,606){\makebox(0,0){\strut{}\footnotesize $2^{12}$}}%
    }%
    \gplgaddtomacro\gplfronttext{%
      \csname LTb\endcsname
      \put(3265,1801){\rotatebox{-270}{\makebox(0,0){\strut{}Memory (GB)}}}%
      \put(5042,298){\makebox(0,0){\strut{}$N$}}%
    }%
    \gplbacktext
    \put(0,0){\includegraphics[width={325.90bp},height={170.00bp}]{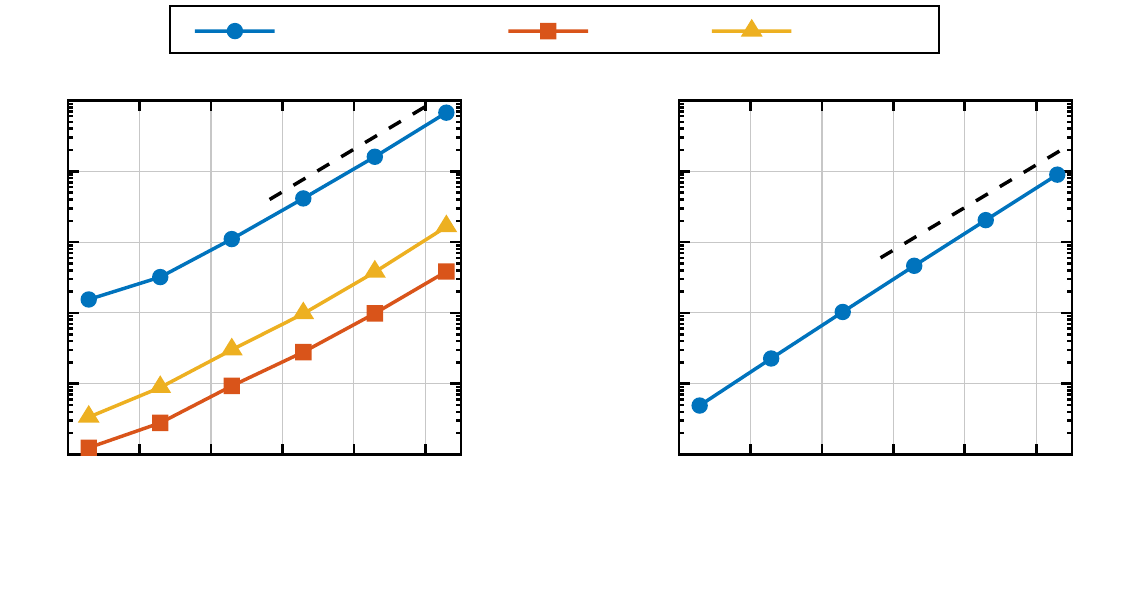}}%
    \gplfronttext
  \end{picture}%
\endgroup

%% file: figures/stellarator_benchmark.pdf_tex.tex
\begingroup
  \makeatletter
  \providecommand\color[2][]{%
    \GenericError{(gnuplot) \space\space\space\@spaces}{%
      Package color not loaded in conjunction with
      terminal option `colourtext'%
    }{See the gnuplot documentation for explanation.%
    }{Either use 'blacktext' in gnuplot or load the package
      color.sty in LaTeX.}%
    \renewcommand\color[2][]{}%
  }%
  \providecommand\includegraphics[2][]{%
    \GenericError{(gnuplot) \space\space\space\@spaces}{%
      Package graphicx or graphics not loaded%
    }{See the gnuplot documentation for explanation.%
    }{The gnuplot epslatex terminal needs graphicx.sty or graphics.sty.}%
    \renewcommand\includegraphics[2][]{}%
  }%
  \providecommand\rotatebox[2]{#2}%
  \@ifundefined{ifGPcolor}{%
    \newif\ifGPcolor
    \GPcolortrue
  }{}%
  \@ifundefined{ifGPblacktext}{%
    \newif\ifGPblacktext
    \GPblacktextfalse
  }{}%
  \let\gplgaddtomacro\g@addto@macro
  \gdef\gplbacktext{}%
  \gdef\gplfronttext{}%
  \makeatother
  \ifGPblacktext
    \def\colorrgb#1{}%
    \def\colorgray#1{}%
  \else
    \ifGPcolor
      \def\colorrgb#1{\color[rgb]{#1}}%
      \def\colorgray#1{\color[gray]{#1}}%
      \expandafter\def\csname LTw\endcsname{\color{white}}%
      \expandafter\def\csname LTb\endcsname{\color{black}}%
      \expandafter\def\csname LTa\endcsname{\color{black}}%
      \expandafter\def\csname LT0\endcsname{\color[rgb]{1,0,0}}%
      \expandafter\def\csname LT1\endcsname{\color[rgb]{0,1,0}}%
      \expandafter\def\csname LT2\endcsname{\color[rgb]{0,0,1}}%
      \expandafter\def\csname LT3\endcsname{\color[rgb]{1,0,1}}%
      \expandafter\def\csname LT4\endcsname{\color[rgb]{0,1,1}}%
      \expandafter\def\csname LT5\endcsname{\color[rgb]{1,1,0}}%
      \expandafter\def\csname LT6\endcsname{\color[rgb]{0,0,0}}%
      \expandafter\def\csname LT7\endcsname{\color[rgb]{1,0.3,0}}%
      \expandafter\def\csname LT8\endcsname{\color[rgb]{0.5,0.5,0.5}}%
    \else
      \def\colorrgb#1{\color{black}}%
      \def\colorgray#1{\color[gray]{#1}}%
      \expandafter\def\csname LTw\endcsname{\color{white}}%
      \expandafter\def\csname LTb\endcsname{\color{black}}%
      \expandafter\def\csname LTa\endcsname{\color{black}}%
      \expandafter\def\csname LT0\endcsname{\color{black}}%
      \expandafter\def\csname LT1\endcsname{\color{black}}%
      \expandafter\def\csname LT2\endcsname{\color{black}}%
      \expandafter\def\csname LT3\endcsname{\color{black}}%
      \expandafter\def\csname LT4\endcsname{\color{black}}%
      \expandafter\def\csname LT5\endcsname{\color{black}}%
      \expandafter\def\csname LT6\endcsname{\color{black}}%
      \expandafter\def\csname LT7\endcsname{\color{black}}%
      \expandafter\def\csname LT8\endcsname{\color{black}}%
    \fi
  \fi
    \setlength{\unitlength}{0.0500bp}%
    \ifx\gptboxheight\undefined%
      \newlength{\gptboxheight}%
      \newlength{\gptboxwidth}%
      \newsavebox{\gptboxtext}%
    \fi%
    \setlength{\fboxrule}{0.5pt}%
    \setlength{\fboxsep}{1pt}%
    \definecolor{tbcol}{rgb}{1,1,1}%
\begin{picture}(5498.00,3798.00)%
    \gplgaddtomacro\gplbacktext{%
      \csname LTb\endcsname
      \put(1671,660){\makebox(0,0)[r]{\strut{}\footnotesize $10^{-12}$}}%
      \csname LTb\endcsname
      \put(1671,1058){\makebox(0,0)[r]{\strut{}\footnotesize $10^{-10}$}}%
      \csname LTb\endcsname
      \put(1671,1456){\makebox(0,0)[r]{\strut{}\footnotesize $10^{-8}$}}%
      \csname LTb\endcsname
      \put(1671,1855){\makebox(0,0)[r]{\strut{}\footnotesize $10^{-6}$}}%
      \csname LTb\endcsname
      \put(1671,2253){\makebox(0,0)[r]{\strut{}\footnotesize $10^{-4}$}}%
      \csname LTb\endcsname
      \put(1671,2651){\makebox(0,0)[r]{\strut{}\footnotesize $10^{-2}$}}%
      \csname LTb\endcsname
      \put(1671,3049){\makebox(0,0)[r]{\strut{}\footnotesize $10^{0}$}}%
      \csname LTb\endcsname
      \put(1751,484){\makebox(0,0){\strut{}\footnotesize $10^{-2}$}}%
      \csname LTb\endcsname
      \put(2402,484){\makebox(0,0){\strut{}\footnotesize $10^{-1}$}}%
      \csname LTb\endcsname
      \put(3053,484){\makebox(0,0){\strut{}\footnotesize $10^{0}$}}%
      \csname LTb\endcsname
      \put(3704,484){\makebox(0,0){\strut{}\footnotesize $10^{1}$}}%
      \csname LTb\endcsname
      \put(4354,484){\makebox(0,0){\strut{}\footnotesize $10^{2}$}}%
    }%
    \gplgaddtomacro\gplfronttext{%
      \csname LTb\endcsname
      \put(1014,1854){\rotatebox{-270}{\makebox(0,0){\strut{}\small $L^\infty$ relative error}}}%
      \put(3078,154){\makebox(0,0){\strut{}\small Runtime (s)}}%
      \csname LTb\endcsname
      \put(1541,3630){\makebox(0,0)[l]{\strut{}\footnotesize Surface HPS}}%
      \csname LTb\endcsname
      \put(1541,3410){\makebox(0,0)[l]{\strut{}\footnotesize Pseudospectral}}%
      \csname LTb\endcsname
      \put(3386,3630){\makebox(0,0)[l]{\strut{}\footnotesize MFEM \texttt{+} UMFPACK}}%
      \csname LTb\endcsname
      \put(3386,3410){\makebox(0,0)[l]{\strut{}\footnotesize Layer potential}}%
    }%
    \gplbacktext
    \put(0,0){\includegraphics[width={274.90bp},height={189.90bp}]{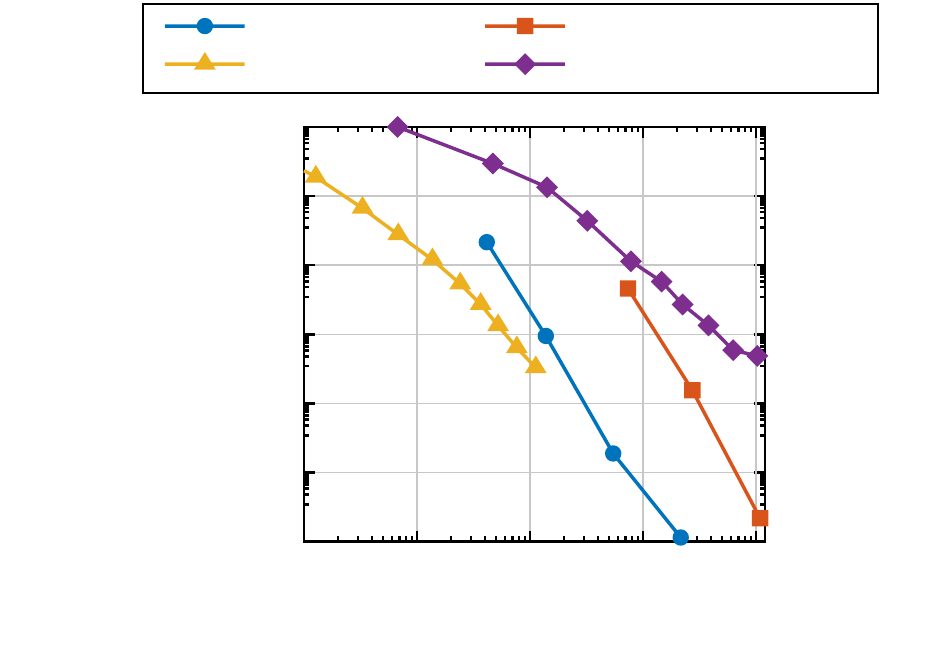}}%
    \gplfronttext
  \end{picture}%
\endgroup

%% file: figures/mfem_comparison.pdf_tex.tex
\begingroup
  \makeatletter
  \providecommand\color[2][]{%
    \GenericError{(gnuplot) \space\space\space\@spaces}{%
      Package color not loaded in conjunction with
      terminal option `colourtext'%
    }{See the gnuplot documentation for explanation.%
    }{Either use 'blacktext' in gnuplot or load the package
      color.sty in LaTeX.}%
    \renewcommand\color[2][]{}%
  }%
  \providecommand\includegraphics[2][]{%
    \GenericError{(gnuplot) \space\space\space\@spaces}{%
      Package graphicx or graphics not loaded%
    }{See the gnuplot documentation for explanation.%
    }{The gnuplot epslatex terminal needs graphicx.sty or graphics.sty.}%
    \renewcommand\includegraphics[2][]{}%
  }%
  \providecommand\rotatebox[2]{#2}%
  \@ifundefined{ifGPcolor}{%
    \newif\ifGPcolor
    \GPcolortrue
  }{}%
  \@ifundefined{ifGPblacktext}{%
    \newif\ifGPblacktext
    \GPblacktextfalse
  }{}%
  \let\gplgaddtomacro\g@addto@macro
  \gdef\gplbacktext{}%
  \gdef\gplfronttext{}%
  \makeatother
  \ifGPblacktext
    \def\colorrgb#1{}%
    \def\colorgray#1{}%
  \else
    \ifGPcolor
      \def\colorrgb#1{\color[rgb]{#1}}%
      \def\colorgray#1{\color[gray]{#1}}%
      \expandafter\def\csname LTw\endcsname{\color{white}}%
      \expandafter\def\csname LTb\endcsname{\color{black}}%
      \expandafter\def\csname LTa\endcsname{\color{black}}%
      \expandafter\def\csname LT0\endcsname{\color[rgb]{1,0,0}}%
      \expandafter\def\csname LT1\endcsname{\color[rgb]{0,1,0}}%
      \expandafter\def\csname LT2\endcsname{\color[rgb]{0,0,1}}%
      \expandafter\def\csname LT3\endcsname{\color[rgb]{1,0,1}}%
      \expandafter\def\csname LT4\endcsname{\color[rgb]{0,1,1}}%
      \expandafter\def\csname LT5\endcsname{\color[rgb]{1,1,0}}%
      \expandafter\def\csname LT6\endcsname{\color[rgb]{0,0,0}}%
      \expandafter\def\csname LT7\endcsname{\color[rgb]{1,0.3,0}}%
      \expandafter\def\csname LT8\endcsname{\color[rgb]{0.5,0.5,0.5}}%
    \else
      \def\colorrgb#1{\color{black}}%
      \def\colorgray#1{\color[gray]{#1}}%
      \expandafter\def\csname LTw\endcsname{\color{white}}%
      \expandafter\def\csname LTb\endcsname{\color{black}}%
      \expandafter\def\csname LTa\endcsname{\color{black}}%
      \expandafter\def\csname LT0\endcsname{\color{black}}%
      \expandafter\def\csname LT1\endcsname{\color{black}}%
      \expandafter\def\csname LT2\endcsname{\color{black}}%
      \expandafter\def\csname LT3\endcsname{\color{black}}%
      \expandafter\def\csname LT4\endcsname{\color{black}}%
      \expandafter\def\csname LT5\endcsname{\color{black}}%
      \expandafter\def\csname LT6\endcsname{\color{black}}%
      \expandafter\def\csname LT7\endcsname{\color{black}}%
      \expandafter\def\csname LT8\endcsname{\color{black}}%
    \fi
  \fi
    \setlength{\unitlength}{0.0500bp}%
    \ifx\gptboxheight\undefined%
      \newlength{\gptboxheight}%
      \newlength{\gptboxwidth}%
      \newsavebox{\gptboxtext}%
    \fi%
    \setlength{\fboxrule}{0.5pt}%
    \setlength{\fboxsep}{1pt}%
    \definecolor{tbcol}{rgb}{1,1,1}%
\begin{picture}(5668.00,3174.00)%
    \gplgaddtomacro\gplbacktext{%
      \csname LTb\endcsname
      \put(260,709){\makebox(0,0)[r]{\strut{}\footnotesize $10^{-1}$}}%
      \csname LTb\endcsname
      \put(260,1206){\makebox(0,0)[r]{\strut{}\footnotesize $10^{0}$}}%
      \csname LTb\endcsname
      \put(260,1703){\makebox(0,0)[r]{\strut{}\footnotesize $10^{1}$}}%
      \csname LTb\endcsname
      \put(260,2200){\makebox(0,0)[r]{\strut{}\footnotesize $10^{2}$}}%
      \csname LTb\endcsname
      \put(260,2697){\makebox(0,0)[r]{\strut{}\footnotesize $10^{3}$}}%
      \csname LTb\endcsname
      \put(340,533){\makebox(0,0){\strut{}$2^5$}}%
      \csname LTb\endcsname
      \put(1077,533){\makebox(0,0){\strut{}$2^7$}}%
      \csname LTb\endcsname
      \put(1813,533){\makebox(0,0){\strut{}$2^9$}}%
      \csname LTb\endcsname
      \put(2550,533){\makebox(0,0){\strut{}$2^{11}$}}%
    }%
    \gplgaddtomacro\gplfronttext{%
      \csname LTb\endcsname
      \put(-331,1703){\rotatebox{-270}{\makebox(0,0){\strut{}Runtime (s)}}}%
      \put(1445,225){\makebox(0,0){\strut{}$N$}}%
      \csname LTb\endcsname
      \put(1415,3015){\makebox(0,0)[l]{\strut{}Surface HPS}}%
      \csname LTb\endcsname
      \put(3194,3015){\makebox(0,0)[l]{\strut{}MFEM \texttt{+} UMFPACK}}%
    }%
    \gplgaddtomacro\gplbacktext{%
      \csname LTb\endcsname
      \put(3094,709){\makebox(0,0)[r]{\strut{}}}%
      \csname LTb\endcsname
      \put(3094,1206){\makebox(0,0)[r]{\strut{}}}%
      \csname LTb\endcsname
      \put(3094,1703){\makebox(0,0)[r]{\strut{}}}%
      \csname LTb\endcsname
      \put(3094,2200){\makebox(0,0)[r]{\strut{}}}%
      \csname LTb\endcsname
      \put(3094,2697){\makebox(0,0)[r]{\strut{}}}%
      \csname LTb\endcsname
      \put(3174,533){\makebox(0,0){\strut{}\footnotesize 4}}%
      \csname LTb\endcsname
      \put(3726,533){\makebox(0,0){\strut{}\footnotesize 8}}%
      \csname LTb\endcsname
      \put(4279,533){\makebox(0,0){\strut{}\footnotesize 12}}%
      \csname LTb\endcsname
      \put(4831,533){\makebox(0,0){\strut{}\footnotesize 16}}%
      \csname LTb\endcsname
      \put(5383,533){\makebox(0,0){\strut{}\footnotesize 20}}%
    }%
    \gplgaddtomacro\gplfronttext{%
      \csname LTb\endcsname
      \put(4278,225){\makebox(0,0){\strut{}$p$}}%
    }%
    \gplbacktext
    \put(0,0){\includegraphics[width={283.40bp},height={158.70bp}]{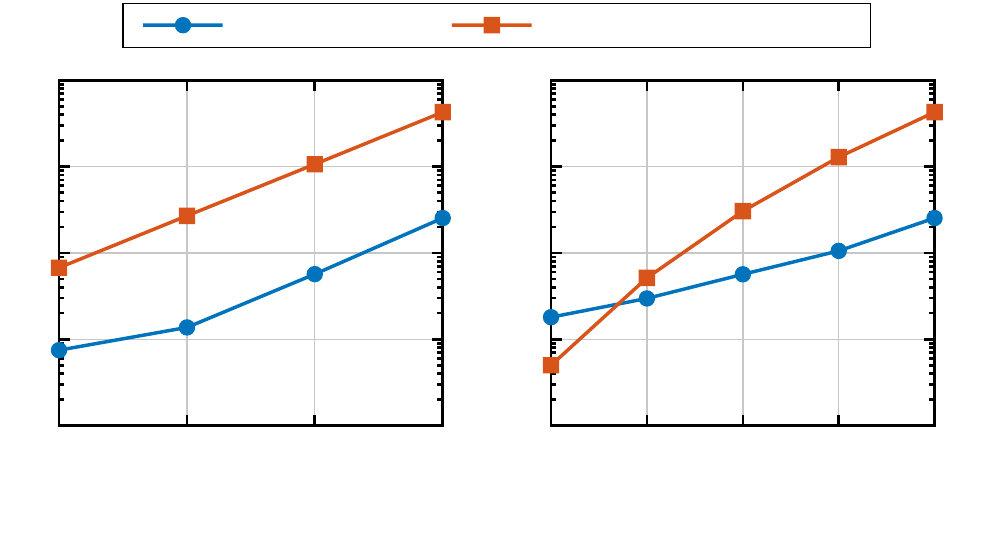}}%
    \gplfronttext
  \end{picture}%
\endgroup

%% file: figures/conv_time.pdf_tex.tex
\begingroup
  \makeatletter
  \providecommand\color[2][]{%
    \GenericError{(gnuplot) \space\space\space\@spaces}{%
      Package color not loaded in conjunction with
      terminal option `colourtext'%
    }{See the gnuplot documentation for explanation.%
    }{Either use 'blacktext' in gnuplot or load the package
      color.sty in LaTeX.}%
    \renewcommand\color[2][]{}%
  }%
  \providecommand\includegraphics[2][]{%
    \GenericError{(gnuplot) \space\space\space\@spaces}{%
      Package graphicx or graphics not loaded%
    }{See the gnuplot documentation for explanation.%
    }{The gnuplot epslatex terminal needs graphicx.sty or graphics.sty.}%
    \renewcommand\includegraphics[2][]{}%
  }%
  \providecommand\rotatebox[2]{#2}%
  \@ifundefined{ifGPcolor}{%
    \newif\ifGPcolor
    \GPcolortrue
  }{}%
  \@ifundefined{ifGPblacktext}{%
    \newif\ifGPblacktext
    \GPblacktextfalse
  }{}%
  \let\gplgaddtomacro\g@addto@macro
  \gdef\gplbacktext{}%
  \gdef\gplfronttext{}%
  \makeatother
  \ifGPblacktext
    \def\colorrgb#1{}%
    \def\colorgray#1{}%
  \else
    \ifGPcolor
      \def\colorrgb#1{\color[rgb]{#1}}%
      \def\colorgray#1{\color[gray]{#1}}%
      \expandafter\def\csname LTw\endcsname{\color{white}}%
      \expandafter\def\csname LTb\endcsname{\color{black}}%
      \expandafter\def\csname LTa\endcsname{\color{black}}%
      \expandafter\def\csname LT0\endcsname{\color[rgb]{1,0,0}}%
      \expandafter\def\csname LT1\endcsname{\color[rgb]{0,1,0}}%
      \expandafter\def\csname LT2\endcsname{\color[rgb]{0,0,1}}%
      \expandafter\def\csname LT3\endcsname{\color[rgb]{1,0,1}}%
      \expandafter\def\csname LT4\endcsname{\color[rgb]{0,1,1}}%
      \expandafter\def\csname LT5\endcsname{\color[rgb]{1,1,0}}%
      \expandafter\def\csname LT6\endcsname{\color[rgb]{0,0,0}}%
      \expandafter\def\csname LT7\endcsname{\color[rgb]{1,0.3,0}}%
      \expandafter\def\csname LT8\endcsname{\color[rgb]{0.5,0.5,0.5}}%
    \else
      \def\colorrgb#1{\color{black}}%
      \def\colorgray#1{\color[gray]{#1}}%
      \expandafter\def\csname LTw\endcsname{\color{white}}%
      \expandafter\def\csname LTb\endcsname{\color{black}}%
      \expandafter\def\csname LTa\endcsname{\color{black}}%
      \expandafter\def\csname LT0\endcsname{\color{black}}%
      \expandafter\def\csname LT1\endcsname{\color{black}}%
      \expandafter\def\csname LT2\endcsname{\color{black}}%
      \expandafter\def\csname LT3\endcsname{\color{black}}%
      \expandafter\def\csname LT4\endcsname{\color{black}}%
      \expandafter\def\csname LT5\endcsname{\color{black}}%
      \expandafter\def\csname LT6\endcsname{\color{black}}%
      \expandafter\def\csname LT7\endcsname{\color{black}}%
      \expandafter\def\csname LT8\endcsname{\color{black}}%
    \fi
  \fi
    \setlength{\unitlength}{0.0500bp}%
    \ifx\gptboxheight\undefined%
      \newlength{\gptboxheight}%
      \newlength{\gptboxwidth}%
      \newsavebox{\gptboxtext}%
    \fi%
    \setlength{\fboxrule}{0.5pt}%
    \setlength{\fboxsep}{1pt}%
    \definecolor{tbcol}{rgb}{1,1,1}%
\begin{picture}(6236.00,3968.00)%
    \gplgaddtomacro\gplbacktext{%
      \csname LTb\endcsname
      \put(1262,950){\makebox(0,0)[r]{\strut{}\footnotesize $10^{-8}$}}%
      \csname LTb\endcsname
      \put(1262,1291){\makebox(0,0)[r]{\strut{}\footnotesize $10^{-7}$}}%
      \csname LTb\endcsname
      \put(1262,1633){\makebox(0,0)[r]{\strut{}\footnotesize $10^{-6}$}}%
      \csname LTb\endcsname
      \put(1262,1974){\makebox(0,0)[r]{\strut{}\footnotesize $10^{-5}$}}%
      \csname LTb\endcsname
      \put(1262,2315){\makebox(0,0)[r]{\strut{}\footnotesize $10^{-4}$}}%
      \csname LTb\endcsname
      \put(1262,2656){\makebox(0,0)[r]{\strut{}\footnotesize $10^{-3}$}}%
      \csname LTb\endcsname
      \put(1262,2997){\makebox(0,0)[r]{\strut{}\footnotesize $10^{-2}$}}%
      \csname LTb\endcsname
      \put(1262,3338){\makebox(0,0)[r]{\strut{}\footnotesize $10^{-1}$}}%
      \csname LTb\endcsname
      \put(1262,3679){\makebox(0,0)[r]{\strut{}\footnotesize $10^{0}$}}%
      \csname LTb\endcsname
      \put(1342,596){\makebox(0,0){\strut{}\footnotesize $2^{3}$}}%
      \csname LTb\endcsname
      \put(1880,596){\makebox(0,0){\strut{}\footnotesize $2^{4}$}}%
      \csname LTb\endcsname
      \put(2418,596){\makebox(0,0){\strut{}\footnotesize $2^{5}$}}%
      \csname LTb\endcsname
      \put(2957,596){\makebox(0,0){\strut{}\footnotesize $2^{6}$}}%
      \csname LTb\endcsname
      \put(3495,596){\makebox(0,0){\strut{}\footnotesize $2^{7}$}}%
      \csname LTb\endcsname
      \put(4033,596){\makebox(0,0){\strut{}\footnotesize $2^{8}$}}%
      \csname LTb\endcsname
      \put(4571,596){\makebox(0,0){\strut{}\footnotesize $2^{9}$}}%
    }%
    \gplgaddtomacro\gplfronttext{%
      \csname LTb\endcsname
      \put(605,2225){\rotatebox{-270}{\makebox(0,0){\strut{}$L^\infty$ relative error}}}%
      \put(2956,222){\makebox(0,0){\strut{}$1/\Delta t$}}%
      \csname LTb\endcsname
      \put(5311,3569){\makebox(0,0)[l]{\strut{}\footnotesize \hspace{-0.15cm} BDF1}}%
      \csname LTb\endcsname
      \put(5311,3349){\makebox(0,0)[l]{\strut{}\footnotesize \hspace{-0.15cm} BDF2}}%
      \csname LTb\endcsname
      \put(5311,3129){\makebox(0,0)[l]{\strut{}\footnotesize \hspace{-0.15cm} BDF3}}%
      \csname LTb\endcsname
      \put(5311,2909){\makebox(0,0)[l]{\strut{}\footnotesize \hspace{-0.15cm} BDF4}}%
    }%
    \gplbacktext
    \put(0,0){\includegraphics[width={311.80bp},height={198.40bp}]{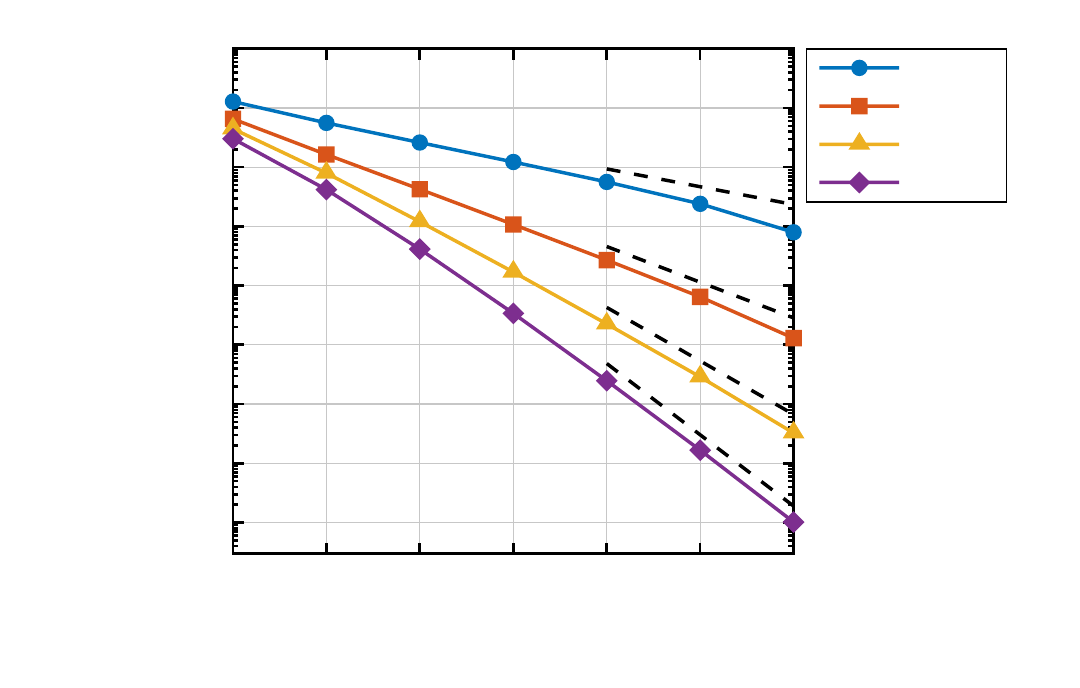}}%
    \gplfronttext
  \end{picture}%
\endgroup